\documentclass[10pt,reqno]{amsart}
\usepackage[T1]{fontenc}
\usepackage{lmodern,microtype,amsmath,amssymb,amsthm,mathtools}
\usepackage[margin=28mm]{geometry}
\usepackage{array,enumerate,xcolor,hyperref}

\hypersetup{hidelinks,pdftitle={Nonlinear stability of planar contact waves for the three-dimensional compressible Navier--Stokes equations}}
\allowdisplaybreaks[2]
\numberwithin{equation}{section}

\newtheorem{theorem}{Theorem}[section]
\newtheorem{proposition}[theorem]{Proposition}
\newtheorem{lemma}[theorem]{Lemma}
\newtheorem{corollary}[theorem]{Corollary}
\theoremstyle{remark}
\newtheorem{remark}[theorem]{Remark}

\newcommand{\R}{\mathbb{R}}
\newcommand{\eps}{\varepsilon}
\DeclareMathOperator{\diver}{div}
\DeclareMathOperator{\curl}{curl}

\begin{document}
	
	\title[Stability of planar contact wave for 3-d Navier-Stokes equations]{Nonlinear stability of planar contact waves
		for the 3-D compressible Navier--Stokes equations in the whole space}
	\author[F.~Huang]{Feimin Huang}
	\address[F.~Huang]{Academy of Mathematics and Systems Science, Chinese Academy of Sciences, Beijing 100190, China, and School of Mathematical Sciences, University of Chinese Academy of Sciences, Beijing 100049, China.}
	\email{fhuang@amt.ac.cn}
	
	\author[R.~Li]{Rui Li}
	\address[R.~Li]{Department of Mathematics, The Chinese University of Hong Kong, Shatin, Hong Kong, China.}
	\email{ruili001@cuhk.edu.hk}

	\author[X.~Wu]{Xiaochun Wu}
	\address[X.~Wu]{School of Mathematics and Statistics,
		Central South University, Changsha, China.}
	\email{xcwu22@csu.edu.cn}
	
	\author[L.~Xu]{Lingda Xu}
	\address[L.~Xu]{Department of Mathematics, The Chinese University of Hong Kong, Hong Kong, China.}
	\email{lingdaxu@cuhk.edu.hk}

	\begin{abstract}
		We prove the global nonlinear stability of small planar contact waves for the three-dimensional full compressible Navier--Stokes equations in $\R^3$ under general $H^4(\R^3)$ perturbations, without any $L^1$ integrability or zero-mass assumptions. Unlike the one-dimensional problem or the periodic channel, the problem in $\R^3$ admits neither a Poincaré inequality in the transverse directions nor a square-integrable anti-derivative. Our proof relies on three main ingredients: a perturbation around the one-dimensional Navier--Stokes solution starting from the smooth contact profile, a formulation of the low-order perturbation system in the velocity, pressure, and entropy variables with separate estimates for the velocity--pressure pair and the entropy, and the construction of characteristic and space-time Gaussian weight functions. Combining these weighted estimates with higher-order derivative bounds and homogeneous Gagliardo--Nirenberg interpolation yields global stability and algebraic $L^\infty$ decay rates in $\R^3$.
	\end{abstract}
	
	\maketitle
	
	\begingroup
	\makeatletter
	\renewcommand{\@makefntext}[1]{\noindent#1}
	\makeatother
	\footnotetext[0]{%
		\textbf{2020 Mathematics Subject Classification:}
		Primary 35Q30, 35B35; Secondary 35B40, 76N10.\par
		\smallskip
		\noindent\textbf{Keywords:}
		compressible Navier--Stokes equations; planar contact wave;
		nonlinear stability; decay rate; relative entropy method;
		space-time Gaussian weight.}
	\endgroup

	\section{Introduction}
	
	\subsection{The model and problem}
	
	We consider the three-dimensional full compressible Navier--Stokes equations for a viscous, heat-conducting ideal polytropic gas in Eulerian coordinates:
	\begin{equation}
		\left\{\begin{aligned}
			&\rho_t+\diver(\rho u)=0, \\[2pt]
			&\rho(\partial_t+u\cdot\nabla)u+\nabla p
			=\mu\Delta u+(\mu+\lambda)\nabla\diver u, \\[2pt]
			&c_v\rho(\partial_t+u\cdot\nabla)\theta+p\diver u
			=\kappa\Delta\theta+\mathcal{Q}(u),
		\end{aligned}\right.
		\label{eq:NSF}
	\end{equation}
	for $(t,x)\in[0,\infty)\times\R^3$, where $\rho>0$, $u=(u_1,u_2,u_3)^{\top}$, and $\theta>0$ denote the fluid density, velocity, and absolute temperature, respectively. The pressure $p$, the viscous dissipation function $\mathcal{Q}(u)$, and the thermodynamic and viscosity parameters are given by
	\begin{gather*}
		p=R\rho\theta,\qquad
		\mathcal{Q}(u)=\frac{\mu}{2}|\nabla u+\nabla u^{\top}|^2+\lambda(\diver u)^2,\\[2pt]
		\gamma=1+\frac{R}{c_v}>1,\qquad c_p=c_v+R,\qquad \nu=2\mu+\lambda>0,
	\end{gather*}
	where the gas constant $R$, the specific heat $c_v$, the shear viscosity $\mu$, and the thermal conductivity $\kappa$ are positive constants, and the bulk viscosity $\lambda$ satisfies $3\lambda+2\mu\ge 0$.
	
	For the inviscid Euler system, a planar contact discontinuity is a fundamental elementary wave across which the pressure and normal velocity remain constant while the density and temperature jump. In the presence of viscosity and heat conductivity, a planar contact discontinuity without tangential shear is smoothed into a diffusive contact wave connecting two constant far-field states $(\rho_\pm,0,\theta_\pm)$ under equal pressure. The aim of this paper is to establish the global nonlinear stability and algebraic $L^\infty$ decay rates of planar contact waves in $\R^3$ under general $H^4(\R^3)$ perturbations.
	
	\subsection{Related work}
	
	The mathematical theory of global classical solutions and time-decay rates for the multidimensional compressible Navier--Stokes equations near a constant equilibrium was initiated by Matsumura and Nishida \cite{MN80,MN83} and developed through spectral, energy, and Green's function methods \cite{Kawashima,Ponce,HZ}. For one-dimensional elementary waves, the anti-derivative method of Matsumura and Nishihara \cite{MN85} and Goodman \cite{Goodman} established the nonlinear stability of viscous shock profiles under zero-mass perturbations, and general perturbations were treated in \cite{Liu85,SX,Liu97,LZ09,LZ15,Zumbrun} using shock shifts and diffusion waves. For one-dimensional rarefaction waves, global asymptotic stability was obtained in \cite{IO,LW17,MN86,MN92} using the monotonicity of the rarefaction profile.
	
	Contact waves belong to the linearly degenerate characteristic field, where neither the compressive structure of shocks nor the expansive monotonicity of rarefaction waves is available, and spatial gradients of the background wave decay at the critical rate $(1+t)^{-1/2}$. For one-dimensional conservation laws with artificial viscosity, Liu and Xin \cite{LX} proved nonlinear stability under the structural conditions
	\begin{equation}
		\nabla l_2\cdot r_2=0,\qquad \nabla r_2\cdot r_2=0,
		\label{eq:LX-conditions}
	\end{equation}
	on the left and right eigenvectors of the contact field; Xin and Zeng \cite{XZ} later removed the left condition in \eqref{eq:LX-conditions} for systems with artificial viscosity. For the one-dimensional compressible Navier--Stokes equations in Lagrangian coordinates, where the left condition in \eqref{eq:LX-conditions} holds automatically, the nonlinear stability of viscous contact waves was established for free-boundary problems by Huang, Matsumura, and Shi \cite{HMS}, for zero-mass perturbations by Huang, Matsumura, and Xin \cite{HMX}, and for general perturbations by Huang, Xin, and Yang \cite{HXY}, with further extensions to composite waves \cite{HM,HLM,KVW23}, decay rates \cite{HWW,Yang20a,Yang20b,LWX}, and kinetic models \cite{HY,HWY,DYY}.
	
	In several space dimensions, stability results have been obtained for planar rarefaction waves \cite{Xin90,LW18,LWW,HXY22}, planar viscous shocks \cite{FS02,FS10,HLZ,LY,WW,Yuan}, compressible vortex sheets \cite{HXXY,YZ}, and inviscid limits of contact discontinuities \cite{KVW}. Recently, Duan, Huang, Li, and Xu \cite{DHLX} proved the asymptotic stability of planar contact waves for the three-dimensional compressible Navier--Stokes equations in a periodic channel in Eulerian coordinates: a linear transformation of the one-dimensional anti-derivative variables recovers the cancellation structure \eqref{eq:LX-conditions} on the transverse zero mode, while the Poincaré inequality on the periodic cross-section controls the nonzero transverse modes. In the whole space $\R^3$, however, the stability of planar contact waves has remained open because there is no Poincaré inequality in the transverse directions, the error of the self-similar profile is not square-integrable over $\R^2$, and a general perturbation in $H^4(\R^3)$ does not admit a square-integrable anti-derivative.
	
	\subsection{Main difficulties and key ideas}
	
	Working directly in the original Eulerian perturbation variables in $\R^3$ without anti-derivatives means that the cancellation conditions \eqref{eq:LX-conditions} used in one dimension \cite{HMX,LWX} and on the channel zero mode \cite{DHLX} are no longer available, while the contact wave provides no expansive velocity gradient to damp the zeroth-order perturbation. To address these difficulties, our proof is built on three main ideas: perturbing around the one-dimensional Navier--Stokes solution starting from the smooth contact profile, formulating the low-order perturbation system in the velocity, pressure, and entropy variables with separate estimates for the velocity--pressure pair and the entropy, and constructing weight functions adapted to the contact layer and the acoustic modes.
	
	First, because the self-similar contact profile leaves a one-dimensional residual that has infinite $L^2$ norm over the unbounded transverse plane $\R^2$, we take as our background state the one-dimensional Navier--Stokes solution $(\bar\rho,\bar u e_1,\bar\theta)$ with initial data $(\tilde\rho_0,\tilde u_0 e_1,\tilde\theta_0)$. Since the initial difference between this one-dimensional solution and the self-similar profile is zero, its one-dimensional anti-derivative has zero initial data, which places no zero-mass or $L^1$ condition on the three-dimensional perturbation. In Sections~\ref{sec:background-global} and \ref{sec:log-spacetime}, we combine the transformed one-dimensional anti-derivative system of \cite{DHLX} with a characteristic Gaussian weight matrix to derive an eighth-order time-weighted energy hierarchy and a Gaussian-weighted space-time bound for the background material derivative $d=\bar u_t+\bar u\bar u_{x_1}$.
	
	Second, we write the low-order three-dimensional perturbation equations in terms of the velocity perturbation $\psi=u-\bar u e_1$, the logarithmic pressure perturbation $q=\log(p/\bar p)$, and the entropy variable $h$, and we treat the estimates for $(u,p)$ and the entropy separately. In the entropy equation, the linear coupling with the background temperature gradient $\bar\theta_{x_1}$ decays at the non-integrable rate $(1+t)^{-1/2}$, while heat conduction alone does not bound the $L^2$ norm of $h$ on $\R^3$. We therefore introduce the weight function
	\begin{equation}
		\omega(t,x_1)=(1+t)^{-1}+\frac{y}{1+y},\qquad
		y=(1+t)^{-1}\exp\!\left(\frac{\alpha_0x_1^2}{1+t}\right),
		\label{eq:intro-omega}
	\end{equation}
	into the zeroth-order relative entropy functional. Inside the contact layer, $\omega(t,x_1)$ provides an additional factor $(1+t)^{-1}$ whose Gaussian growth in $x_1$ is absorbed by the Gaussian decay of $\tilde\theta_{x_1}$, making $\sqrt{\omega}\,\bar\theta_{x_1}$ {\color{black}square-integrable in time, uniformly in $x_1$}. Moreover, $\omega_t<0$ everywhere, so the time derivative of the weight yields a positive term controlling $|\partial_{x_1}\omega|^2/\omega$. {\color{black}In the outer region $|x_1|\ge\sqrt{(1+t)\log(1+t)/\alpha_0}$}, $\omega$ is bounded below by $1/2$, which avoids any loss of time weight along the acoustic cones and allows us to prove a one-dimensional weighted Sobolev inequality to control the coupling with $d$.
	
	Third, we keep the energy for the velocity--pressure pair $(\psi,q)$ and the higher-order derivatives unweighted in space. Because $\omega$ is included only in the zeroth-order entropy functional, the heat-conduction term in the pressure equation produces an unweighted inner product $\int_{\R^3}\nabla q\cdot\nabla h\,dx$. We cancel this term, as well as its first-order counterpart, by introducing cross functionals between the velocity and the entropy gradient via the momentum equation, and we group the background heat-conduction terms in the pressure estimate before integrating over $\R^3$ so that no $L^1(\R^3)$ assumption is needed. Combining the low-order estimates for $(\psi,q)$ and $h$ with the higher-order energy bounds and homogeneous Gagliardo--Nirenberg interpolation across spatial derivatives closes the a priori estimates and yields the global $L^\infty$ decay rates in Theorem~\ref{thm:main}.
	
	\subsection{The planar viscous contact profile and main result}
	
	Let $(\rho_\pm,0,\theta_\pm)$ be constant far-field states with $\rho_\pm,\theta_\pm>0$ satisfying the equal-pressure relation
	\[
	p_*=R\rho_-\theta_-=R\rho_+\theta_+, \qquad
	\delta=|\rho_+-\rho_-|+|\theta_+-\theta_-|>0.
	\]
	Following \cite{HMX,DHLX}, let $\tilde\rho_0(x_1)$ be the smooth, positive, monotone solution to the boundary value problem
	\begin{equation}
		\left\{\begin{aligned}
			&\frac{\kappa}{c_p}\left(\frac{\tilde\rho_0'}{\tilde\rho_0}\right)'+\frac{x_1}{2}\tilde\rho_0'=0, \\[2pt]
			&\tilde\rho_0(-\infty)=\rho_-,\qquad \tilde\rho_0(+\infty)=\rho_+,
		\end{aligned}\right.
		\label{eq:initial-density}
	\end{equation}
	and define the self-similar viscous contact profile $(\tilde\rho,\tilde u,\tilde\theta)(t,x_1)$ with initial value $(\tilde\rho_0,\tilde u_0,\tilde\theta_0)(x_1)$ by
	\begin{equation}
		\tilde\rho(t,x_1)=\tilde\rho_0\!\left(\frac{x_1}{\sqrt{1+t}}\right),\qquad
		\tilde u=-\frac{\kappa}{c_p}\frac{\tilde\rho_{x_1}}{\tilde\rho^2},\qquad
		\tilde\theta=\frac{p_*}{R\tilde\rho}-\frac{\tilde u^2}{2c_v}.
		\label{eq:self-similar-profile}
	\end{equation}
	Throughout the paper, $(\bar\rho,\bar u e_1,\bar\theta)(t,x_1)$ denotes the global one-dimensional solution to \eqref{eq:NSF} with initial data $(\tilde\rho_0,\tilde u_0 e_1,\tilde\theta_0)$, and $\bar p=R\bar\rho\bar\theta$. For a three-dimensional solution $(\rho,u,\theta)$ to \eqref{eq:NSF}, we define the perturbation $U=(\phi,\psi,\zeta)$ and the pressure and entropy variables $(q,h)$ by
	\begin{equation}
		\begin{aligned}
			&U=(\phi,\psi,\zeta)=(\rho-\bar\rho,u-\bar u e_1,\theta-\bar\theta),\\[2pt]
			&q=\log(p/\bar p), \qquad
			h=\gamma\bar\theta\left[
			(\theta/\bar\theta)^{1/\gamma}
			(\bar\rho/\rho)^{(\gamma-1)/\gamma}-1\right], \\[2pt]
			&\theta=e^{(\gamma-1)q/\gamma}(\bar\theta+h/\gamma), \qquad
			\rho=\frac{\bar p}{R(\bar\theta+h/\gamma)}e^{q/\gamma}.
		\end{aligned}
		\label{eq:qh}
	\end{equation}
	
	Our main result is stated as follows.
	
	\begin{theorem}
		\label{thm:main}
		Fix a sufficiently small relatively compact neighborhood $\mathcal K\subset(0,\infty)^2$
		of a positive equilibrium state, and let $(\rho_\pm,\theta_\pm)\in\mathcal K$
		satisfy $R\rho_-\theta_-=R\rho_+\theta_+$.
		There exist positive constants $\delta_0,\eps_0,c,C>0$, depending only
		on $\mathcal K$ and the physical parameters, such that for every $0<\delta\le\delta_0$, if the initial data
		\begin{equation}
			(\rho_0,u_0,\theta_0)(x)=(\tilde\rho_0,\tilde u_0 e_1,\tilde\theta_0)(x_1)+(\phi_0,\psi_0,\zeta_0)(x)
			\label{eq:3d-data}
		\end{equation}
		satisfy
		\begin{equation}
			U_0=(\phi_0,\psi_0,\zeta_0)\in H^4(\mathbb{R}^3),\qquad
			\eps=\|U_0\|_{H^4(\mathbb{R}^3)}\le\eps_0,
			\label{eq:3d-smallness}
		\end{equation}
		then the Cauchy problem \eqref{eq:NSF}, \eqref{eq:3d-data} admits a unique global strong solution $(\rho,u,\theta)$ satisfying
		\[
		\inf_{(t,x)\in[0,\infty)\times\R^3}\rho(t,x)>0,\qquad
		\inf_{(t,x)\in[0,\infty)\times\R^3}\theta(t,x)>0,
		\]
		and, on every finite time interval $[0,T]$,
		\begin{align*}
			&U\in C([0,T];H^4(\mathbb{R}^3)), \qquad
			(\psi,\zeta)\in L^2(0,T;H^5(\mathbb{R}^3)), \\[2pt]
			&\nabla\phi\in L^2(0,T;H^3(\mathbb{R}^3)), \qquad
			U_t\in L^2(0,T;H^3(\mathbb{R}^3)).
		\end{align*}
		Moreover, with $\eta=C\sqrt{\delta}<1/8$, the perturbation satisfies the global decay estimates
		\begin{align}
			&\|U(t)\|_\infty+\|h(t)\|_\infty
			\le C\eps(1+t)^{-1/4+\eta},\notag \\[2pt]
			&\|(q,\psi)(t)\|_\infty+\|(p-\bar p)(t)\|_\infty
			+\|\nabla U(t)\|_\infty+\|\nabla^2U(t)\|_\infty
			\le C\eps(1+t)^{-1/2+\eta},
			\label{eq:main-decay}
		\end{align}
		for all $t\ge 0$, together with the energy inequality
		\begin{equation}
			(1+t)^{-2\eta}\mathcal{E}(t)
			+c\int_0^t(1+s)^{-2\eta}\mathcal{D}(s)\,ds
			\le C\eps^2,
			\label{eq:main-energy}
		\end{equation}
		where $\mathcal N$ and $\mathcal D$ are defined in \eqref{eq:Nstar} and \eqref{eq:dissipation}, and $\mathcal{E}\sim\mathcal N$ is constructed in Proposition~\ref{prop:apriori-closure}. Here the constants $c,C$ and $\eps_0$ are independent of $\delta$ and $t$.
	\end{theorem}
	
	\begin{remark}
		In Theorem~\ref{thm:main}, the initial smallness threshold $\eps_0$ is independent of the wave strength $\delta$, and no $L^1(\R^3)$ integrability, spatial moment condition, or zero-mass assumption is imposed on $U_0$. Moreover, combining \eqref{eq:main-decay} with the background convergence bound \eqref{eq:bg-diff-decay} in Corollary~\ref{lem:bg-high-order} shows that the three-dimensional solution converges globally toward the self-similar contact wave in $L^\infty(\R^3)$:
		\[
		\|(\rho-\tilde\rho,u-\tilde u e_1,\theta-\tilde\theta)(t)\|_{L^\infty(\R^3)}
		\le C\eps(1+t)^{-1/4+\eta}+C\sqrt{\delta}\,(1+t)^{-11/16}\qquad \text{for all }t\ge 0.
		\]
	\end{remark}
	
	\subsection*{Notations and organization of the paper}
	
	Throughout the paper, $C$ and $c$ denote generic positive constants depending only on $\mathcal{K}$ and the physical parameters; dependence on auxiliary parameters is indicated by subscripts, and $A\sim B$ means that $C^{-1}B\le A\le CB$. Spatial norms of the three-dimensional perturbation are taken over $\mathbb{R}^3$, whereas norms of the one-dimensional background solution and its anti-derivatives are taken over $\mathbb{R}$. For scalar or vector-valued functions we write
	\[
	\|f\|_p=\|f\|_{L^p},\qquad
	\|(f,g)\|_p=\|f\|_p+\|g\|_p,\qquad
	\|\nabla^j f\|_2^2=\sum_{|\alpha|=j}\|\partial^\alpha f\|_2^2.
	\]
	
	The rest of the paper is organized as follows. In Section~\ref{sec:background-global}, we construct the one-dimensional background solution and establish its eighth-order time-weighted energy hierarchy. In Section~\ref{sec:log-spacetime}, we prove the Gaussian-weighted space-time estimate for the background material derivative $d$ and the one-dimensional weighted Sobolev inequality associated with $\omega$. In Section~\ref{sec:3d-low-order}, we derive the low-order estimates for the entropy variable $h$ and the velocity--pressure variables $(\psi,q)$. Section~\ref{sec:3d-high-order} contains the higher-order derivative estimates and the closure of the a priori energy inequality, and Section~\ref{sec:proof-main} completes the proof of Theorem~\ref{thm:main}.

	\section{The one-dimensional background solution}
	\label{sec:background-global}
	
	In this section, we construct the one-dimensional Navier--Stokes background solution $(\bar\rho,\bar u,\bar\theta)$ and establish its eighth-order time-weighted energy estimates. Throughout Sections~\ref{sec:background-global} and \ref{sec:log-spacetime}, we write $x=x_1\in\R$ and $\partial_x=\partial_{x_1}$, and set
	\[
	D_{\bar u}=\partial_t+\bar u\partial_x,\qquad
	D_{\tilde u}=\partial_t+\tilde u\partial_x,\qquad
	\bar p=R\bar\rho\bar\theta,\qquad
	d=D_{\bar u}\bar u=\bar u_t+\bar u\bar u_x.
	\]
	
	\subsection{The self-similar contact profile and the one-dimensional solution}
	
	By \eqref{eq:initial-density}--\eqref{eq:self-similar-profile}, the self-similar contact profile $(\tilde\rho,\tilde u,\tilde\theta)(t,x)$ satisfies
	\[
	\tilde\rho_t=\frac{\kappa}{c_p}\left(\frac{\tilde\rho_x}{\tilde\rho}\right)_x,\qquad
	\tilde\rho\Bigl(c_v\tilde\theta+\frac{1}{2}\tilde u^2\Bigr)=\frac{c_vp_*}{R}.
	\]
	
	\begin{lemma}
		\label{lem:profile-construction}
		For sufficiently small $\delta>0$, the boundary value problem \eqref{eq:initial-density} admits a unique smooth, positive, monotone solution $\tilde\rho_0(x)$ whose associated temperature $\tilde\theta(t,x)$ is strictly positive. Moreover, for every fixed finite range of derivative orders, there exist constants $C_{i,j},c_0>0$, uniform for $(\rho_\pm,\theta_\pm)\in\mathcal{K}$, such that
		\begin{equation}
			\left\{\begin{aligned}
				&|\tilde\rho-\rho_\pm|+|\tilde\theta-\theta_\pm|
				\le C\delta e^{-c_0x^2/(1+t)}\qquad (\pm x>0),\\[2pt]
				&|\partial_t^i\partial_x^j(\tilde\rho,\tilde\theta)|
				\le C_{i,j}\delta(1+t)^{-i-j/2}e^{-c_0x^2/(1+t)}
				\qquad (i+j\ge1),\\[2pt]
				&|\partial_t^i\partial_x^j\tilde u|
				\le C_{i,j}\delta(1+t)^{-i-(j+1)/2}e^{-c_0x^2/(1+t)}.
			\end{aligned}\right.
			\label{eq:profile-spacetime-decay}
		\end{equation}
	\end{lemma}
	\begin{proof}
		The existence and Gaussian decay of the initial profile follow from the nonlinear diffusion analysis in \cite{HMX,DHLX}. In terms of the logarithmic derivative $G=\tilde\rho_0'/\tilde\rho_0$, the profile ODE \eqref{eq:initial-density} reads
		\[
		\tilde\rho_0'=\tilde\rho_0G,\qquad
		G'=-\frac{c_p}{2\kappa}y\tilde\rho_0G,\qquad
		|G(y)|\le C\delta e^{-2c_0y^2}.
		\]
		{\color{black}By induction, every positive-order derivative of $\tilde\rho_0$ is $G$ times a polynomial in $y$, $\tilde\rho_0$ and $G$, where the latter two functions are bounded.} Absorbing the polynomial factors by slightly reducing $c_0>0$ and evaluating at the self-similar variable $y=x/\sqrt{1+t}$ yields \eqref{eq:profile-spacetime-decay}.
	\end{proof}
	
	\begin{proposition}
		\label{thm:background-global}
		Let $(\rho_\pm,\theta_\pm)\in\mathcal{K}$ satisfy $R\rho_-\theta_-=R\rho_+\theta_+$. For sufficiently small wave strength $\delta>0$, the one-dimensional compressible Navier--Stokes system with initial data
		\[
		(\bar\rho,\bar u,\bar\theta)(0,x)=(\tilde\rho_0,\tilde u_0,\tilde\theta_0)(x)
		\]
		admits a unique global strong solution $(\bar\rho,\bar u,\bar\theta)(t,x)$ satisfying
		\[
		\rho_*\le\bar\rho(t,x)\le\rho^*,\qquad
		\theta_*\le\bar\theta(t,x)\le\theta^*\qquad \text{for all }(t,x)\in[0,\infty)\times\R,
		\]
		and, on every finite time interval $[0,T]$,
		\[
		(\bar\rho-\tilde\rho,\bar u-\tilde u,\bar\theta-\tilde\theta)\in C([0,T];H^3(\R)),\qquad
		(\bar u-\tilde u,\bar\theta-\tilde\theta)\in L^2(0,T;H^4(\R)).
		\]
	\end{proposition}
	\begin{proof}
		Consider the problem in the periodic channel studied in \cite[Theorem~2.4]{DHLX} with initial data $(\tilde\rho_0,\tilde u_0 e_1,\tilde\theta_0)$ depending only on $x_1$. Since the initial perturbation relative to the contact profile is zero, its Sobolev norms, excess mass, and initial anti-derivative all vanish, so the assumptions of \cite[Theorem~2.4]{DHLX} are satisfied. By uniqueness, the global strong solution remains one-dimensional with zero transverse velocity, which gives $(\bar\rho,\bar u,\bar\theta)$ and the uniform upper and lower bounds for $\bar\rho$ and $\bar\theta$.
	\end{proof}
	
	To derive higher-order estimates for $(\bar\rho,\bar u,\bar\theta)$, we write the one-dimensional conserved quantities and fluxes as
	\[
	Q=\begin{pmatrix}\bar\rho\\ \bar\rho\bar u\\ \bar\rho(c_v\bar\theta+\bar u^2/2)\end{pmatrix},\quad
	\mathcal F(Q)=\begin{pmatrix}\bar\rho\bar u\\ \bar\rho\bar u^2+\bar p\\ \bar u\bigl(\bar\rho(c_v\bar\theta+\bar u^2/2)+\bar p\bigr)\end{pmatrix},\quad
	\mathcal V(Q,Q_x)=\begin{pmatrix}0\\ \nu\bar u_x\\ \kappa\bar\theta_x+\nu\bar u\bar u_x\end{pmatrix}.
	\]
	While $Q$ satisfies the homogeneous system, the self-similar profile $\tilde Q=\bigl(\tilde\rho,\tilde\rho\tilde u,\tilde\rho(c_v\tilde\theta+\tilde u^2/2)\bigr)^{\top}$ satisfies
	\begin{equation}
		Q_t+\partial_x\mathcal{F}(Q)-\partial_x\mathcal{V}(Q,Q_x)=0,\qquad
		\tilde Q_t+\partial_x\mathcal F(\tilde Q)
		-\partial_x\mathcal V(\tilde Q,\tilde Q_x)
		=\partial_x\mathcal S,
		\label{eq:tildeQ-approx}
	\end{equation}
	where $\mathcal S=(0,r_1,r_2)^{\top}$ with
	\begin{equation}
		\begin{aligned}
			r_1&=-\frac{\kappa}{c_p}\frac{\tilde\rho_t}{\tilde\rho}
			+\frac{3-\gamma}{2}\tilde\rho\tilde u^2-\nu\tilde u_x,\\
			r_2&=-\frac{\gamma-1}{2}\tilde\rho\tilde u^3
			+\left(\frac{\kappa}{c_v}-\nu\right)\tilde u\tilde u_x.
		\end{aligned}
		\label{eq:Q1-Q2}
	\end{equation}
	We define the conservative anti-derivative vector $Y=(\Phi,\Psi,\mathcal W)^{\top}$ by integrating the flux difference in time:
	\begin{equation}
		Y(t,x)=\int_0^t\!\left\{-\mathcal F(Q)+\mathcal F(\tilde Q)
		+\mathcal V(Q,Q_x)-\mathcal V(\tilde Q,\tilde Q_x)-\mathcal S\right\}(s,x)\,ds.
		\label{eq:Y-def}
	\end{equation}
	By Proposition~\ref{thm:background-global} and Lemma~\ref{lem:profile-construction}, differentiating \eqref{eq:Y-def} with respect to $x$ and using \eqref{eq:tildeQ-approx} and the vanishing initial difference between $Q$ and $\tilde Q$ yields
	\begin{equation}
		Y_x=Q-\tilde Q,\qquad Y(0)=0,\qquad
		Y\in C([0,T];H^4(\R))\cap H^1(0,T;H^3(\R)).
		\label{eq:Y-reg}
	\end{equation}
	
	\subsection{Transformed anti-derivative variables and structural cancellation}
	
	Set
	\[
	\mathcal{J}(Q)=\Bigl(0,\nu u,\kappa\theta+\frac{1}{2}\nu u^2\Bigr)^{\top},\qquad
	\mathcal{B}(Q)=D_Q\mathcal{J}(Q),\qquad
	\mathcal{B}(Q)Q_x=\mathcal{V}(Q,Q_x).
	\]
	By \eqref{eq:Y-def}--\eqref{eq:Y-reg}, the conservative anti-derivative $Y=(\Phi,\Psi,\mathcal{W})^{\top}$ satisfies the integrated system
	\begin{equation}
		Y_t+\mathcal{F}(\tilde Q+Y_x)-\mathcal{F}(\tilde Q)
		=\mathcal{B}(Q)(\tilde Q_x+Y_{xx})-\mathcal{B}(\tilde Q)\tilde Q_x
		-\mathcal{S}.
		\label{eq:Y-system}
	\end{equation}
	In Eulerian coordinates, the system \eqref{eq:Y-system} does not satisfy the left structural condition in \eqref{eq:LX-conditions}, so differentiating \eqref{eq:Y-system} directly produces non-integrable profile-gradient couplings. Following \cite{DHLX}, we introduce the transformed anti-derivative vector $V=(F,Z,H)^{\top}=MY$ with inverse $Y=NV$, where
	\begin{equation}
		M=\begin{pmatrix}
			c_v\tilde\theta & 0 & 0\\[2pt]
			-\tilde u & 1 & 0\\[2pt]
			\frac{1}{2}\tilde u^2-c_v\tilde\theta & -\tilde u & 1
		\end{pmatrix}, \qquad
		N=M^{-1}=\begin{pmatrix}
			(c_v\tilde\theta)^{-1} & 0 & 0\\[2pt]
			\tilde u(c_v\tilde\theta)^{-1} & 1 & 0\\[2pt]
			1+\frac{1}{2}\tilde u^2(c_v\tilde\theta)^{-1} & \tilde u & 1
		\end{pmatrix},
		\label{eq:M-N-matrices}
	\end{equation}
	or in component form,
	\begin{equation}
		F=c_v\tilde\theta\Phi, \qquad
		Z=\Psi-\tilde u\Phi, \qquad
		H=\mathcal{W}-\tilde u Z-\left(c_v\tilde\theta+\frac{1}{2}\tilde u^2\right)\Phi.
		\label{eq:FZH-def}
	\end{equation}
	
	Let $A(\tilde Q)=D_Q\mathcal F(\tilde Q)$ be the inviscid Jacobian matrix evaluated at $\tilde Q$. Substituting $Y=NV$ into \eqref{eq:Y-system} and multiplying on the left by $M$ gives
	\begin{equation}
		V_t+A_1V_x=A_2V_{xx}+C_1V_x+C_0V+R_V+\partial_x D_V+\mathcal P,
		\label{eq:V-system}
	\end{equation}
	where
	\begin{equation}
		\begin{aligned}
			A_1&=MA(\tilde Q)N-\tilde u I
			=\begin{pmatrix}
				0&c_v\tilde\theta&0\\
				\gamma-1&0&\gamma-1\\
				0&(\gamma-1)c_v\tilde\theta&0
			\end{pmatrix},\\[2pt]
			A_2&=M\mathcal B(\tilde Q)N
			=\operatorname{diag}\!\left(0,\frac{\nu}{\tilde\rho},
			\frac{\kappa}{c_v\tilde\rho}\right).
		\end{aligned}
		\label{eq:MAcN}
	\end{equation}
	Denoting the Taylor remainders of $\mathcal{F}$ and $\mathcal{J}$ around $\tilde Q$ by
	\begin{equation}
		\mathcal{F}_2=\mathcal{F}(\tilde Q+Y_x)-\mathcal{F}(\tilde Q)-A(\tilde Q)Y_x, \qquad
		\mathcal{J}_2=\mathcal{J}(\tilde Q+Y_x)-\mathcal{J}(\tilde Q)-\mathcal{B}(\tilde Q)Y_x,
		\label{eq:F2-J2-def}
	\end{equation}
	the remaining terms in \eqref{eq:V-system} are given by
	\begin{equation}
		\left\{\begin{aligned}
			&C_1=-\tilde u I+2M\mathcal{B}(\tilde Q)N_x
			+M[\mathcal{B}(\tilde Q)]_xN, \\[2pt]
			&C_0=M_tN-MA(\tilde Q)N_x+M\mathcal{B}(\tilde Q)N_{xx}
			+M[\mathcal{B}(\tilde Q)]_xN_x, \\[2pt]
			&D_V=M \mathcal{J}_2, \qquad
			R_V=-M \mathcal{F}_2-M_x \mathcal{J}_2, \qquad
			\mathcal{P}=(0,-r_1,\tilde ur_1-r_2)^{\top}.
		\end{aligned}\right.
		\label{eq:V-coeffs}
	\end{equation}
	
	Under the transformation \eqref{eq:FZH-def}, the terms involving $\tilde\theta_x$ cancel in the matrix $C_0$. Indeed, from \eqref{eq:M-N-matrices} and \eqref{eq:MAcN} we have
	\[
	M_tN-MA(\tilde Q)N_x=-M(N_t+\tilde uN_x)-A_1MN_x,\qquad
	MN_x=\begin{pmatrix}
		-\tilde\theta_x/\tilde\theta&0&0\\
		\tilde u_x/(c_v\tilde\theta)&0&0\\
		\tilde\theta_x/\tilde\theta&\tilde u_x&0
	\end{pmatrix}.
	\]
	Since the first and third entries in the second row of $A_1$ are both equal to $\gamma-1$, the entries $\mp\tilde\theta_x/\tilde\theta$ in the first column of $MN_x$ cancel in the product $A_1MN_x$, and all remaining entries involve $\tilde u_x$, $D_{\tilde u}\tilde\theta$ {\color{black}or $D_{\tilde u}\tilde u$ (the latter is $O(\delta\tau^{-3/2})$ with Gaussian decay, where $\tau=1+t$)}. By Lemma~\ref{lem:profile-construction}, $C_1$ and $C_0$ satisfy
	\begin{equation}
		\begin{aligned}
			&|\partial_x^j C_1(t,x)|\le C\delta (1+t)^{-(j+1)/2}e^{-c_0 x^2/(1+t)}, \\[2pt]
			&|\partial_x^j C_0(t,x)|\le C\delta (1+t)^{-(j+2)/2}e^{-c_0 x^2/(1+t)},
		\end{aligned}
		\label{eq:C1-C0-bounds}
	\end{equation}
	{\color{black}for each fixed finite range of $j\ge0$, with constants depending on that range} (after slightly reducing $c_0>0$).
	
	\subsection{Time-weighted energy hierarchy and physical background bounds}
	\label{subsec:background-high-order}
	
	We next establish time-weighted energy estimates for $V=(F,Z,H)^{\top}$ up to order eight. Set $\tau=1+t$. The matrix $A_1$ in \eqref{eq:V-system} is diagonalized as
	\[
	LA_1L^{-1}=\Lambda_0=\operatorname{diag}(-\tilde c,0,\tilde c),
	\qquad \tilde c=\sqrt{\gamma R\tilde\theta},
	\]
	where the right and left eigenvector matrices $L^{-1}$ and $L$ are given explicitly by
	\begin{align}
		L^{-1}&=\begin{pmatrix}
			(2\gamma)^{-1/2} & \sqrt{(\gamma-1)/\gamma} & (2\gamma)^{-1/2}\\[2pt]
			-\sqrt{(\gamma-1)/(2c_v\tilde\theta)} & 0 & \sqrt{(\gamma-1)/(2c_v\tilde\theta)}\\[2pt]
			(\gamma-1)(2\gamma)^{-1/2} & -\sqrt{(\gamma-1)/\gamma} & (\gamma-1)(2\gamma)^{-1/2}
		\end{pmatrix}, \label{eq:Linv-matrix}\\[4pt]
		L&=\begin{pmatrix}
			(2\gamma)^{-1/2} & -\sqrt{c_v\tilde\theta/(2(\gamma-1))} & (2\gamma)^{-1/2}\\[2pt]
			\sqrt{(\gamma-1)/\gamma} & 0 & -(\gamma(\gamma-1))^{-1/2}\\[2pt]
			(2\gamma)^{-1/2} & \sqrt{c_v\tilde\theta/(2(\gamma-1))} & (2\gamma)^{-1/2}
		\end{pmatrix}. \label{eq:L-matrix}
	\end{align}
	In terms of $B=(B_1,B_2,B_3)^{\top}=LV$, the components $B_1$ and $B_3$ correspond to the acoustic modes with speeds $\mp\tilde c$, and $B_2$ corresponds to the contact mode. Since the second row of $L$ and the second column of $L^{-1}$ are constant, the matrices $L_xL^{-1}$ and $L_tL^{-1}$ have zero second row and second column, which reflects both conditions in \eqref{eq:LX-conditions}.
	
	Before defining the weighted energy functionals, we record the algebraic relations connecting $(F,Z,H)$ and $B=LV$ to the physical differences $\hat u=\bar u-\tilde u$ and $\bar p-\tilde p$.
	
	\begin{lemma}
		\label{lem:anti-identities}
		Denoting $\hat u=\bar u-\tilde u$, the transformed anti-derivative variables $(F,Z,H)$ defined in \eqref{eq:FZH-def} satisfy the identities
		\begin{align}
			&\bar\rho\hat u=Z_x+\tilde u_x\Phi,\qquad \Phi=\frac{F}{c_v\tilde\theta},\label{eq:uhat-identity} \\[2pt]
			&\bar p-\tilde p
			=(\gamma-1)\left[(F+H)_x+\tilde u_xZ-\frac{1}{2}\bar\rho\hat u^2\right].
			\label{eq:phat-identity}
		\end{align}
		Moreover, in terms of the diagonalized vector $B=(B_1,B_2,B_3)^{\top}=LV$, we have
		\begin{equation}
			F+H=\sqrt{\gamma/2}\,(B_1+B_3).
			\label{eq:acoustic-sum}
		\end{equation}
	\end{lemma}
	\begin{proof}
		Differentiating the second relation in \eqref{eq:FZH-def} with respect to $x$ yields \eqref{eq:uhat-identity}. Next, differentiating the third relation in \eqref{eq:FZH-def} and expanding the kinetic energy difference via \eqref{eq:uhat-identity} gives
		\begin{align*}
			\mathcal W_x&=(F+H)_x+(\tilde uZ)_x+\frac{1}{2}(\tilde u^2\Phi)_x,\\[2pt]
			\frac{1}{2}(\bar\rho\bar u^2-\tilde\rho\tilde u^2)
			&=\tilde u\bar\rho\hat u+\frac{1}{2}\tilde u^2\Phi_x+\frac{1}{2}\bar\rho\hat u^2
			=\tilde u Z_x+\frac{1}{2}(\tilde u^2\Phi)_x+\frac{1}{2}\bar\rho\hat u^2.
		\end{align*}
		Subtracting these two identities establishes \eqref{eq:phat-identity}, and multiplying $V=L^{-1}B$ on the left by $(1,0,1)$ via \eqref{eq:Linv-matrix} yields \eqref{eq:acoustic-sum}.
	\end{proof}
	
	Fix constants $0<\beta_0<\beta/4$ and $4\beta<c_0$. For $0\le k\le 8$, we define the symmetrizer $S$ and the weight matrix $W$ by
	\begin{equation}
		\begin{aligned}
			&S=\operatorname{diag}\!\left(\frac{\gamma-1}{c_v\tilde\theta},1,
			\frac1{c_v\tilde\theta}\right),\qquad
			g_\beta=\tau^{-1/2}e^{-\beta x^2/\tau},\\
			&\vartheta_\beta(t,x)=\int_0^xg_\beta(t,y)\,dy,\qquad
			H_\beta=\operatorname{diag}(\vartheta_\beta,0,-\vartheta_\beta),
			\qquad W=S+\kappa_1L^{\top}H_\beta L.
		\end{aligned}
		\label{eq:W-matrix-def}
	\end{equation}
	Here $S$ is chosen so that
	\[
	SA_1=\begin{pmatrix}
		0 & \gamma-1 & 0\\[2pt]
		\gamma-1 & 0 & \gamma-1\\[2pt]
		0 & \gamma-1 & 0
	\end{pmatrix}
	\]
	is symmetric and constant. Since $\partial_x\vartheta_\beta=g_\beta>0$, the weight matrix $H_\beta$ uses the opposite signs of $\mp\tilde c$ in $\Lambda_0$ to produce a negative-definite quadratic form on the acoustic components:
	\[
	\partial_x(H_\beta\Lambda_0)
	=-g_\beta\tilde c\operatorname{diag}(1,0,1)+{\color{black}H_\beta\partial_x\Lambda_0}
	\le -(c_*-C_\beta\delta)g_\beta\operatorname{diag}(1,0,1),\qquad
	c_*=\inf\tilde c>0.
	\]
	To obtain dissipation for $\partial_x^{k+1}F$, we also set
	\begin{equation}
		\begin{aligned}
			&K_k=\|\partial_x^{k+1}V\|_2^2,\qquad
			\mathcal A_k=\int_{\R}g_\beta\bigl(|(L\partial_x^k V)_1|^2+|(L\partial_x^k V)_3|^2\bigr)\,dx,\qquad
			a=\frac{\nu}{\bar\rho c_v\tilde\theta},\\
			&\mathcal C_k=\int_{\R}(\partial_x^k Z)(\partial_x^{k+1}F)\,dx+\frac12\int_{\R}a|\partial_x^{k+1}F|^2\,dx,
			\qquad
			\mathcal E_k=\frac12\int_{\R}(\partial_x^k V)^{\top}W(\partial_x^k V)\,dx+\kappa_2\mathcal C_k.
		\end{aligned}
		\label{eq:Ek-def}
	\end{equation}
	For fixed small constants $0<\kappa_1\ll\kappa_2\ll1$, we have
	\begin{equation}
		\mathcal E_k\sim\|\partial_x^k V\|_2^2+\|\partial_x^{k+1}F\|_2^2,
		\qquad 0\le k\le8.
		\label{eq:Ek-equiv}
	\end{equation}
	
	\begin{lemma}
		\label{lem:1d-unified}
		For the one-dimensional solution in Proposition~\ref{thm:background-global}
		with $V(0)=0$ and sufficiently small $\delta>0$, the functionals in \eqref{eq:Ek-def} satisfy
		\begin{equation}
			\tau^r\sum_{k=r}^8\mathcal E_k(t)
			+\int_0^t(1+s)^r\sum_{k=r}^8
			\bigl(K_k+\kappa_1\mathcal A_k\bigr)(s)\,ds
			\le C\delta\log(2+t)\le C\delta\tau^{1/8},
			\qquad r=0,\ldots,8.
			\label{eq:1d-hierarchy}
		\end{equation}
	\end{lemma}
	
	\begin{proof}
		{\color{black}
			Apply the closed estimates in \cite[proof of Theorem~2.4,
			p.~31, (4.17)--(4.18)]{DHLX} with zero initial difference. Replacing $\theta$ by $c_v\theta$ and $\kappa$
			by $\kappa/c_v$ gives its normalization, and the nonzero transverse modes vanish.
			Their energy equivalences give
			\begin{gather*}
				\|V(t)\|_2^2+\int_0^t(K_0+K_1)(s)\,ds\le C\delta\log(2+t),\\
				\|\partial_x^jV(t)\|_2^2\le C\delta\tau^{-j}\log(2+t),\qquad j=1,2.
			\end{gather*}
			Hence one-dimensional interpolation gives
			$\|V\|_\infty+\|V_x\|_\infty\le C\sqrt\delta$,
			without using the estimates below.
			For the fixed Gaussian exponent $\beta_0$, the heat-kernel argument of
			\cite[Lemmas~4.4--4.5]{DHLX} gives
			\begin{equation}
				I_0(t)=\int_{\R}e^{-\beta_0x^2/\tau}|V(t,x)|^2\,dx,\qquad
				\int_0^t(1+s)^{-1}I_0(s)\,ds\le C\delta\log(2+t).
				\label{eq:dhlx-local-energy}
			\end{equation}
			Indeed, the preceding bound makes the endpoint and weight-time terms
			$O(\delta)$, since $\tau^{-1/2}\log(2+t)$ is bounded and
			$\int_0^\infty(1+s)^{-3/2}\log(2+s)\,ds<\infty$.
			The acoustic slow source contributes
			$C_{\beta_0}\delta^2\tau^{-3/2}$ in the weighted source estimate.
			The contact projection removes this source and leaves
			$\tilde u r_1-r_2=(\kappa/c_p-\kappa/c_v)\tilde u\tilde u_x$;
			its weighted squared source term in the heat-kernel estimate is bounded by
			$C_{\beta_0}\delta^4\tau^{-3/2}$.
			The remaining pairings use these low-order bounds, giving
			\[
			\int_0^t(1+s)^{-1}I_0(s)\,ds
			\le C_{\beta_0}\delta\log(2+t)
			+C_{\beta_0}\delta\int_0^t(1+s)^{-1}I_0(s)\,ds.
			\]
			Fix $\beta_0$ and take $C_{\beta_0}\delta\le1/2$ to obtain
			\eqref{eq:dhlx-local-energy}, with $C$ independent of $\delta$.
			
			The finite-time bounds in Proposition~\ref{thm:background-global}, the smooth
			initial profile and regularization justify the following differentiations.
			Under the a priori assumption
			$\max_{0\le j\le2}\tau^{j/2}\|\partial_x^jV\|_\infty\le\varepsilon_1$,
			differentiating \eqref{eq:uhat-identity} gives
			\begin{equation}
				\begin{aligned}
					\nu\partial_x^n\hat u
					&=\frac{\nu}{\bar\rho}\partial_x^{n+1}Z
					-\frac{\nu\hat u}{\bar\rho c_v\tilde\theta}\partial_x^{n+1}F+f_n,
					\qquad 1\le n\le9,\\
					\|f_n\|_2&\le C(\delta+\varepsilon_1)
					\sum_{j=1}^n\tau^{-(n+1-j)/2}\|\partial_x^jV\|_2
					+C\delta\tau^{-(n+1)/2}I_0^{1/2}.
				\end{aligned}
				\label{eq:viscous-top}
			\end{equation}
			The remainder follows from Leibniz' rule and one-dimensional
			Gagliardo--Nirenberg interpolation after scaling space by $\sqrt\tau$.
			The coefficient commutators are estimated similarly, retaining a profile
			Gaussian on undifferentiated $V$. The derivatives
			$\partial_x^k(R_V,D_V)$ contain at most $\partial_x^{k+1}V$;
			after integrating the $\partial_xD_V$ term by parts once, their highest-order
			terms are absorbed by $K_k$ under the same smallness assumption.
			
			Differentiate \eqref{eq:V-system} $k$ times, test with $W\partial_x^kV$,
			and add $\kappa_2\mathcal C_k'$.
			The mass and momentum equations cancel the
			$\partial_x^{k+2}Z\,\partial_x^{k+1}F$ terms because
			$a c_v\tilde\theta=\nu/\bar\rho$.
			The remaining highest density product is integrated by parts;
			the pressure term and $SA_2$ provide the full $K_k$ dissipation.
			The transport coefficient satisfies
			$\|a_t+(a\bar u)_x\|_\infty\le C(\delta+\varepsilon_1)\tau^{-1}$.
			The sign of $\partial_x(H_\beta\Lambda_0)$ gives $\kappa_1\mathcal A_k$,
			while $|\partial_tH_\beta|\le C\tau^{-1}e^{-\beta_0x^2/\tau}$
			and the parabolic cross term is absorbed for $\kappa_1^2\ll\kappa_2$.
			The slow source is acoustic and satisfies
			$\int_{\R}g_\beta^{-1}|\partial_x^kr_1|^2\,dx
			\le C\delta^2\tau^{-k-1}$.
			Derivatives of $A_1$ also couple only the acoustic components.
			Leibniz' rule, \eqref{eq:C1-C0-bounds} and Young's inequality therefore give
			\begin{align}
				\mathcal E_0'+cK_0+c\kappa_1\mathcal A_0
				&\le C(\delta+\kappa_1)\tau^{-1}I_0
				+C_{\kappa_1}\delta\tau^{-1},\label{eq:1d-E0}\\
				\mathcal E_k'+cK_k+c\kappa_1\mathcal A_k
				&\le C(\delta+\varepsilon_1+\kappa_1)
				\sum_{j<k}\tau^{-(k-j)}K_j
				+\frac{C\delta^2}{\kappa_1}
				\sum_{1\le j<k}\tau^{-(k-j)}\mathcal A_j\notag\\
				&\quad+C\delta\tau^{-k-1}I_0
				+C_{\kappa_1}\delta\tau^{-k-1},\qquad 1\le k\le8.
				\label{eq:1d-Ek}
			\end{align}
			Choose $\kappa_2$, then $\kappa_1$, and then $\varepsilon_1,\delta$
			small. Summation and \eqref{eq:dhlx-local-energy} give
			\eqref{eq:1d-hierarchy} for $r=0$.
			For $r=1,\ldots,8$, multiply the equations with $k\ge r$ by $\tau^r$.
			Terms with $j\ge r$ are absorbed; for $j<r$, use
			$\tau^{r-(k-j)}\le\tau^j$ and the preceding weight $j$.
			The weight derivative uses $\mathcal E_k\le C(K_{k-1}+K_k)$
			and weight $r-1$, while $\tau^{r-k-1}\le\tau^{-1}$ handles the sources.
			This proves \eqref{eq:1d-hierarchy}. Finally,
			\begin{equation}
				\|\partial_x^jV(t)\|_\infty
				\le\sqrt2\,\|\partial_x^jV(t)\|_2^{1/2}
				\|\partial_x^{j+1}V(t)\|_2^{1/2}
				\le C\sqrt\delta\,\tau^{-j/2-1/4}[\log(2+t)]^{1/2},
				\qquad 0\le j\le7,
				\label{eq:V-Linfty}
			\end{equation}
			which closes the a priori assumption for sufficiently small $\delta$.
		}
	\end{proof}
	
	\begin{corollary}
		\label{lem:bg-high-order}
		For sufficiently small wave strength $\delta>0$, the one-dimensional background solution $(\bar\rho,\bar u,\bar\theta)$ satisfies, with $\sigma=C\sqrt\delta$,
		\begin{align}
			&\|\partial_x^j(\bar\rho,\bar\theta)(t)\|_\infty
			\le C\sigma\tau^{-j/2},\qquad 1\le j\le6,
			\label{eq:bg-thermal}\\
			&\|\partial_x^j\bar u(t)\|_\infty
			\le C\sigma\tau^{-(j+1)/2},\qquad 0\le j\le6,
			\label{eq:bg-velocity}\\
			&\|\partial_x^jD_{\bar u}(\bar\rho,\bar\theta)(t)\|_\infty
			+\|\partial_x^j(\bar p_x,d)(t)\|_\infty
			\le C\sigma\tau^{-1-j/2},\qquad 0\le j\le4.
			\label{eq:bg-target}
		\end{align}
		Moreover, the difference between $(\bar\rho,\bar u,\bar\theta)$ and the self-similar contact profile satisfies
		\begin{equation}
			\|\partial_x^j(\bar\rho-\tilde\rho,\bar u-\tilde u,
			\bar\theta-\tilde\theta)(t)\|_\infty
			\le C\sqrt\delta\,\tau^{-j/2-11/16},\qquad 0\le j\le6.
			\label{eq:bg-diff-decay}
		\end{equation}
	\end{corollary}
	
	\begin{proof}
		Combining \eqref{eq:V-Linfty} with $Y_x=NV_x+N_xV$, \eqref{eq:uhat-identity}--\eqref{eq:phat-identity}, and Lemma~\ref{lem:profile-construction}, we find
		\[
		\|\partial_x^j(\bar\rho-\tilde\rho,\bar u-\tilde u,
		\bar\theta-\tilde\theta)(t)\|_\infty
		\le C\sqrt\delta\,\tau^{-j/2-3/4}\bigl[\log(2+t)\bigr]^{1/2},\qquad 0\le j\le6.
		\]
		Since $[\log(2+t)]^{1/2}\le C\tau^{1/16}$, this gives \eqref{eq:bg-diff-decay}, and adding \eqref{eq:profile-spacetime-decay} yields \eqref{eq:bg-thermal} and \eqref{eq:bg-velocity}. To prove \eqref{eq:bg-target}, {\color{black}first use \eqref{eq:phat-identity} and \eqref{eq:V-Linfty} to retain the pressure cancellation:
			\[
			\|\partial_x^{j+1}(\bar p-\tilde p)\|_\infty
			\le C\sqrt\delta\,\tau^{-j/2-5/4}[\log(2+t)]^{1/2}
			\le C\sqrt\delta\,\tau^{-j/2-19/16},\qquad 0\le j\le4.
			\]
			Together with the profile pressure relation, this yields the pressure bound below.}
		We then use the one-dimensional Navier--Stokes equations:
		\[
		\begin{aligned}
			&D_{\bar u}\bar\rho=-\bar\rho\bar u_x,\qquad
			\bar\rho d=\nu\bar u_{xx}-\bar p_x,\qquad
			c_v\bar\rho D_{\bar u}\bar\theta=\kappa\bar\theta_{xx}+\nu\bar u_x^2-\bar p\bar u_x,\\[2pt]
			&\tilde p=p_*-\frac{\gamma-1}{2}\tilde\rho\tilde u^2,\qquad
			\|\partial_x^{j+1}\bar p(t)\|_\infty\le C\sigma\tau^{-1-j/2}\quad (0\le j\le 4).
		\end{aligned}
		\]
		Differentiating the three evolution equations up to four times with respect to $x$ and substituting \eqref{eq:bg-thermal}--\eqref{eq:bg-velocity} yields \eqref{eq:bg-target}.
	\end{proof}
	
	\section{Logarithmic space-time estimates and weighted Sobolev inequalities}
	\label{sec:log-spacetime}
	
	\subsection{Gaussian-weighted acoustic dissipation}
	
	In the three-dimensional stability analysis, controlling the coupling with the background material derivative $d=D_{\bar u}\bar u$ in the contact layer requires a Gaussian-weighted space-time integral bound for $d$ with a Gaussian width parameter $0<\beta<c_0/4$. We deduce this estimate from the acoustic dissipation $\mathcal{A}_k(t)$ established in Lemma~\ref{lem:1d-unified}.
	
	\begin{corollary}
		\label{cor:1d-gaussian-B}
		For each fixed $0<\beta<c_0/4$, there exists $\delta_0(\beta)>0$ such that for all $0<\delta\le\delta_0(\beta)$ and $t\ge 0$,
		\begin{equation}
			\int_0^t(1+s)^{3/2}\int_{\R}e^{-\beta x^2/(1+s)}
			\bigl(|(B_1)_{xx}|^2+|(B_3)_{xx}|^2\bigr)\,dx\,ds
			\le C_\beta\delta\log(2+t).
			\label{eq:B-gaussian-bound}
		\end{equation}
	\end{corollary}
	\begin{proof}
		Set $\tau=1+t$. Differentiating $B=LV$ twice with respect to $x$ and using \eqref{eq:L-matrix} gives
		\[
		B_{xx}=LV_{xx}+2L_xV_x+L_{xx}V,\qquad
		Z=-\sqrt{\frac{\gamma-1}{2c_v\tilde\theta}}\,(B_1-B_3).
		\]
		By \eqref{eq:L-matrix}, only the middle column of $L_x$ and $L_{xx}$ is nonzero, and that column multiplies the second component $Z$ of $V$, which is a linear combination of the acoustic components $B_1$ and $B_3$. Using the Gaussian decay of the profile temperature derivatives in Lemma~\ref{lem:profile-construction}, we obtain
		\[
		\tau^2\int_{\R}g_\beta\bigl(|(B_1)_{xx}|^2+|(B_3)_{xx}|^2\bigr)\,dx
		\le C_\beta\bigl(\tau^2\mathcal A_2
		+\delta^2\tau\mathcal A_1
		+\delta^2\mathcal A_0\bigr).
		\]
		Integrating this inequality over $[0,t]$ and applying \eqref{eq:1d-hierarchy} with $r=0,1,2$ yields \eqref{eq:B-gaussian-bound}.
	\end{proof}
	
	\subsection{Logarithmic space-time bounds for the background material derivative}
	
	We now combine Lemma~\ref{lem:1d-unified} and Corollary~\ref{cor:1d-gaussian-B} to establish the space-time bounds for $d$. Substituting \eqref{eq:phat-identity} and \eqref{eq:acoustic-sum} from Lemma~\ref{lem:anti-identities} into the one-dimensional momentum equation, we decompose $d=(\nu\bar u_{xx}-\bar p_x)/\bar\rho$ into its principal acoustic part and a remainder $f_d$:
	\begin{align}
		d&=-\frac{(\gamma-1)\sqrt{\gamma/2}}{\bar\rho}(B_1+B_3)_{xx}+f_d,
		\label{eq:d-decomp} \\[2pt]
		f_d&=\frac{1}{\bar\rho}\left[
		\nu \hat u_{xx}+\nu\tilde u_{xx}-\tilde p_x
		-(\gamma-1)(\tilde u_{xx}Z+\tilde u_xZ_x)
		+\frac{\gamma-1}{2}(\bar\rho\hat u^2)_x\right].
		\label{eq:fd-def}
	\end{align}
	Using \eqref{eq:1d-hierarchy} with $r=0,1,2$ together with \eqref{eq:uhat-identity} and one-dimensional Gagliardo--Nirenberg interpolation, we have
	\begin{equation}
		(1+t)^{1/2}\|\hat u(t)\|_2+(1+t)\|\hat u_x(t)\|_2
		+(1+t)^{3/4}\|\hat u(t)\|_\infty
		\le C\sqrt{\delta\log(2+t)}.
		\label{eq:uhat-rates}
	\end{equation}
	
	\begin{proposition}
		\label{prop:1d-d-estimates}
		Under the bounds \eqref{eq:1d-hierarchy} with $r=0,1,2$ and \eqref{eq:bg-thermal}--\eqref{eq:bg-velocity}, the background material derivative $d$ and the remainder $f_d$ in \eqref{eq:fd-def} satisfy
		\begin{align}
			&\int_0^t(1+s)\|d(s)\|_2^2\,ds
			\le C\delta\log(2+t),
			\label{eq:d-L2-bound} \\[2pt]
			&\int_0^t(1+s)^{3/2}\|f_d(s)\|_2^2\,ds
			\le C\delta\log(2+t).
			\label{eq:fd-L2-bound}
		\end{align}
	\end{proposition}
	\begin{proof}
		Differentiating the velocity identity \eqref{eq:uhat-identity} twice in $x$ and using the density derivative bounds in \eqref{eq:bg-thermal}, we obtain
		\begin{equation}
			\|\hat u_{xx}\|_2
			\le C\Bigl(
			\|V_{xxx}\|_2+\sigma (1+t)^{-1/2}\|V_{xx}\|_2
			+\sigma (1+t)^{-1}\|V_x\|_2+\delta (1+t)^{-2}\|V\|_2
			\Bigr).
			\label{eq:uhat-xx}
		\end{equation}
		Squaring \eqref{eq:uhat-xx}, multiplying by $(1+s)^{3/2}$, and integrating over $[0,t]$ via \eqref{eq:1d-hierarchy} with $r=0,1,2$ gives
		\begin{align*}
			\int_0^t (1+s)^{3/2}\|\hat u_{xx}(s)\|_2^2\,ds
			&\le C\int_0^t (1+s)^2\|V_{xxx}(s)\|_2^2\,ds
			+C\sigma^2\delta\int_0^\infty (1+s)^{-3/2}\log(2+s)\,ds \\[2pt]
			&\quad +C\delta^3\int_0^\infty (1+s)^{-5/2}\log(2+s)\,ds
			\le C\delta\log(2+t).
		\end{align*}
		Next, by Lemma~\ref{lem:profile-construction} and \eqref{eq:uhat-rates}, the remaining terms in $f_d$ satisfy
		\begin{align*}
			&\|\nu\tilde u_{xx}-\tilde p_x\|_2 \le C\delta (1+t)^{-5/4},\\[2pt]
			&\|\tilde u_{xx}Z\|_2+\|\tilde u_xZ_x\|_2
			\le C\delta^{3/2} (1+t)^{-3/2}\sqrt{\log(2+t)},\\[2pt]
			&\|(\bar\rho\hat u^2)_x\|_2
			\le C\|\bar\rho_x\|_\infty\|\hat u\|_\infty\|\hat u\|_2
			+C\|\hat u\|_\infty\|\hat u_x\|_2
			\le C\delta (1+t)^{-7/4}\log(2+t).
		\end{align*}
		Squaring these three bounds, multiplying by $(1+s)^{3/2}$, and integrating in time proves \eqref{eq:fd-L2-bound}. Combining \eqref{eq:d-decomp}, \eqref{eq:acoustic-sum}, \eqref{eq:1d-hierarchy} with $r=1$, and \eqref{eq:fd-L2-bound} then yields \eqref{eq:d-L2-bound}.
	\end{proof}
	
	We now summarize all properties of the planar background solution $(\bar\rho,\bar u,\bar\theta)$ required for the three-dimensional stability analysis. For a fixed small constant $\alpha_0>0$, recall the space-time Gaussian weight $\omega(t,x_1)$:
	\begin{equation}
		y=(1+t)^{-1}\exp\!\left(\frac{\alpha_0x_1^2}{1+t}\right),\qquad
		\omega(t,x_1)=(1+t)^{-1}+\frac{y}{1+y}.
		\label{eq:gaussian-weight}
	\end{equation}
	
	\begin{lemma}
		\label{lem:background-bounds}
		Fix $0<\alpha_0<c_0/2$ sufficiently small. There exists $\delta_0=\delta_0(\alpha_0)>0$ such that for any $0<\delta\le\delta_0(\alpha_0)$, the one-dimensional Navier--Stokes solution $(\bar\rho,\bar u,\bar\theta)$ from Proposition~\ref{thm:background-global} satisfies \eqref{eq:bg-thermal}--\eqref{eq:bg-target} with $\sigma=C\sqrt{\delta}$, together with the weighted estimates
		\begin{align}
			&\sqrt{\omega(t,x_1)}\,|\bar\theta_{x_1}(t,x_1)|\le C\sqrt{\delta}\,(1+t)^{-1},
			\label{eq:weighted-thetax} \\[2pt]
			&\int_0^t (1+s)^{3/2}\int_{\mathbb{R}}
			e^{-\alpha_0x_1^2/(2(1+s))}|d(s,x_1)|^2\,dx_1\,ds
			\le C_{\alpha_0}\delta\log(2+t).
			\label{eq:d-spacetime}
		\end{align}
	\end{lemma}
	\begin{proof}
		By Corollary~\ref{lem:bg-high-order} and Lemma~\ref{lem:profile-construction}, we have
		\[
		\|\partial_{x_1}(\bar\theta-\tilde\theta)(t)\|_\infty\le C\sqrt{\delta}\,(1+t)^{-19/16},\qquad
		|\tilde\theta_{x_1}(t,x_1)|\le C\delta(1+t)^{-1/2}e^{-c_0x_1^2/(1+t)},
		\]
		while \eqref{eq:gaussian-weight} gives
		\[
		\omega(t,x_1)\le 2,\qquad
		\omega(t,x_1)\le 2(1+t)^{-1}e^{\alpha_0 x_1^2/(1+t)}.
		\]
		Combining these bounds yields \eqref{eq:weighted-thetax} whenever $\alpha_0\le 2c_0$. Finally, combining \eqref{eq:d-decomp}, Corollary~\ref{cor:1d-gaussian-B} with $\beta=\alpha_0/2<c_0/4$, and \eqref{eq:fd-L2-bound} gives \eqref{eq:d-spacetime}.
	\end{proof}
	
	\subsection{Properties of the weight \texorpdfstring{$\omega$}{omega} and weighted Sobolev inequalities}
	
	Recalling \eqref{eq:gaussian-weight}, we write
	\begin{equation}
		\xi=\frac{\alpha_0x_1^2}{1+t},\qquad
		\vartheta=\frac{y}{1+y},\qquad
		\omega(t,x_1)=(1+t)^{-1}+\vartheta,\qquad
		y=(1+t)^{-1}e^\xi.
		\label{eq:weight-def}
	\end{equation}
	By direct differentiation of \eqref{eq:weight-def}, the weight $\omega(t,x_1)$ satisfies, for all $(t,x_1)\in[0,\infty)\times\R$,
	\begin{align}
		&(1+t)^{-1}\le\omega(t,x_1)\le 2,\notag \\[2pt]
		&\omega_t(t,x_1)=-(1+t)^{-2}-\frac{1+\xi}{1+t}\vartheta(1-\vartheta)
		\le-(1+t)^{-2}<0,\notag \\[2pt]
		&\partial_{x_1}\omega(t,x_1)=\frac{2\alpha_0x_1}{1+t}\vartheta(1-\vartheta),\notag \\[2pt]
		&\frac{|\partial_{x_1}\omega(t,x_1)|^2}{\omega(t,x_1)}\le 4\alpha_0\bigl(-\omega_t(t,x_1)\bigr),\notag \\[2pt]
		&\bigl|\partial_{x_1}\bigl(\omega(t,x_1)^{-1/2}\bigr)\bigr|\le\sqrt{\alpha_0/e}.
		\label{eq:weight-props}
	\end{align}
	Indeed, since $\omega\ge\vartheta$ and $1-\vartheta\le 1$, the last two inequalities in \eqref{eq:weight-props} follow from
	\begin{align*}
		\frac{|\partial_{x_1}\omega|^2}{\omega}
		&\le\frac{4\alpha_0^2x_1^2}{(1+t)^2}\vartheta(1-\vartheta)^2
		\le 4\alpha_0\frac{\xi}{1+t}\vartheta(1-\vartheta)
		\le 4\alpha_0(-\omega_t),\\[4pt]
		\bigl|\partial_{x_1}\bigl(\omega^{-1/2}\bigr)\bigr|
		&=\frac{|\partial_{x_1}\omega|}{2\omega^{3/2}}
		\le\frac{\sqrt{\alpha_0\xi}}{\sqrt{(1+t)y(1+y)}}
		=\frac{\sqrt{\alpha_0\xi}\,e^{-\xi/2}}{\sqrt{1+y}}
		\le\sqrt{\alpha_0/e}.
	\end{align*}
	
	Using the uniform Lipschitz bound for $\omega^{-1/2}$ and the lower bound $\omega>1/2$ on $\{y>1\}$, we now establish the one-dimensional weighted Sobolev inequality that controls $d$ against the weighted zeroth-order norm and unweighted first-order derivative.
	
	\begin{lemma}
		\label{lem:weighted-multiplier}
		Let $(\bar\rho,\bar u,\bar\theta)$ be the planar background solution from Lemma~\ref{lem:background-bounds}. There exists a nonnegative function $\Lambda(t)$ satisfying
		\begin{equation}
			\int_0^t\Lambda(s)\,ds\le C_{\alpha_0}\delta\log(1+t)\qquad \text{for all }t\ge 0,
			\label{eq:Lambda-int}
		\end{equation}
		such that for every $f\in H^1(\mathbb{R})$ and all $t\ge 0$,
		\begin{equation}
			(1+t)\int_{\mathbb{R}}|d(t,x_1)|^2|f(x_1)|^2\,dx_1
			\le\Lambda(t)\left(\int_{\mathbb{R}}\omega(t,x_1)|f(x_1)|^2\,dx_1
			+\int_{\mathbb{R}}|f'(x_1)|^2\,dx_1\right).
			\label{eq:1d-multiplier}
		\end{equation}
	\end{lemma}
	
	\begin{proof}
		Fix $t\ge 0$ and $x_1\in\R$, and set $\tau=1+t$. By the last bound in \eqref{eq:weight-props}, $\omega(t,\cdot)^{-1/2}$ is Lipschitz on $\R$ with Lipschitz constant at most $\sqrt{\alpha_0/e}$. On the interval $I=[x_1-r,x_1+r]$ with radius
		\[
		r=\frac{\omega(t,x_1)^{-1/2}}{2(1+\sqrt{\alpha_0/e})},
		\]
		we have $\frac{1}{2}\omega(t,x_1)^{-1/2}\le\omega(t,z)^{-1/2}\le\frac{3}{2}\omega(t,x_1)^{-1/2}$ for all $z\in I$. Averaging the fundamental theorem of calculus for $f$ over $z\in I$ yields
		\begin{equation}
			\begin{aligned}
				|f(x_1)|^2
				&\le C\left(r^{-1}\int_{I}|f(z)|^2\,dz
				+r\int_{I}|f'(z)|^2\,dz\right)\\
				&\le C_{\alpha_0}\omega(t,x_1)^{-1/2}\int_{\R}\bigl(\omega(t,z)|f(z)|^2+|f'(z)|^2\bigr)\,dz.
			\end{aligned}
			\label{eq:local-f}
		\end{equation}
		
		In the regions $\{y\le 1\}$ and $\{y>1\}$, we have
		\[
		\omega(t,x_1)^{-1/2}\le\sqrt{2\tau}\,e^{-\alpha_0x_1^2/(2\tau)}\quad (y\le 1),\qquad
		\omega(t,x_1)>\frac{1}{2}\quad (y>1).
		\]
		Applying \eqref{eq:local-f} on $\{y\le 1\}$ and using \eqref{eq:bg-target} for $d$ on $\{y>1\}$ gives
		\begin{equation}
			\begin{aligned}
				&\tau\int_{\R}|d(t,x_1)|^2|f(x_1)|^2\,dx_1\\
				&\quad\le\left(C\sigma^2\tau^{-1}
				+C_{\alpha_0}\tau^{3/2}\int_{\R}e^{-\alpha_0 x_1^2/(2\tau)}|d(t,x_1)|^2\,dx_1\right)\int_{\R}\bigl(\omega|f|^2+|f'|^2\bigr)\,dx_1.
			\end{aligned}
			\label{eq:1d-mult-raw}
		\end{equation}
		Moreover, {\color{black}for all $t\ge0$}, the lower bound $\omega\ge\tau^{-1}$ and \eqref{eq:bg-target} directly bound the left-hand side of \eqref{eq:1d-mult-raw} by the weighted $L^2$ norm of $f$. Hence \eqref{eq:1d-multiplier} holds with
		\[
		\Lambda(t)=\min\left\{C\sigma^2,\,
		C\sigma^2\tau^{-1}+C_{\alpha_0}\tau^{3/2}\int_{\R}e^{-\alpha_0 x_1^2/(2\tau)}|d(t,x_1)|^2\,dx_1\right\}.
		\]
		Using the first bound in the minimum on $[0,\min\{t,1\}]$ and integrating the second bound over $[\min\{t,1\},t]$ via \eqref{eq:d-spacetime} yields \eqref{eq:Lambda-int}.
	\end{proof}

	\section{Three-dimensional perturbation system and low-order energy estimates}
	\label{sec:3d-low-order}
	
	\subsection{Reformulated perturbation systems and energy functionals}
	\label{subsec:3d-framework}
	
	We first record the governing equations and energy-dissipation functionals for the three-dimensional perturbation
	\[
	U=(\phi,\psi,\zeta)=(\rho-\bar\rho,u-\bar u e_1,\theta-\bar\theta)
	\]
	on $\R^3$. Unless otherwise stated, all spatial integrals and $L^p$ norms in Sections~\ref{sec:3d-low-order}--\ref{sec:proof-main} are taken over $\R^3$.
	
	Let $(\bar\rho,\bar u,\bar\theta)(t,x_1)$ be the one-dimensional background solution from Lemma~\ref{lem:background-bounds} with $\sigma=C\sqrt{\delta}$. On any time interval $[0,T]$ on which a smooth solution $(\rho,u,\theta)$ to \eqref{eq:NSF} exists, we make the a priori assumption
	\begin{equation}
		\frac{1}{2}\rho_*\le\rho(t,x)\le 2\rho^*,\qquad
		\frac{1}{2}\theta_*\le\theta(t,x)\le 2\theta^*,\qquad
		\mathcal{M}(t)\le \delta_1,
		\label{eq:apriori-assumption}
	\end{equation}
	where $\delta_1>0$ is a small constant to be determined, and
	\begin{equation}
		\mathcal{M}(t)=\|U(t)\|_{W^{2,\infty}}+\|\psi(t)\|_3
		+\|\nabla U(t)\|_3+\|\nabla^2U(t)\|_3
		+\sqrt{1+t}\,\|\psi(t)\|_\infty\|\nabla U(t)\|_3.
		\label{eq:M-def}
	\end{equation}
	Under \eqref{eq:apriori-assumption}, $U=(\phi,\psi,\zeta)$ and the variables $(q,\psi,h)$ defined in \eqref{eq:qh} are related by a smooth, invertible change of variables with $U=0$ if and only if $(q,\psi,h)=0$.
	
	Denote the material derivative along the velocity field $u$ and the Lam\'e viscosity operator by
	\[
	D_u=\partial_t+u\cdot\nabla,\qquad
	\mathcal{L}\psi=\mu\Delta\psi+(\mu+\lambda)\nabla\diver\psi.
	\]
	Using \eqref{eq:qh}, \eqref{eq:NSF}, and the one-dimensional background equations, direct computation yields the $(q,h,\psi)$ system
	\begin{equation}
		\left\{\begin{aligned}
			&D_uq+\gamma\diver\psi=\mathcal G, \\[2pt]
			&D_uh=\mathcal F, \\[2pt]
			&D_u\psi+R\theta\nabla q=\rho^{-1}\mathcal{L}\psi-\mathcal B,
		\end{aligned}\right.
		\label{eq:qhu-system}
	\end{equation}
	where the source terms $\mathcal G,\mathcal B,\mathcal F$ are given by
	\begin{align}
		&\mathcal G=(\gamma-1)\left(\frac{\kappa\Delta\theta+\mathcal Q(u)}p-\frac{\kappa\bar\theta_{x_1x_1}+\nu\bar u_{x_1}^2}{\bar p}\right)
		-\psi_1\frac{\bar p_{x_1}}{\bar p}, \notag\\[2pt]
		&\mathcal B=\left(\bar u_{x_1}\psi_1+\frac{\phi}{\rho} d+
		\frac{e^q-1}{\rho}\bar p_{x_1}\right)e_1, \label{eq:GBF-def}\\[2pt]
		&\mathcal F=\frac{\kappa\Delta\theta+\mathcal Q(u)
			-e^q(\kappa\bar\theta_{x_1x_1}+\nu\bar u_{x_1}^2)}
		{c_v\rho e^{(\gamma-1)q/\gamma}}
		+\frac{hD_{\bar u}\bar\theta}{\bar\theta}
		+\left[\frac{(\gamma-1)(\bar\theta+h/\gamma)\bar p_{x_1}}{\bar p}
		-\gamma\bar\theta_{x_1}\right]\psi_1.\notag
	\end{align}
	Similarly, subtracting the one-dimensional background equations from \eqref{eq:NSF} in terms of the physical perturbation $U=(\phi,\psi,\zeta)$ yields
	\begin{equation}
		\left\{\begin{aligned}
			&D_u\phi+\rho\diver\psi=-\bar\rho_{x_1}\psi_1-\bar u_{x_1}\phi, \\[2pt]
			&D_u\psi+\frac{R\theta}{\rho}\nabla \phi+R\nabla \zeta
			=\rho^{-1}\mathcal{L}\psi-
			\left[\bar u_{x_1}\psi_1+\frac{\phi}{\rho} d+
			\frac{R}{\rho}(\zeta \bar\rho_{x_1}+\phi\bar\theta_{x_1})\right]e_1, \\[2pt]
			&D_u\zeta+(\gamma-1)\theta\diver\psi=\frac{\kappa}{c_v\rho}\Delta \zeta
			-\bar\theta_{x_1}\psi_1-\frac{\phi}{\rho} D_{\bar u}\bar\theta
			-\frac{p-\bar p}{c_v\rho}\bar u_{x_1}
			+\frac{J_b(\psi)+\mathcal{Q}(\psi)}{c_v\rho},
		\end{aligned}\right.
		\label{eq:physical-pert-system}
	\end{equation}
	where the linear viscous heating term is given by
	\[
	J_b(\psi)=2\bar u_{x_1}(2\mu\partial_{x_1}\psi_1+\lambda\diver \psi).
	\]
	
	Recalling the definitions of $q,h$ in \eqref{eq:qh} and the space-time Gaussian weight $\omega(t,x_1)$ in \eqref{eq:gaussian-weight}, we denote the zeroth-order weighted energy and entropy dissipation functionals by
	\begin{gather*}
		E_0(t)=\|(q,\psi)(t)\|_2^2+\int_{\R^3}\omega(t,x_1) h(t,x)^2\,dx,\\[2pt]
		D_h(t)=\int_{\R^3}\omega(t,x_1)|\nabla h(t,x)|^2\,dx,\qquad
		D_t(t)=\int_{\R^3}\bigl(-\omega_t(t,x_1)\bigr)h(t,x)^2\,dx.
	\end{gather*}
	We define the energy norms $\mathcal{N}_0(t)$ and $\mathcal{N}(t)$ by
	\begin{align}
		&\mathcal{N}_0(t)
		=\|(q,\psi)(t)\|_2^2+\int_{\R^3}\omega(t,x_1) h(t,x)^2\,dx
		+\sum_{j=1}^4\|\nabla^jU(t)\|_2^2,\notag \\[2pt]
		&\mathcal{N}(t)
		=\mathcal{N}_0(t)+(1+t)\left(\|\nabla(q,\psi)(t)\|_2^2
		+\sum_{j=2}^4\|\nabla^jU(t)\|_2^2\right),
		\label{eq:Nstar}
	\end{align}
	and the dissipation functionals $\mathcal{D}_0(t)$, $\mathcal{D}_1(t)$, and $\mathcal{D}(t)$ by
	\begin{align}
		&\mathcal{D}_0(t)
		=\|\nabla(q,\psi)(t)\|_2^2+D_h(t)+D_t(t)
		+\sum_{j=2}^4\|\nabla^j\phi(t)\|_2^2
		+\sum_{j=2}^5\|\nabla^j(\psi,\zeta)(t)\|_2^2,\notag \\[2pt]
		&\mathcal{D}_1(t)
		=\|\nabla^2(q,\psi)(t)\|_2^2
		+\sum_{j=3}^4\|\nabla^j\phi(t)\|_2^2
		+\sum_{j=3}^5\|\nabla^j(\psi,\zeta)(t)\|_2^2,\notag \\[2pt]
		&\mathcal{D}(t)=\mathcal{D}_0(t)+(1+t)\mathcal{D}_1(t).
		\label{eq:dissipation}
	\end{align}
	By \eqref{eq:weight-props}, the entropy dissipation functionals satisfy the global lower bounds
	\[
	D_h(t)\ge(1+t)^{-1}\|\nabla h(t)\|_2^2,\qquad
	D_t(t)\ge(1+t)^{-2}\|h(t)\|_2^2.
	\]
	
	Applying the one-dimensional weighted Sobolev inequality \eqref{eq:1d-multiplier} from Lemma~\ref{lem:weighted-multiplier} to $\phi(t,\cdot,x')$ for each fixed transverse coordinate $x'=(x_2,x_3)\in\R^2$ and integrating over $x'\in\R^2$, we obtain
	\[
	(1+t)\|\phi d\|_2^2
	\le\Lambda(t)\left(\|\sqrt\omega\,\phi\|_2^2
	+\|\partial_{x_1}\phi\|_2^2\right).
	\]
	By \eqref{eq:qh} and the positive state bounds in \eqref{eq:apriori-assumption}, the density perturbation satisfies
	\[
	\phi=\bar\rho\,
	\frac{\bar\theta(e^{q/\gamma}-1)-h/\gamma}{\bar\theta+h/\gamma},
	\qquad |\phi|\le C(|q|+|h|).
	\]
	Since $\omega(t,x_1)\le 2$, the right-hand side is bounded by $C\Lambda(t)\mathcal{N}_0(t)$, and therefore
	\begin{equation}
		(1+t)\|\phi d\|_2^2\le C\Lambda(t)\mathcal N_0(t).
		\label{eq:coupling-L2}
	\end{equation}
	
	\subsection{The Gaussian-weighted relative entropy estimate}
	
	We define the relative entropy function $\mathcal H$ and the weighted entropy functional $E_h(t)$ by
	\[
	\mathcal H=h-\gamma\bar\theta\log\!\left(1+\frac{h}{\gamma\bar\theta}\right),\qquad
	E_h(t)=c_v\int_{\R^3}\omega(t,x_1)\rho \mathcal H\,dx,
	\]
	where $\mathcal{H}\sim h^2$ under \eqref{eq:apriori-assumption}.
	
	\begin{lemma}
		\label{lem:entropy-estimate}
		Under the assumptions of Lemma~\ref{lem:background-bounds} and \eqref{eq:apriori-assumption}, if $\alpha_0>0$ is fixed sufficiently small, there exist constants $c,C>0$ such that
		\begin{equation}
			E_h'(t)+c\bigl(D_h(t)+D_t(t)\bigr)
			\le C\|\nabla q(t)\|_2^2
			+C\bigl(\sigma+\mathcal{M}(t)\bigr)\mathcal{D}_0(t)
			+\frac{C\sigma}{1+t}\mathcal{N}_0(t).
			\label{eq:entropy-ineq}
		\end{equation}
	\end{lemma}
	
	\begin{proof}
		Direct differentiation of $\mathcal{H}(h,\bar\theta)$ gives
		\begin{align*}
			\partial_h\mathcal H
			&=\frac{h}{\gamma\bar\theta+h}
			=\frac{h}{\gamma(\bar\theta+h/\gamma)},\\
			\partial_{\bar\theta}\mathcal H
			&=-\gamma\log\!\left(1+\frac{h}{\gamma\bar\theta}\right)
			+\frac{h}{\bar\theta+h/\gamma}.
		\end{align*}
		Using the continuity equation and the entropy equation in \eqref{eq:qhu-system}, we compute
		\[
		\rho D_u\mathcal H
		=\rho(\partial_h\mathcal H)\mathcal F
		+\rho(\partial_{\bar\theta}\mathcal H)
		(D_{\bar u}\bar\theta+\psi_1\bar\theta_{x_1}).
		\]
		Substituting $\mathcal{F}$ from \eqref{eq:GBF-def} and using \eqref{eq:qh}, we have
		\begin{align*}
			&\frac{\rho h}{\gamma(\bar\theta+h/\gamma)}
			\frac{\kappa\Delta\theta+\mathcal Q(u)-e^q(\kappa\bar\theta_{x_1x_1}+\nu\bar u_{x_1}^2)}
			{c_v\rho e^{(\gamma-1)q/\gamma}}
			=\frac{h}{c_v\gamma\theta}
			\left[\kappa\Delta\theta+\mathcal Q(u)
			-e^q(\kappa\bar\theta_{x_1x_1}+\nu\bar u_{x_1}^2)\right],\\[2pt]
			&(\partial_h\mathcal H)\frac{h}{\bar\theta}
			+\partial_{\bar\theta}\mathcal H
			=\frac{\mathcal H}{\bar\theta}.
		\end{align*}
		Multiplying by $c_v\omega(t,x_1)$ and integrating by parts over $\R^3$ yields the weighted relative entropy identity
		\begin{align}
			&E_h'(t)+c_v\int_{\R^3}\rho(-\omega_t)\mathcal H\,dx\notag \\[2pt]
			&=\int_{\R^3}\frac{\omega h}{\gamma\theta}
			\bigl[\kappa\Delta\theta+\mathcal{Q}(u)
			-e^q(\kappa\bar\theta_{x_1x_1}+\nu \bar u_{x_1}^2)\bigr]\,dx +c_v\int_{\R^3}\omega\rho \mathcal H\frac{D_{\bar u}\bar\theta}{\bar\theta}\,dx
			+c_v\int_{\R^3}\rho \mathcal H u_1\partial_{x_1}\omega\,dx\notag \\[2pt]
			&\qquad +c_v\int_{\R^3}\omega\rho
			\left[-\gamma\log\!\left(1+\frac{h}{\gamma\bar\theta}\right)\bar\theta_{x_1}
			+\frac{\gamma-1}{\gamma} h\frac{\bar p_{x_1}}{\bar p}\right]\psi_1\,dx.
			\label{eq:entropy-identity}
		\end{align}
		
		By the three-dimensional Sobolev embedding and the fourth property in \eqref{eq:weight-props}, we have the weighted Sobolev inequality
		\begin{equation}
			\|\sqrt{\omega}\,h\|_6^2
			\le C\|\nabla(\sqrt{\omega}\,h)\|_2^2
			\le C\bigl(D_h(t)+\alpha_0 D_t(t)\bigr).
			\label{eq:weighted-sobolev}
		\end{equation}
		We now estimate the integrals on the right-hand side of \eqref{eq:entropy-identity} in four steps.
		
		\medskip\noindent
		\emph{Step 1.}
		We first estimate the convective term in \eqref{eq:entropy-identity} by writing $u_1=\bar u+\psi_1$. For the $\bar u$ part, using \eqref{eq:bg-velocity} and $2\sqrt{\xi}\le 1+\xi$ gives
		\begin{equation}
			|\bar u\partial_{x_1}\omega|
			\le C\sigma(1+t)^{-1/2}\frac{2\sqrt{\alpha_0\xi}}{\sqrt{1+t}}\vartheta(1-\vartheta)
			\le C\sigma\sqrt{\alpha_0}\,(-\omega_t),
			\label{eq:bg-transport-pt}
		\end{equation}
		and hence
		\[
		c_v\left|\int_{\R^3}\rho \mathcal H\bar u\partial_{x_1}\omega\,dx\right|
		\le C\sigma\sqrt{\alpha_0}\,D_t(t),
		\]
		which is absorbed by the left-hand side of \eqref{eq:entropy-identity} for small $\sigma>0$. For the perturbation velocity contribution, H\"older's inequality ($L^3\times L^6\times L^2$), \eqref{eq:weight-props}, and \eqref{eq:weighted-sobolev} yield
		\begin{equation}
			\begin{aligned}
				\int_{\R^3}|\psi_1|\,|\partial_{x_1}\omega|h^2\,dx
				&\le\|\psi_1\|_3\|\sqrt\omega\,h\|_6\left(\int_{\R^3}\frac{|\partial_{x_1}\omega|^2h^2}{\omega}\,dx\right)^{1/2}\\
				&\le C\sqrt{\alpha_0}\,\|\psi\|_3\bigl(D_h(t)+D_t(t)\bigr)
				\le C\mathcal{M}(t)\mathcal{D}_0(t).
			\end{aligned}
			\label{eq:pert-transport-bound}
		\end{equation}
		
		\medskip\noindent
		\emph{Step 2.}
		In the heat-conduction term of \eqref{eq:entropy-identity}, writing $\Delta\theta-e^q\bar\theta_{x_1x_1}=\Delta\zeta-(e^q-1)\bar\theta_{x_1x_1}$ and integrating $\Delta\zeta$ by parts gives
		\[
		\int_{\R^3}\frac{\omega h}{\gamma\theta}\kappa\Delta\zeta\,dx
		=-\frac{\kappa}{\gamma}\int_{\R^3}\omega\nabla\!\left(\frac{h}{\theta}\right)\cdot\nabla\zeta\,dx
		-\frac{\kappa}{\gamma}\int_{\R^3}\partial_{x_1}\omega\frac{h}{\theta}\partial_{x_1}\zeta\,dx.
		\]
		Using \eqref{eq:qh}, we expand the gradients $\nabla(h/\theta)$ and $\nabla\zeta$ as
		\begin{align*}
			&\nabla\!\left(\frac{h}{\theta}\right)
			=e^{-\frac{\gamma-1}{\gamma} q}\left(\frac{\bar\theta}{(\bar\theta+h/\gamma)^2}\nabla h-\frac{h}{(\bar\theta+h/\gamma)^2}\nabla\bar\theta-\frac{\gamma-1}{\gamma}\frac{h}{\bar\theta+h/\gamma}\nabla q\right),\\[2pt]
			&\nabla\zeta
			=e^{\frac{\gamma-1}{\gamma} q}\left(\frac{\nabla h}{\gamma}+\frac{\gamma-1}{\gamma} (\bar\theta+h/\gamma)\nabla q\right)+(e^{\frac{\gamma-1}{\gamma} q}-1)\nabla\bar\theta.
		\end{align*}
		The main quadratic term in $\nabla h$ satisfies
		\[
		-\frac{\kappa}{\gamma^2}\int_{\R^3}\omega\frac{\bar\theta}{(\bar\theta+h/\gamma)^2}|\nabla h|^2\,dx
		\le -2\kappa_0 D_h(t),
		\]
		for a constant $\kappa_0>0$. By Cauchy--Schwarz, \eqref{eq:weighted-thetax}, and \eqref{eq:bg-thermal}, the remaining terms in the first integral and the $\bar\theta_{x_1x_1}$ term satisfy
		\begin{align*}
			&\left|-\frac{\kappa(\gamma-1)}{\gamma^2}\int_{\R^3}\omega\frac{\bar\theta}{\bar\theta+h/\gamma}\nabla h\cdot\nabla q\,dx\right|
			\le\frac{1}{4}\kappa_0 D_h(t)+C\|\nabla q\|_2^2,\\[2pt]
			&\int_{\R^3}\omega|h|\bigl(|\nabla h|\,|\nabla q|+|\nabla q|^2\bigr)\,dx
			\le C\|h\|_\infty\bigl(D_h(t)+\|\nabla q\|_2^2\bigr)\le C\mathcal{M}(t)\mathcal{D}_0(t),\\[2pt]
			&\int_{\R^3}\omega|\bar\theta_{x_1}|(|q|+|h|)\bigl(|\nabla q|+|\nabla h|+|\bar\theta_{x_1}|(|q|+|h|)\bigr)\,dx\\*
			&\quad\le C\sigma\bigl(D_h(t)+D_t(t)+\|\nabla q\|_2^2+(1+t)^{-2}\|q\|_2^2\bigr),\\[2pt]
			&\left|\int_{\R^3}\frac{\omega h}{\gamma\theta}\kappa(e^q-1)\bar\theta_{x_1x_1}\,dx\right|
			\le\frac{C\sigma}{1+t}\|\sqrt{\omega}\,h\|_2\|q\|_2
			\le\frac{C\sigma}{1+t}E_0(t).
		\end{align*}
		
		For the integral involving $\partial_{x_1}\omega$, substituting the expansion of $\partial_{x_1}\zeta$ and applying Young's inequality with \eqref{eq:weight-props} and \eqref{eq:weighted-thetax} gives, for any $\eta_0>0$,
		\begin{equation}
			\left|\frac{\kappa}{\gamma}\int_{\R^3}\partial_{x_1}\omega\frac{h}{\theta}\partial_{x_1}\zeta\,dx\right|
			\le\eta_0\bigl(D_h(t)+\|\nabla q\|_2^2\bigr)
			+(C_{\eta_0}\alpha_0+C\sigma\sqrt{\alpha_0})D_t(t)
			+\frac{C\sigma\sqrt{\alpha_0}}{(1+t)^2}\|q\|_2^2.
			\label{eq:omega1-heat-bound}
		\end{equation}
		Fixing $\eta_0=\frac{1}{4}\kappa_0$ and then choosing $\alpha_0>0$ sufficiently small ensures that the $D_t(t)$ term on the right-hand side of \eqref{eq:omega1-heat-bound} is absorbed by the left-hand side of \eqref{eq:entropy-identity}.
		
		\medskip\noindent
		\emph{Step 3.}
		For the last two integrals in \eqref{eq:entropy-identity}, \eqref{eq:weighted-thetax} and \eqref{eq:bg-target} together with Cauchy--Schwarz give
		\[
		c_v\left|\int_{\R^3}\omega\rho \mathcal H\frac{D_{\bar u}\bar\theta}{\bar\theta}\,dx\right|
		+c_v\left|\int_{\R^3}\omega\rho
		\left[-\gamma\log\!\left(1+\frac{h}{\gamma\bar\theta}\right)\bar\theta_{x_1}
		+\frac{\gamma-1}{\gamma} h\frac{\bar p_{x_1}}{\bar p}\right]\psi_1\,dx\right|
		\le\frac{C\sigma}{1+t}E_0(t).
		\]
		
		\medskip\noindent
		\emph{Step 4.}
		Finally, expanding $\mathcal{Q}(u)-e^q\nu\bar u_{x_1}^2$ and using \eqref{eq:bg-velocity}, we bound the viscous dissipation term in \eqref{eq:entropy-identity} by
		\[
		\left|\int_{\R^3}\frac{\omega h}{\gamma\theta}\bigl[\mathcal{Q}(u)-e^q\nu\bar u_{x_1}^2\bigr]\,dx\right|
		\le C\|h\|_\infty\|\nabla \psi\|_2^2
		+C\sigma\bigl(D_t(t)+\|\nabla \psi\|_2^2\bigr)
		+\frac{C\sigma}{(1+t)^2}\|q\|_2^2.
		\]
		Combining the estimates from Steps~1--4 completes the proof of \eqref{eq:entropy-ineq}.
	\end{proof}
	
	\begin{remark}
		\label{rem:weight-role}
		In Lemma~\ref{lem:entropy-estimate}, the weight $\omega(t,x_1)$ plays two roles: it provides the extra decay needed to bound the $\bar\theta_{x_1}$ term in Step~3 via \eqref{eq:weighted-thetax}, and its negative time derivative gives the term $D_t(t)$, which controls the $\partial_{x_1}\omega$ terms in Steps~1 and 2 via \eqref{eq:weight-props} once $\alpha_0>0$ is chosen small.
	\end{remark}
	
	\subsection{Zeroth-order estimates for \texorpdfstring{$(\psi,q)$}{(psi,q)}}
	
	We define the zeroth-order energy $E_p(t)$ for $(\psi,q)$, the cross functional $Y_0(t)$, and the viscous dissipation $\mathcal{D}_u(\psi)$ by
	\begin{gather*}
		E_p(t)=\int_{\R^3}\left(\frac{1}{2}\rho|\psi|^2
		+\frac{\bar p}{\gamma-1}f(q)\right)\,dx,\qquad
		f(q)=e^q+\gamma-1-\gamma e^{q/\gamma},\\[2pt]
		Y_0(t)=\int_{\R^3}\frac{\psi}{\theta}\cdot\nabla h\,dx,\qquad
		\mathcal{D}_u(\psi)=\mu\|\curl \psi\|_2^2+\nu\|\diver \psi\|_2^2,
	\end{gather*}
	where $f(q)\sim q^2$ under \eqref{eq:apriori-assumption} and $\mathcal{D}_u(\psi)\sim\|\nabla\psi\|_2^2$.
	
	\begin{lemma}
		\label{lem:zero-acoustic}
		Under the assumptions of Lemma~\ref{lem:background-bounds} and \eqref{eq:apriori-assumption}, for any given $\eta_0>0$, there exist constants $c,C_{\eta_0}>0$, with $c$ independent of $\eta_0$, such that
		\begin{equation}
			\begin{aligned}
				&\frac{d}{dt}\left(E_p-\frac{\kappa}{\gamma c_p}Y_0\right)
				+c\|\nabla(q,\psi)\|_2^2\\
				&\quad\le C_{\eta_0}\|\nabla^2U\|_2^2+\eta_0 D_h(t)
				+C_{\eta_0}\bigl(\sigma+\mathcal{M}(t)\bigr)\mathcal{D}_0(t)
				+C_{\eta_0}\Gamma(t)\mathcal{N}_0(t),
			\end{aligned}
			\label{eq:zero-acoustic-ineq}
		\end{equation}
		where
		\[
		\Gamma(t)=\frac{\sigma}{1+t}+\sqrt{\frac{\Lambda(t)}{1+t}}+\Lambda(t).
		\]
	\end{lemma}
	
	\begin{proof}
		Differentiating $\frac{\bar p}{\gamma-1}f(q)$ along the pressure equation in \eqref{eq:qhu-system} and using $f'(q)=e^q-e^{q/\gamma}$ together with the background pressure relation
		\[
		D_{\bar u}\bar p=-\gamma\bar p\bar u_{x_1}+(\gamma-1)(\kappa\bar\theta_{x_1x_1}+\nu\bar u_{x_1}^2)
		\]
		gives
		\begin{align*}
			&\partial_t\!\left(\frac{\bar p}{\gamma-1}f(q)\right)+\diver\!\left(\frac{\bar p}{\gamma-1}f(q)u\right) \\[2pt]
			&\quad=-\bar p(e^q-1)\diver\psi
			+\frac{\bar p f'(q)\mathcal G+f(q)
				(D_{\bar u}\bar p+\bar p\bar u_{x_1}+\psi_1\bar p_{x_1})}{\gamma-1}.
		\end{align*}
		Adding this relation to the kinetic energy identity
		\[
		\frac{1}{2}\partial_t(\rho|\psi|^2)+\frac{1}{2}\diver(\rho|\psi|^2 u)=\rho \psi\cdot D_u\psi,
		\]
		and using the momentum equation in \eqref{eq:qhu-system}, we observe that the pressure work cancels against the pressure-gradient integral upon integration by parts:
		\[
		-\int_{\R^3}\bigl[\rho R\theta\nabla q+(e^q-1)\bar p_{x_1}e_1\bigr]\cdot \psi\,dx
		=-\int_{\R^3}\nabla(p-\bar p)\cdot \psi\,dx
		=\int_{\R^3}\bar p(e^q-1)\diver\psi\,dx.
		\]
		Combining the $\bar\theta_{x_1x_1}$ terms in $\bar p f'(q)\mathcal{G}$ and $f(q)D_{\bar u}\bar p$ before integrating over $\R^3$ gives
		\[
		\kappa\int_{\R^3}\left[\bigl(1-e^{-\frac{\gamma-1}{\gamma} q}\bigr)\Delta\zeta
		+g(q)\bar\theta_{x_1x_1}\right]dx,
		\qquad
		g(q)=1-e^{-\frac{\gamma-1}{\gamma} q}+(\gamma-1)(1-e^{q/\gamma}),
		\]
		where $|g(q)|\le Cq^2$. {\color{black}Expanding $\nabla\zeta$ via \eqref{eq:qh} and integrating both the $\Delta\zeta$ and $g(q)\bar\theta_{x_1x_1}$ terms by parts yields}
		\[
		-\kappa\frac{(\gamma-1)^2}{\gamma^2}\int_{\R^3}(\bar\theta+h/\gamma)|\nabla q|^2\,dx
		-\frac{R\kappa}{\gamma c_p}\int_{\R^3}\nabla q\cdot\nabla h\,dx
		+\kappa\frac{\gamma-1}{\gamma}\int_{\R^3}\bigl(e^{q/\gamma}-1\bigr)\nabla q\cdot\nabla\bar\theta\,dx.
		\]
		Hence $E_p(t)$ satisfies
		\begin{align}
			&E_p'(t)+\mathcal{D}_u(\psi)+\kappa\frac{(\gamma-1)^2}{\gamma^2}\int_{\R^3} (\bar\theta+h/\gamma)|\nabla q|^2\,dx\notag \\[2pt]
			&\quad=-\frac{R\kappa}{\gamma c_p}\int_{\R^3}\nabla q\cdot\nabla h\,dx
			+\kappa\frac{\gamma-1}{\gamma}\int_{\R^3}\bigl(e^{q/\gamma}-1\bigr)\nabla q\cdot\nabla\bar\theta\,dx\notag \\[2pt]
			&\qquad +\int_{\R^3}\bigl(1-e^{-\frac{\gamma-1}{\gamma} q}\bigr)\bigl[J_b(\psi)+\mathcal{Q}(\psi)\bigr]\,dx
			+\int_{\R^3}\nu \bar u_{x_1}^2 g(q)\,dx\notag \\[2pt]
			&\qquad -\int_{\R^3} \bar p f(q)\bar u_{x_1}\,dx-\int_{\R^3}\rho \bar u_{x_1}|\psi_1|^2\,dx +\int_{\R^3}\bigl(1-e^{q/\gamma}\bigr)\bar p_{x_1}\psi_1\,dx
			-\int_{\R^3} \phi d\, \psi_1\,dx.
			\label{eq:Ep-identity}
		\end{align}
		
		To cancel the inner product
		\[
		-\frac{R\kappa}{\gamma c_p}\int_{\R^3}\nabla q\cdot\nabla h\,dx
		\]
		on the right-hand side of \eqref{eq:Ep-identity}, we differentiate $Y_0(t)$ in time along \eqref{eq:qhu-system}. Using
		\[
		\partial_t\nabla h=\nabla(D_u h-u\cdot\nabla h)=\nabla \mathcal F-(\nabla u)^{\top}\nabla h-(u\cdot\nabla)\nabla h
		\]
		and integrating by parts over $\R^3$, we find
		\begin{align}
			&Y_0'(t) =-R\int_{\R^3}\nabla q\cdot\nabla h\,dx+\mathcal{I}_0,\notag \\[2pt]
			&\mathcal{I}_0 =\int_{\R^3}\frac{\mathcal{L}\psi}{\rho\theta}\cdot\nabla h\,dx
			-\int_{\R^3}\frac{\mathcal B}{\theta}\cdot\nabla h\,dx
			-\int_{\R^3}\diver\left(\frac{\psi}{\theta}\right)\mathcal F\,dx\notag \\[2pt]
			&\quad +\int_{\R^3}\left(\diver u-\frac{D_u\theta}{\theta}\right)
			\frac{\psi}{\theta}\cdot\nabla h\,dx
			-\int_{\R^3}\frac{\psi}{\theta}\cdot(\nabla u)^{\top}\nabla h\,dx.
			\label{eq:Y0-identity}
		\end{align}
		Subtracting $\frac{\kappa}{\gamma c_p}Y_0'(t)$ from \eqref{eq:Ep-identity} cancels the $\nabla q\cdot\nabla h$ integral.
		
		We now estimate the remaining integrals in \eqref{eq:Ep-identity} and the five integrals $I_{0,1},\dots,I_{0,5}$ in $\mathcal{I}_0$. In \eqref{eq:Ep-identity}, using \eqref{eq:bg-thermal}--\eqref{eq:bg-target} and Cauchy--Schwarz gives
		\[
		\left|\int_{\R^3}\bigl(e^{q/\gamma}-1\bigr)\nabla q\cdot\nabla\bar\theta\,dx\right|
		\le C\sigma\|\nabla q\|_2^2+\frac{C\sigma}{1+t}\|q\|_2^2,
		\]
		while \eqref{eq:coupling-L2} bounds the coupling with $d$ by
		\begin{equation}
			\left|\int_{\R^3}\phi d\,\psi_1\,dx\right|
			\le\|\phi d\|_2\|\psi_1\|_2
			\le C\sqrt{\frac{\Lambda(t)}{1+t}}\,\mathcal N_0(t).
			\label{eq:coupling-zero}
		\end{equation}
		All other integrals on the right-hand side of \eqref{eq:Ep-identity} are bounded by
		\[
		C\bigl(\sigma+\mathcal{M}(t)\bigr)\mathcal{D}_0(t)+\frac{C\sigma}{1+t}\mathcal{N}_0(t).
		\]
		
		To bound $I_{0,1},\dots,I_{0,5}$ in $\mathcal{I}_0$, set $\tau=1+t$. By \eqref{eq:qh} and \eqref{eq:bg-thermal}, the chain rule gives, for $1\le\ell\le 4$,
		\begin{equation}
			\|\nabla^\ell(q,h)\|_2
			\le C\|\nabla^\ell U\|_2
			+C\sigma\sum_{i=0}^{\ell-1}\tau^{-(\ell-i)/2}\|\nabla^iU\|_2,
			\label{eq:qh-chain}
		\end{equation}
		while Gagliardo--Nirenberg interpolation yields the homogeneous product bound
		\begin{equation}
			\left\|\prod_{i=1}^n\nabla^{r_i}U\right\|_p
			\le C\|U\|_\infty^{n-1}\|\nabla^rU\|_p,
			\qquad r_i\ge1,\quad \sum_{i=1}^n r_i=r,\quad p\in\{2,6\}.
			\label{eq:homogeneous-product}
		\end{equation}
		Combining \eqref{eq:qh-chain} with the lower bounds for $D_h(t)$ and $D_t(t)$ gives
		\begin{equation}
			\begin{aligned}
				&\tau^{-1}\|\nabla h\|_2^2\le D_h(t),
				\qquad \tau^{-2}\|U\|_2^2\le C\bigl(D_t(t)+\tau^{-2}\|(q,\psi)\|_2^2\bigr),\\[2pt]
				&\tau^{-1}\|\nabla U\|_2^2\le C\bigl(\|\nabla(q,\psi)\|_2^2+D_h(t)+\sigma^2\tau^{-2}\|U\|_2^2\bigr),\\[2pt]
				&\|\nabla h\|_6^2+\|\nabla^2(q,h)\|_2^2
				\le C\|\nabla^2U\|_2^2+C\sigma^2\mathcal D_0(t)
				+C\sigma^2\tau^{-1}\mathcal N_0(t).
			\end{aligned}
			\label{eq:zero-order-controls}
		\end{equation}
		
		Using \eqref{eq:zero-order-controls}, we estimate $I_{0,1},\dots,I_{0,5}$ in \eqref{eq:Y0-identity} as follows:
		\begin{enumerate}[(i)]
			\item In $I_{0,1}$, integrating one spatial derivative in $\mathcal{L}\psi$ by parts and applying Young's inequality together with \eqref{eq:zero-order-controls} gives, for any $a>0$,
			\[
			|I_{0,1}|\le a\|\nabla(q,\psi)\|_2^2+C_a \|\nabla^2U\|_2^2
			+C_a\bigl(\sigma+\mathcal{M}(t)\bigr)\mathcal D_0(t)
			+C_a\sigma\tau^{-1}\mathcal N_0(t).
			\]
			\item In $I_{0,2}$, the term involving $\phi d$ inside $\mathcal{B}$ is bounded using Cauchy--Schwarz, $\omega^{-1}\le 1+t$, and \eqref{eq:coupling-L2}: for any $b>0$,
			\begin{equation}
				\begin{aligned}
					\left|\int_{\R^3}\frac{1}{\rho\theta}\phi d\,\partial_{x_1}h\,dx\right|
					&\le b D_h(t)+C_b\|\omega^{-1/2}\phi d\|_2^2\\
					&\le b D_h(t)+C_b\Lambda(t)\mathcal N_0(t).
				\end{aligned}
				\label{eq:coupling-cross}
			\end{equation}
			Using \eqref{eq:bg-velocity}--\eqref{eq:bg-target} for the remaining terms in $\mathcal{B}$, we have, for any $b>0$,
			\[
			|I_{0,2}|\le bD_h(t)
			+C_b\bigl(\Lambda(t)+\sigma^2\tau^{-1}\bigr)\mathcal N_0(t).
			\]
			\item In $I_{0,3}$, we decompose $\mathcal F=\mathcal F_\kappa+\mathcal F_\nu+\mathcal F_b$ according to \eqref{eq:GBF-def}:
			\begin{align*}
				\mathcal F_\kappa&=\frac{\kappa e^{-\frac{\gamma-1}{\gamma} q}}{c_v\rho}
				\left[\Delta\zeta-(e^q-1)\bar\theta_{x_1x_1}\right],\\
				\mathcal F_\nu&=\frac{e^{-\frac{\gamma-1}{\gamma} q}}{c_v\rho}
				\left[J_b(\psi)+\mathcal Q(\psi)+(1-e^q)\nu\bar u_{x_1}^2\right],\\
				\mathcal F_b&=\frac{D_{\bar u}\bar\theta}{\bar\theta}h
				-\left(\gamma\bar\theta_{x_1}-(\gamma-1)\bar\theta
				\frac{\bar p_{x_1}}{\bar p}\right)\psi_1
				+\frac{\gamma-1}{\gamma} h\frac{\bar p_{x_1}}{\bar p}\psi_1.
			\end{align*}
			Taking $L^2(\R^3)$ norms gives
			\begin{align*}
				&\|\mathcal F_\kappa\|_2 \le C\bigl(\|\nabla^2U\|_2+\sigma\tau^{-1}\|(q,\psi)\|_2\bigr),\\
				&\|\mathcal F_\nu\|_2 \le C\Bigl[\bigl(\mathcal{M}(t)+\sigma\tau^{-1}\bigr)\|\nabla(q,\psi)\|_2
				+\sigma^2\tau^{-2}\|(q,\psi)\|_2\Bigr],\\
				&\|\mathcal F_b\|_2 \le C\sigma\bigl(\tau^{-1}\|h\|_2
				+\tau^{-1/2}\|(q,\psi)\|_2\bigr),\\
				&\|\diver(\psi/\theta)\|_2 \le
				C\bigl(\|\nabla(q,\psi)\|_2+\mathcal{M}(t)\|\nabla^2U\|_2+\sigma\tau^{-1/2}\|(q,\psi)\|_2\bigr).
			\end{align*}
			Using the entropy dissipation bounds
			\begin{align*}
				\tau^{-1}\|\nabla(q,\psi)\|_2\|h\|_2&\le C\bigl(\|\nabla(q,\psi)\|_2^2+D_t(t)\bigr),\\[2pt]
				\tau^{-3/2}\|(q,\psi)\|_2\|h\|_2&\le C\bigl(\tau^{-1}\|(q,\psi)\|_2^2+D_t(t)\bigr),
			\end{align*}
			Cauchy--Schwarz and Young's inequality yield, for any $a>0$,
			\[
			|I_{0,3}|\le a\|\nabla(q,\psi)\|_2^2+C_a\|\nabla^2U\|_2^2
			+C_a\bigl(\sigma+\mathcal{M}(t)\bigr)\mathcal D_0(t)
			+C_a\sigma\tau^{-1}\mathcal N_0(t).
			\]
			\item In $I_{0,4}+I_{0,5}$, using the temperature equation in \eqref{eq:NSF} to write
			\[
			\diver u-\frac{D_u\theta}{\theta}
			=\gamma\diver u-\frac{\kappa\Delta\theta+\mathcal Q(u)}{c_v\rho\theta},
			\]
			and applying H\"older's inequality ($L^3\times L^2\times L^6$) together with \eqref{eq:zero-order-controls}, we obtain, for any $b>0$,
			\begin{align*}
				|I_{0,4}|+|I_{0,5}|
				&\le C\mathcal{M}(t)\bigl(\|\nabla(q,\psi)\|_2+\|\nabla^2U\|_2\bigr)\|\nabla h\|_6
				+C\sigma\tau^{-1}\|(q,\psi)\|_2\|\nabla h\|_2\\
				&\le bD_h(t)+C_b\bigl(\sigma+\mathcal{M}(t)\bigr)\mathcal D_0(t)+C_b\sigma\tau^{-1}\mathcal N_0(t).
			\end{align*}
		\end{enumerate}
		Choosing $a>0$ small relative to the left-hand side of \eqref{eq:Ep-identity} and choosing $b>0$ so that $\frac{2\kappa}{\gamma c_p}b\le\eta_0$ yields \eqref{eq:zero-acoustic-ineq}.
	\end{proof}
	
	\subsection{First-order divergence and curl estimates}
	
	At the first-derivative level, we define $E_d(t)$, $E_c(t)$, and the cross functional $Y_1(t)$ by
	\begin{align*}
		E_d(t)&=\frac{1}{2}\left(\|\diver\psi\|_2^2+\frac{R}{\gamma}\|\sqrt\theta\,\nabla q\|_2^2\right),\qquad
		E_c(t)=\frac{1}{2}\|\curl(\rho\psi)\|_2^2,\\[2pt]
		Y_1(t)&=R\int_{\R^3}\frac{e^{-q/\gamma}}{\bar p}(\diver\psi)\Delta h\,dx.
	\end{align*}
	
	\begin{lemma}
		\label{lem:first-order-div-curl}
		Under the assumptions of Lemma~\ref{lem:background-bounds} and \eqref{eq:apriori-assumption}, there exist constants $c,C,C_1>0$ such that
		\begin{equation}
			\begin{aligned}
				&(1+t)\frac{d}{dt}\left(E_d+E_c-\frac{\kappa}{\gamma c_p}Y_1\right)
				+c(1+t)\|\nabla^2(q,\psi)\|_2^2 \\[2pt]
				&\quad\le C_1(1+t)\|\nabla^3 U\|_2^2
				+C\bigl(\sigma+\mathcal{M}(t)\bigr)\mathcal{D}(t)
				+C\Gamma(t)\mathcal{N}(t).
			\end{aligned}
			\label{eq:div-curl-ineq}
		\end{equation}
	\end{lemma}
	
	\begin{proof}
		Taking the divergence of the momentum equation in \eqref{eq:qhu-system} gives
		\[
		D_u(\diver\psi)+\diver(R\theta\nabla q)=N_\psi,\qquad\text{where}\quad
		N_\psi=\diver(\rho^{-1}\mathcal{L}\psi-\mathcal B)-\partial_i u_j\partial_j \psi_i.
		\]
		Differentiating $E_d(t)$ in time along \eqref{eq:qhu-system} and integrating by parts, we use
		\[
		-\int_{\R^3} (\diver\psi)\diver(R\theta\nabla q)\,dx-\int_{\R^3} R\theta\nabla q\cdot\nabla (\diver\psi)\,dx=0
		\]
		to obtain
		\begin{align}
			E_d'(t)&=\int_{\R^3} (\diver\psi)N_\psi\,dx+\frac{1}{\gamma}\int_{\R^3} R\theta\nabla q\cdot\nabla \mathcal G\,dx
			+\frac{1}{2}\int_{\R^3}(\diver u)(\diver\psi)^2\,dx\notag \\[2pt]
			&\quad +\frac{1}{2\gamma}\int_{\R^3}\bigl(D_u(R\theta)+R\theta\diver u\bigr)|\nabla q|^2\,dx
			-\frac{1}{\gamma}\int_{\R^3} R\theta\nabla q\cdot(\nabla u)^{\top}\nabla q\,dx,
			\label{eq:Ed-identity}
		\end{align}
		where
		\[
		D_u(R\theta)+R\theta\diver u=(2-\gamma)R\theta\diver u+\frac{R}{c_v\rho}\bigl(\kappa\Delta\theta+\mathcal Q(u)\bigr).
		\]
		
		In the first integral of \eqref{eq:Ed-identity}, writing $\varpi=\curl\psi$ and integrating by parts via $\mathcal{L}\psi=\nu\nabla (\diver\psi)-\mu\curl\varpi$ gives
		\begin{equation}
			\int_{\R^3} (\diver\psi)\diver(\rho^{-1}\mathcal{L}\psi)\,dx
			=-\nu\int_{\R^3}\rho^{-1}|\nabla (\diver\psi)|^2\,dx
			-\mu\int_{\R^3}\varpi\cdot\bigl(\nabla (\diver\psi)\times\nabla\rho^{-1}\bigr)\,dx.
			\label{eq:w-viscous}
		\end{equation}
		In the second integral of \eqref{eq:Ed-identity}, integrating by parts once and writing the heat-conduction parts of $\mathcal G$ and $\mathcal F$ as
		\begin{align*}
			\mathcal G_\kappa&=\frac{(\gamma-1)\kappa}{\bar p}e^{-q/\gamma}\mathcal{T}, \qquad
			\mathcal F_\kappa=\frac{(\gamma-1)\kappa}{R\rho}\mathcal{T},\\[2pt]
			\mathcal{T}&=\frac{\Delta h}{\gamma}+\frac{\gamma-1}{\gamma} (\bar\theta+h/\gamma)\Delta q
			+2\frac{\gamma-1}{\gamma}\nabla (\bar\theta+h/\gamma)\cdot\nabla q\\[2pt]
			&\quad +\frac{(\gamma-1)^2}{\gamma^2}(\bar\theta+h/\gamma)|\nabla q|^2
			+\bar\theta_{x_1x_1}\bigl(1-e^{q/\gamma}\bigr),
		\end{align*}
		expanding $\diver(R\theta\nabla q)$ shows that the $\Delta q$ term in $\mathcal{T}$ gives
		\begin{equation}
			-\frac{(\gamma-1)^2\kappa}{\gamma^2}\int_{\R^3}\frac{\theta}{\rho}|\Delta q|^2\,dx
			\label{eq:pressure-hessian}
		\end{equation}
		together with the lower-order term
		\[
		-\frac{(\gamma-1)^2\kappa}{\gamma^2R}\int_{\R^3}\rho^{-1}\nabla(R\theta)\cdot\nabla q\,\Delta q\,dx,
		\]
		while the $\Delta h/\gamma$ term in $\mathcal{T}$ gives
		\[
		-\frac{\kappa}{\gamma c_p}\int_{\R^3} \frac{R e^{-q/\gamma}}{\bar p}\diver(R\theta\nabla q)\Delta h\,dx.
		\]
		To cancel this $\Delta h$ integral, we differentiate $Y_1(t)$ in time along \eqref{eq:qhu-system}:
		\begin{align}
			&Y_1'(t) =-R\int_{\R^3}\frac{e^{-q/\gamma}}{\bar p}\diver(R\theta\nabla q)\Delta h\,dx+\mathcal I_1,\notag\\[2pt]
			&\mathcal I_1 =R\int_{\R^3}\frac{e^{-q/\gamma}}{\bar p}
			\bigl[N_\psi\Delta h+(\diver\psi)(\Delta \mathcal F-2\partial_i u_j\partial_{ij}h-\Delta u\cdot\nabla h)\bigr]dx\notag\\[2pt]
			&\quad+R\int_{\R^3}\frac{e^{-q/\gamma}}{\bar p}
			\bigl[\diver u-D_u(q/\gamma+\log\bar p)\bigr](\diver\psi)\Delta h\,dx.
			\label{eq:Y1-identity}
		\end{align}
		Subtracting $\frac{\kappa}{\gamma c_p}Y_1'(t)$ from $E_d'(t)$ cancels the mixed integral involving $\Delta h$.
		
		Next, for $E_c(t)$, we set $K_b=(\rho\bar u_{x_1}\psi_1+\phi d)e_1$. Writing the momentum equation in conservative form and taking its curl eliminates $\nabla(p-\bar p)$:
		\begin{align}
			&(\rho\psi)_t+\partial_{x_1}\bigl(\bar u (\rho\psi)\bigr)+\diver\bigl((\rho\psi)\otimes \psi\bigr)+\nabla(p-\bar p)
			=\mathcal{L}\psi-K_b,\notag \\[2pt]
			&D_{\bar u}\curl(\rho\psi)-\mu\Delta\varpi
			=-\curl\diver\bigl((\rho\psi)\otimes \psi\bigr)
			-\curl K_b-\bar u_{x_1}\curl(\rho\psi)
			-e_1\times\partial_{x_1}(\bar u_{x_1}\rho\psi).
			\label{eq:curl-eq}
		\end{align}
		Taking the $L^2(\R^3)$ inner product of \eqref{eq:curl-eq} with $\curl(\rho\psi)$ and integrating by parts yields
		\begin{align}
			E_c'(t)+\mu\int_{\R^3}\rho|\nabla\varpi|^2\,dx
			&=-\mu\int_{\R^3}\bigl(\nabla \curl(\rho\psi)-\rho\nabla\varpi\bigr):\nabla\varpi\,dx\notag \\[2pt]
			&\quad -\int_{\R^3}\Bigl[\curl\curl(\rho\psi)\cdot\bigl(\diver(\rho\psi\otimes\psi)+K_b\bigr)
			+\frac{1}{2}\bar u_{x_1}|\curl(\rho\psi)|^2\Bigr]dx\notag \\[2pt]
			&\quad -\int_{\R^3}\curl(\rho\psi)\cdot\bigl[e_1\times\partial_{x_1}(\bar u_{x_1}\rho\psi)\bigr]dx.
			\label{eq:Ec-identity}
		\end{align}
		Since
		\[
		\|\nabla^2 \psi\|_2^2{\color{black}\approx}\|\nabla (\diver\psi)\|_2^2+\|\nabla\varpi\|_2^2,\qquad
		\|\nabla^2 q\|_2^2{\color{black}\approx}\|\Delta q\|_2^2,
		\]
		the main dissipative terms in $(1+t)\bigl(E_d'+E_c'-\frac{\kappa}{\gamma c_p}Y_1'\bigr)$ bound $2c(1+t)\|\nabla^2(q,\psi)\|_2^2$.
		
		We now estimate the remaining terms in \eqref{eq:Ed-identity}, \eqref{eq:Y1-identity}, and \eqref{eq:Ec-identity}. Set $\tau=1+t$. By \eqref{eq:qh-chain} and the definitions of $\mathcal{D}_0(t)$ and $\mathcal{D}(t)$, we have
		\begin{equation}
			\begin{aligned}
				&{\tau}^{-1}\|\nabla h\|_2^2+{\tau}^{-2}\|h\|_2^2\le\mathcal D_0(t),\\
				&\|\nabla^2(q,h)\|_2^2
				\le C\mathcal D_0(t)+C\sigma^2{\tau}^{-1}\mathcal N_0(t),\\
				&{\tau}\|\nabla^3(q,h)\|_2^2
				\le C{\tau}\|\nabla^3U\|_2^2+C\sigma^2\mathcal D(t)
				+C\sigma^2{\tau}^{-1}\mathcal N_0(t),\\
				&|\nabla\rho^{-1}|+|\nabla\theta|+\bigl|\nabla(e^{-q/\gamma}/\bar p)\bigr|
				\le C|\nabla U|+C\sigma\tau^{-1/2}.
			\end{aligned}
			\label{eq:first-order-chain}
		\end{equation}
		
		\medskip\noindent
		\emph{Step 1.}
		We first bound the remaining terms in \eqref{eq:Ed-identity}. After extracting \eqref{eq:pressure-hessian} and cancelling the $\Delta h$ term via $Y_1'(t)$, the nonlinear terms in \eqref{eq:Ed-identity} {\color{black}include the following types}
		\[
		\nabla\psi\,|\nabla q|^2,\quad (\diver\psi)^3,\quad
		\Delta\zeta\,|\nabla q|^2,\quad
		\nabla(q,h)\,\nabla q\,\Delta q,\quad
		|\nabla\psi|^2\Delta q,\quad
		|\nabla\psi|^2|\nabla q|^2.
		\]
		For the terms of the first three types in \eqref{eq:Ed-identity}, integrating by parts once to move a derivative from $\nabla\psi$ or $\Delta\zeta$ onto the remaining factors and applying H\"older's inequality ($L^3\times L^6\times L^2$) gives
		\begin{align*}
			&\tau\left|\int_{\R^3}\left(\frac{2-\gamma}{2\gamma}R\theta(\diver\psi)|\nabla q|^2-\frac{1}{\gamma}R\theta\nabla q\cdot(\nabla\psi)^{\top}\nabla q+\frac{1}{2}(\diver\psi)^3+\frac{R\kappa}{2\gamma c_v\rho}\Delta\zeta\,|\nabla q|^2\right)dx\right|\\
			&\quad\le C\tau\bigl(\|\psi\|_3+\|\nabla\zeta\|_3\bigr)\|\nabla q\|_6\|\nabla^2 q\|_2
			+C\tau\|\psi\|_3\|\diver\psi\|_6\|\nabla\diver\psi\|_2\\
			&\qquad +C\sigma\sqrt{\tau}\,\|\nabla(q,\psi)\|_2\|\nabla^2(q,\psi)\|_2
			\le C\mathcal{M}(t)\tau\|\nabla^2(q,\psi)\|_2^2+C\sigma\mathcal{D}(t).
		\end{align*}
		{\color{black}For the velocity product in $N_\psi$, integration by parts gives
			\begin{align*}
				&\int_{\R^3}(\diver\psi)\partial_i\psi_j\partial_j\psi_i\,dx=-\int_{\R^3}\psi_j\partial_i(\diver\psi)\partial_j\psi_i\,dx
				-\int_{\R^3}(\diver\psi)\psi\cdot\nabla(\diver\psi)\,dx.
			\end{align*}
			H\"older's inequality with exponents $(3,6,2)$ therefore yields
			\begin{align*}
				\tau\left|\int_{\R^3}(\diver\psi)\partial_i\psi_j\partial_j\psi_i\,dx\right|
				\le C\tau\|\psi\|_3\|\nabla\psi\|_6\|\nabla(\diver\psi)\|_2 \le C\mathcal M(t)\tau\|\nabla^2\psi\|_2^2.
			\end{align*}
			The coefficient derivatives are bounded using $(3,3,6,6)$ for the perturbations and $L^\infty$ for background gradients:
			\begin{align*}
				&C\tau(\|\psi\|_3+\|\nabla\zeta\|_3)\|\nabla U\|_3\|\nabla q\|_6^2+C\sigma\sqrt\tau(\|\psi\|_3+\|\nabla\zeta\|_3)\|\nabla q\|_2\|\nabla q\|_6\\
				&\le C\mathcal M(t)^2\tau\|\nabla^2q\|_2^2
				+C\sigma\mathcal M(t)\mathcal D(t).
		\end{align*}}
		{\color{black}By \eqref{eq:qh} and \eqref{eq:bg-thermal},
			$|\nabla(q,h)|\le C|\nabla U|+C\sigma\tau^{-1/2}|U|$. Hence
			\begin{align*}
				\tau\int_{\R^3}|\nabla(q,h)|\,|\nabla q|\,|\Delta q|\,dx
				&\le C\tau\|\nabla U\|_3\|\nabla q\|_6\|\Delta q\|_2
				+C\sigma\sqrt\tau\|U\|_\infty\|\nabla q\|_2\|\Delta q\|_2\\
				&\le C\mathcal M(t)\tau\|\nabla^2q\|_2^2+C\sigma\mathcal M(t)\mathcal D(t).
			\end{align*}
			For $|\nabla\psi|^2|\Delta q|$ and $|\nabla\psi|^2|\nabla q|^2$, H\"older's inequality with exponents $(3,6,2)$ and $(3,3,6,6)$ gives the bound
			$C(\mathcal M(t)+\mathcal M(t)^2)\tau\|\nabla^2(q,\psi)\|_2^2$.
			Products with two coefficient gradients are estimated by the same decomposition. In particular,
			\[
			\tau\int_{\R^3}|\nabla h|^2|\nabla q|^2\,dx
			\le C\tau\|\nabla U\|_3^2\|\nabla q\|_6^2
			+C\sigma^2\|U\|_\infty^2\|\nabla q\|_2^2 .
			\]}
		The linear background terms in \eqref{eq:Ed-identity} and $\mathcal{T}$ are bounded by
		\[
		C\sigma\bigl(\tau\|\nabla^2(q,\psi)\|_2^2+\tau^{-1}\|(q,\psi)\|_2^2\bigr)
		\le C\sigma\mathcal{D}(t)+C\sigma\tau^{-1}\mathcal{N}_0(t),
		\]
		and the viscosity commutator in \eqref{eq:w-viscous} satisfies
		\begin{align*}
			\tau\left|\int_{\R^3}\varpi\cdot\bigl(\nabla(\diver\psi)\times\nabla\rho^{-1}\bigr)\,dx\right|
			&\le C\tau\|\nabla\phi\|_3\|\varpi\|_6\|\nabla(\diver\psi)\|_2
			+C\sigma\sqrt{\tau}\,\|\varpi\|_2\|\nabla(\diver\psi)\|_2\\
			&\le C\bigl(\mathcal{M}(t)+\sigma\bigr)\mathcal{D}(t).
		\end{align*}
		
		\medskip\noindent
		\emph{Step 2.}
		In \eqref{eq:Ed-identity}, \eqref{eq:Y1-identity}, and \eqref{eq:Ec-identity}, the vectors $\mathcal{B}$ and $K_b$ appear in the three integrals
		\[
		-\int_{\R^3}(\diver\psi)\diver\mathcal B\,dx,\qquad
		-\int_{\R^3}\frac{R e^{-q/\gamma}}{\bar p}(\diver\mathcal B)\Delta h\,dx,\qquad
		-\int_{\R^3}\curl(\rho\psi)\cdot\curl K_b\,dx.
		\]
		Integrating by parts once in each integral transfers the divergence or curl onto the multiplier. When the gradient in the second integral hits the coefficient $R e^{-q/\gamma}/\bar p$, we use \eqref{eq:bg-target} and \eqref{eq:first-order-chain} to bound
		\begin{align*}
			\tau\left|\int_{\R^3}\phi d\,\nabla\!\left(\frac{R e^{-q/\gamma}}{\bar p}\right)\cdot e_1\,\Delta h\,dx\right|
			&\le C\sigma\|\phi\|_\infty\|\nabla q\|_2\|\Delta h\|_2
			+C\sigma^2{\tau}^{-1}\|\phi\|_2\|\Delta h\|_2\\
			&\le C\bigl(\sigma+\mathcal{M}(t)\bigr)\mathcal{D}_0(t)+C\sigma\tau^{-1}\mathcal{N}_0(t).
		\end{align*}
		When the derivative hits $\diver\psi$, $\Delta h$, or $\curl(\rho\psi)$, Cauchy--Schwarz and \eqref{eq:coupling-L2} give, for any $\eta_0>0$,
		\begin{align*}
			&\tau\int_{\R^3}|\phi d|\,\bigl(|\nabla\diver\psi|+|\nabla\Delta h|+|\curl\curl(\rho\psi)|\bigr)\,dx\\
			&\quad\le\eta_0\tau\|\nabla^2(q,\psi)\|_2^2+C_{\eta_0}\tau\|\nabla^3 U\|_2^2
			+C_{\eta_0}\tau\|\nabla\curl(\rho\psi)-\rho\nabla\varpi\|_2^2\\
			&\qquad +C_{\eta_0}\bigl(\sigma+\mathcal{M}(t)\bigr)\mathcal{D}(t)+C_{\eta_0}\Lambda(t)\mathcal{N}_0(t).
		\end{align*}
		{\color{black}The remainder $C\sigma^2\tau^{-1}\mathcal N_0(t)$ from \eqref{eq:first-order-chain} is included in $C\Lambda(t)\mathcal N_0(t)$, since both entries in the minimum defining $\Lambda$ are bounded below by $c\sigma^2\tau^{-1}$.}
		The remaining parts of $\mathcal{B}$ and $K_b$ satisfy the same bound with $\sigma\tau^{-1}$ in place of $\Lambda(t)$.
		
		\medskip\noindent
		\emph{Step 3.}
		For the remaining terms in \eqref{eq:Ec-identity}, expanding $\curl(\rho\psi)$ and $\diver(\rho\psi\otimes\psi)$ and applying H\"older's and Sobolev's inequalities gives
		\begin{align*}
			&\|\nabla\curl(\rho\psi)-\rho\nabla\varpi\|_2
			\le C\bigl(\|\nabla\phi\|_3\|\nabla\psi\|_6+\|\psi\|_3\|\nabla^2\phi\|_6
			+\sigma\tau^{-1/2}\|\nabla\psi\|_2+\sigma\tau^{-1}\|\psi\|_2\bigr)\\
			&\quad\le C\mathcal M(t)\bigl(\|\nabla^2\psi\|_2+\|\nabla^3\phi\|_2\bigr)
			+C\sigma\bigl({\tau}^{-1/2}\|\nabla\psi\|_2+{\tau}^{-1}\|\psi\|_2\bigr),\\[2pt]
			&\|\diver(\rho\psi\otimes\psi)\|_2
			\le C\|\psi\|_3\|\nabla\psi\|_6
			+C\bigl(\|\psi\|_\infty\|\nabla\phi\|_3
			+\sigma {\tau}^{-1/2}\|\psi\|_3\bigr)\|\nabla\psi\|_2\\
			&\quad\le C\mathcal M(t)\|\nabla^2\psi\|_2
			+C\bigl(\|\psi\|_\infty\|\nabla\phi\|_3
			+\sigma {\tau}^{-1/2}\|\psi\|_3\bigr)\|\nabla\psi\|_2.
		\end{align*}
		By the definition of $\mathcal{M}(t)$ in \eqref{eq:M-def}, we have
		\[
		\tau\|\psi\|_\infty^2\|\nabla\phi\|_3^2\|\nabla\psi\|_2^2
		\le\mathcal{M}(t)^2\|\nabla\psi\|_2^2
		\le\mathcal{M}(t)^2\mathcal{D}_0(t).
		\]
		Multiplying \eqref{eq:Ec-identity} by $\tau=1+t$ and using \eqref{eq:bg-velocity} for the last two integrals in \eqref{eq:Ec-identity}, we bound the right-hand side of $\tau E_c'(t)$ by
		\[
		\eta_0\tau\|\nabla^2(q,\psi)\|_2^2+C_{\eta_0}\bigl(\sigma+\mathcal{M}(t)\bigr)\mathcal{D}(t)+C_{\eta_0}\Gamma(t)\mathcal{N}(t).
		\]
		
		\medskip\noindent
		\emph{Step 4.}
		Finally, we estimate the terms in $\mathcal{I}_1$ from \eqref{eq:Y1-identity}:
		\begin{enumerate}[(i)]
			\item In the first and third terms of $\mathcal{I}_1$ (with $\diver\mathcal{B}$ already bounded in Step~2), {\color{black}integration by parts in the viscous term gives the bound
				\[
				C\tau\|\nabla^2\psi\|_2\bigl(\|\nabla^3 h\|_2+
				\|\nabla U\|_3\|\Delta h\|_6+\sigma\tau^{-1/2}\|\Delta h\|_2\bigr)
				\]
				by \eqref{eq:first-order-chain}. Integrating one derivative off $\psi$ in the quadratic velocity products gives}
			\begin{align*}
				&\tau\left|\int_{\R^3}\frac{R e^{-q/\gamma}}{\bar p}\bigl(\partial_i\psi_j\partial_j\psi_i\Delta h+2(\diver\psi)\partial_i\psi_j\partial_{ij}h\bigr)\,dx\right|\\
				&\quad\le C\tau\|\psi\|_3\bigl(\|\nabla^2\psi\|_2\|\nabla^2 h\|_6+\|\nabla\psi\|_6\|\nabla^3 h\|_2\bigr)
				+C\sigma\mathcal{M}(t)\mathcal{D}(t)\\
				&\quad\le C\mathcal M(t){\tau}\bigl(\|\nabla^2\psi\|_2^2+\|\nabla^3h\|_2^2\bigr)
				+C\sigma\mathcal M(t)\mathcal D(t)+C\sigma^2{\tau}^{-1}\mathcal N_0(t).
			\end{align*}
			\item In the second term of $\mathcal{I}_1$, integrating by parts once gives
			\[
			-\tau R\int_{\R^3}\nabla\!\left(\frac{e^{-q/\gamma}}{\bar p}\diver\psi\right)\cdot\nabla\mathcal{F}\,dx.
			\]
			Using the decomposition $\mathcal{F}=\mathcal{F}_\kappa+\mathcal{F}_\nu+\mathcal{F}_b$ from the proof of Lemma~\ref{lem:zero-acoustic}, we have
			\begin{align*}
				&{\color{black}
					\nabla\mathcal F_\kappa=
					\frac{\kappa e^{-(\gamma-1)q/\gamma}}{c_v\rho}
					\big[\nabla\Delta\zeta-e^q\bar\theta_{x_1x_1}\nabla q
					-(e^q-1)\nabla\bar\theta_{x_1x_1}\big]
					-\Big(\frac{\gamma-1}{\gamma}\nabla q+\frac{\nabla\rho}{\rho}\Big)\mathcal F_\kappa,}\\
				&{\color{black}
					\nabla\mathcal F_\nu=
					\frac{e^{-(\gamma-1)q/\gamma}}{c_v\rho}
					\nabla\big[J_b(\psi)+\mathcal Q(\psi)+(1-e^q)\nu\bar u_{x_1}^2\big]
					-\Big(\frac{\gamma-1}{\gamma}\nabla q+\frac{\nabla\rho}{\rho}\Big)\mathcal F_\nu,}\\[2pt]
				&\|\nabla\mathcal{F}_b\|_2
				\le C\sigma\bigl(\tau^{-1/2}\|\nabla\psi\|_2+\tau^{-1}\|\psi\|_2+\tau^{-1}\|\nabla(q,h)\|_2+\tau^{-3/2}\|(q,h)\|_2\bigr).
			\end{align*}
			{\color{black}Splitting $\nabla\rho=\nabla\phi+\bar\rho_{x_1}e_1$, we obtain, for the coefficient derivative in the viscous heating term,
				\[
				\|\nabla\rho\,\mathcal Q(\psi)\|_2
				\le C\|\nabla\phi\|_3\|\nabla\psi\|_6\|\nabla\psi\|_\infty
				+C\sigma\tau^{-1/2}\|\nabla\psi\|_\infty\|\nabla\psi\|_2 .
				\]
				The terms with $\Delta\zeta$ are bounded by
				$C\|\nabla U\|_3\|\Delta\zeta\|_6+C\sigma\tau^{-1/2}\|\Delta\zeta\|_2$.
				Using \eqref{eq:bg-thermal}--\eqref{eq:bg-target} for the remaining background derivatives and \eqref{eq:first-order-chain} for $h$ gives
				\[
				\tau\|\nabla\mathcal F\|_2^2
				\le C\tau\|\nabla^3 U\|_2^2+C(\sigma+\mathcal M(t))\mathcal D(t)
				+C\sigma\tau^{-1}\mathcal N_0(t).
				\]
				Moreover,
				$\|\nabla(e^{-q/\gamma}\diver\psi/\bar p)\|_2
				\le C(1+\mathcal M(t))\|\nabla^2\psi\|_2+
				C\sigma\tau^{-1/2}\|\nabla\psi\|_2$.
				Young's inequality then bounds the integral by}
			\[
			\eta_0\tau\|\nabla^2(q,\psi)\|_2^2+C_{\eta_0}\tau\|\nabla^3 U\|_2^2+C_{\eta_0}\bigl(\sigma+\mathcal{M}(t)\bigr)\mathcal{D}(t)+C_{\eta_0}\Gamma(t)\mathcal{N}(t).
			\]
			\item In the fourth term of $\mathcal{I}_1$, using $\Delta u=\Delta\psi+\bar u_{x_1x_1}e_1$ and \eqref{eq:first-order-chain} gives
			\begin{align*}
				\tau\left|\int_{\R^3}\frac{e^{-q/\gamma}}{\bar p}(\diver\psi)(\Delta u\cdot\nabla h)\,dx\right|
				&\le C\tau{\color{black}\|\nabla U\|_3}\|\diver\psi\|_6\|\Delta\psi\|_2
				+C\sigma\sqrt{\tau}\,\|\diver\psi\|_2\|\Delta\psi\|_2\\
				&\quad +C\sigma\tau^{-1/2}\|\diver\psi\|_2\|\nabla h\|_2\\
				&\le\eta_0\tau\|\nabla^2(q,\psi)\|_2^2+C_{\eta_0}\bigl(\sigma+\mathcal{M}(t)\bigr)\mathcal{D}(t).
			\end{align*}
			\item In the fifth term of $\mathcal{I}_1$, we have
			\[
			\diver u-D_u(q/\gamma+\log\bar p)
			= (\gamma+1)\bar u_{x_1}+2\diver\psi
			-\frac{(\gamma-1)(\kappa\bar\theta_{x_1x_1}+\nu\bar u_{x_1}^2)}{\bar p}
			-\psi_1\frac{\bar p_{x_1}}{\bar p}-\frac{\mathcal G}{\gamma}.
			\]
			The $2\diver\psi$ part is bounded as in (i), {\color{black}and the background terms are bounded via \eqref{eq:bg-thermal}--\eqref{eq:bg-target}. For the $\mathcal G$ part, the decomposition in the proof of Lemma~\ref{lem:zero-acoustic} gives
				\[
				\mathcal G=\frac{\mathcal F_\kappa+\mathcal F_\nu}{\bar\theta+h/\gamma}
				-\psi_1\frac{\bar p_{x_1}}{\bar p}.
				\]
				Differentiating this identity and using the estimates in (ii), with coefficient products bounded by
				$C\|\nabla U\|_3\|\Delta\zeta\|_6+
				C\sigma\tau^{-1/2}\|\Delta\zeta\|_2$ for the heat term, gives
				\begin{align*}
					&\|\mathcal G\|_2^2
					\le C\mathcal D_0(t)+C\sigma^2\tau^{-1}\mathcal N_0(t),\\
					&\tau\|\nabla\mathcal G\|_2^2
					\le C\tau\|\nabla^3U\|_2^2
					+C\bigl(\sigma+\mathcal M(t)\bigr)\mathcal D(t)
					+C\sigma\tau^{-1}\mathcal N_0(t).
				\end{align*}
				Integrating by parts,
				\begin{align*}
					&\tau R\int_{\R^3}\frac{e^{-q/\gamma}}{\bar p}
					(\diver\psi)\mathcal G\Delta h\,dx=-\tau R\int_{\R^3}
					\nabla\left(\frac{e^{-q/\gamma}}{\bar p}(\diver\psi)\mathcal G\right)
					\cdot\nabla h\,dx.
				\end{align*}
				Expand the derivative and use \eqref{eq:first-order-chain} together with
				$|\nabla h|\le C|\nabla U|+C\sigma\tau^{-1/2}|U|$.
				H\"older's inequality with exponents $(3,6,2)$, $(3,3,6,6)$ and $(\infty,2,2)$, followed by Sobolev and Young, yields
				\begin{align*}
					&\tau R\left|\int_{\R^3}\frac{e^{-q/\gamma}}{\bar p}
					(\diver\psi)\mathcal G\Delta h\,dx\right|\\
					&\le C\mathcal M(t)
					\bigl(\sqrt\tau\|\nabla^2(q,\psi)\|_2+\sigma\|\nabla\psi\|_2\bigr)
					\bigl(\sqrt\tau\|\nabla\mathcal G\|_2+\sigma\|\mathcal G\|_2\bigr)\\
					&\le\eta_0\tau\|\nabla^2(q,\psi)\|_2^2
					+C_{\eta_0}\tau\|\nabla^3U\|_2^2+C_{\eta_0}\bigl(\sigma+\mathcal M(t)\bigr)\mathcal D(t)
					+C_{\eta_0}\Gamma(t)\mathcal N(t).
			\end{align*}}
		\end{enumerate}
		Choosing $\eta_0>0$ sufficiently small so that the $\eta_0\tau\|\nabla^2(q,\psi)\|_2^2$ terms are absorbed by the left-hand side completes the proof of \eqref{eq:div-curl-ineq}.
	\end{proof}

	\section{Higher-order derivative estimates and closure of the a priori inequality}
	\label{sec:3d-high-order}
	
	In this section, we derive the unweighted and $(1+t)$-weighted higher-order energy estimates for $E_j(t)$ ($1\le j\le 4$) and the density-dissipation estimates for $X_j(t)$ ($1\le j\le 3$) from the perturbation system \eqref{eq:physical-pert-system}, and combine them with the low-order bounds from Section~\ref{sec:3d-low-order} to close the a priori energy inequality.
	
	\subsection{Higher-order energy estimates}
	
	For each integer $1\le j\le 4$, we define $E_j(t)$ by
	\[
	E_j(t)=\frac{1}{2}\sum_{|\alpha|=j}\int_{\R^3}\left(
	\frac{R\theta}{\rho}|\partial^\alpha \phi |^2
	+\rho|\partial^\alpha \psi |^2
	+\frac{c_v\rho}{\theta}|\partial^\alpha \zeta |^2\right)\,dx.
	\]
	
	\begin{lemma}
		\label{lem:higher-order-E}
		Under the assumptions of Lemma~\ref{lem:background-bounds} and \eqref{eq:apriori-assumption}, there exist constants $c,C>0$ such that
		\begin{align}
			&E_j'(t)+c\|\nabla^{j+1}(\psi,\zeta)\|_2^2
			\le C\bigl(\sigma+\mathcal{M}(t)\bigr)\mathcal{D}_0(t)
			+C\Gamma(t)\mathcal{N}_0(t),\qquad 1\le j\le 4,\label{eq:Ej-unweighted} \\[2pt]
			&(1+t)E_j'(t)+c(1+t)\|\nabla^{j+1}(\psi,\zeta)\|_2^2
			\le C\bigl(\sigma+\mathcal{M}(t)\bigr)\mathcal{D}(t)
			+C\Gamma(t)\mathcal{N}(t),\qquad 2\le j\le 4.\label{eq:Ej-weighted}
		\end{align}
	\end{lemma}
	
	\begin{proof}
		Let $(f_\rho,f_u,f_\theta)$ denote the right-hand sides of the physical perturbation system \eqref{eq:physical-pert-system} excluding the principal parabolic diffusion terms. For any multi-index $\alpha$ with $|\alpha|=j$ ($1\le j\le 4$), we define the Leibniz commutator operators and differentiated sources by
		\begin{equation}
			\begin{aligned}
				\mathcal C_\alpha(f,g)
				&=\sum_{0<\beta\le\alpha}\binom{\alpha}{\beta}
				(\partial^\beta f)\partial^{\alpha-\beta}g, \qquad
				\mathcal C_\alpha^u f
				=\sum_{0<\beta\le\alpha}\binom{\alpha}{\beta}
				(\partial^\beta u)\cdot\nabla\partial^{\alpha-\beta}f, \\[2pt]
				F_{\rho,\alpha}&=\partial^\alpha f_\rho-\mathcal C_\alpha^u\phi 
				-\mathcal C_\alpha(\rho,\diver\psi), \\[2pt]
				F_{u,\alpha}&=\partial^\alpha f_u-\mathcal C_\alpha^u \psi 
				-\mathcal C_\alpha\!\left(\frac{R\theta}{\rho},\nabla \phi\right)
				+\mathcal C_\alpha(\rho^{-1},\mathcal{L}\psi), \\[2pt]
				F_{\theta,\alpha}&=\partial^\alpha f_\theta-\mathcal C_\alpha^u\zeta 
				-(\gamma-1)\mathcal C_\alpha(\theta,\diver\psi)
				+\frac{\kappa}{c_v}\mathcal C_\alpha(\rho^{-1},\Delta \zeta).
			\end{aligned}
			\label{eq:comm-source-def}
		\end{equation}
		Applying $\partial^\alpha$ to \eqref{eq:physical-pert-system}, multiplying the three equations by the symmetrizing multipliers
		\[
		\left(\frac{R\theta}{\rho}\partial^\alpha\phi,\,\rho \partial^\alpha\psi,\,\frac{c_v\rho}{\theta}\partial^\alpha\zeta\right),
		\]
		and integrating by parts over $\R^3$, we use $(\gamma-1)c_v=R$ so that the principal hyperbolic gradient terms combine into a lower-order gradient commutator:
		\begin{align*}
			&\int_{\R^3}\Bigl(R\theta\partial^\alpha\phi\diver \partial^\alpha\psi
			+R\theta \partial^\alpha\psi\cdot\nabla\partial^\alpha\phi
			+R\rho \partial^\alpha\psi\cdot\nabla\partial^\alpha\zeta
			+R\rho\partial^\alpha\zeta\diver \partial^\alpha\psi\Bigr)\,dx\\[2pt]
			&\quad =-\int_{\R^3}\bigl(R\partial^\alpha\phi\,\partial^\alpha\psi\cdot\nabla\theta+R\partial^\alpha\zeta\,\partial^\alpha\psi\cdot\nabla\rho\bigr)\,dx.
		\end{align*}
		This gives the $j$-th order energy identity
		\begin{align}
			&E_j'(t)+\sum_{|\alpha|=j}
			\left[\mathcal{D}_u(\partial^\alpha\psi)
			+\kappa\int_{\R^3}\theta^{-1}|\nabla\partial^\alpha\zeta|^2\,dx\right]\notag \\[2pt]
			&\quad =\sum_{|\alpha|=j}\int_{\R^3}\Bigl[
			\frac{R\theta}{\rho}\partial^\alpha\phi F_{\rho,\alpha}
			+\rho\partial^\alpha\psi\cdot F_{u,\alpha}
			+\frac{c_v\rho}{\theta}\partial^\alpha\zeta F_{\theta,\alpha}
			+\frac12 B_\rho|\partial^\alpha\phi|^2\notag \\[2pt]
			&\qquad\qquad
			+\frac12 B_\theta|\partial^\alpha\zeta|^2
			+R\partial^\alpha\phi\,\partial^\alpha\psi\cdot\nabla\theta
			+R\partial^\alpha\zeta\,\partial^\alpha\psi\cdot\nabla\rho
			+\kappa\theta^{-2}\partial^\alpha\zeta\nabla\theta\cdot\nabla\partial^\alpha\zeta
			\Bigr]\,dx,\label{eq:Ej-identity}
		\end{align}
		where the material-derivative coefficients $B_\rho$ and $B_\theta$ are given by
		\begin{align*}
			B_\rho&=D_u\!\left(\frac{R\theta}{\rho}\right)+\frac{R\theta}{\rho}\diver u=(3-\gamma) \frac{R\theta}{\rho} \diver u+\frac{R\bigl(\kappa\Delta\theta+\mathcal Q(u)\bigr)}{c_v\rho^2}, \\[2pt]
			B_\theta&=D_u\!\left(\frac{c_v\rho}{\theta}\right)+\frac{c_v\rho}{\theta}\diver u=(\gamma-1) \frac{c_v\rho}{\theta} \diver u-\frac{\kappa\Delta\theta+\mathcal Q(u)}{\theta^2}.
		\end{align*}
		
		To treat both \eqref{eq:Ej-unweighted} and \eqref{eq:Ej-weighted}, set $\tau=1+t$ and let $(a,j)$ range over
		\[
		(a,j)\in\bigl\{(0,1),(0,2),(0,3),(0,4),(1,2),(1,3),(1,4)\bigr\}.
		\]
		We estimate the terms on the right-hand side of $\tau^a E_j'(t)$ in \eqref{eq:Ej-identity} in four steps.
		{\color{black}We take $\sigma,\mathcal M(t)\le1$ below.}
		
		\medskip\noindent
		\emph{Step 1.}
		By \eqref{eq:zero-order-controls}, we have
		\begin{equation}
			\tau^{-2}\|U\|_2^2
			+\tau^{-1}\|\nabla U\|_2^2
			\le C\mathcal D_0(t)+C\tau^{-1}\mathcal N_0(t).
			\label{eq:high-order-low-bound}
		\end{equation}
		Combining \eqref{eq:high-order-low-bound} with the definitions of $\mathcal{D}_0(t)$ and $\mathcal{D}_1(t)$ in \eqref{eq:dissipation} gives
		\begin{equation}
			\begin{aligned}
				&\sum_{r=0}^{j}\tau^{a-1-r}\|\nabla^{j-r}U\|_2^2
				+\sum_{r=0}^{j}\tau^{a-r}\|\nabla^{j-r}\psi\|_2^2
				+\sum_{r=1}^{j}\tau^{a+1-r}
				\|\nabla^{j+2-r}(\psi,\zeta)\|_2^2\\
				&\quad\le C\bigl(\mathcal D_0(t)+a\tau\mathcal D_1(t)\bigr)
				+C\tau^{-1}\bigl((1-a)\mathcal N_0(t)+a\mathcal N(t)\bigr).
			\end{aligned}
			\label{eq:high-order-sum-bound}
		\end{equation}
		When $r\ge 1$ derivatives in a commutator or source term fall on $(\bar\rho,\bar\theta)$ or when $r\ge 0$ derivatives fall on $(\bar u_{x_1},\bar p_{x_1},d,D_{\bar u}\bar\theta)$, Cauchy--Schwarz and Young's inequality yield
		\begin{align*}
			&\sigma\tau^{a-(r+1)/2}\|\nabla^jU\|_2\|\nabla^{j-r}\psi\|_2
			\le C\sigma\bigl(\tau^{a-1}\|\nabla^jU\|_2^2+
			\tau^{a-r}\|\nabla^{j-r}\psi\|_2^2\bigr),\\
			&\sigma\tau^{a-(r+1)/2}\|\nabla^j\psi\|_2\|\nabla^{j-r}U\|_2
			\le C\sigma\bigl(\tau^a\|\nabla^j\psi\|_2^2+
			\tau^{a-r-1}\|\nabla^{j-r}U\|_2^2\bigr),\\
			&\sigma\tau^{a-1-r/2}\|\nabla^jU\|_2\|\nabla^{j-r}U\|_2
			\le C\sigma\bigl(\tau^{a-1}\|\nabla^jU\|_2^2+
			\tau^{a-1-r}\|\nabla^{j-r}U\|_2^2\bigr),\\
			&\sigma\tau^{a-r/2}\|\nabla^j(\psi,\zeta)\|_2\|\nabla^{j-r+2}(\psi,\zeta)\|_2
			\le C\sigma\bigl(\tau^{a-1}\|\nabla^j(\psi,\zeta)\|_2^2+
			\tau^{a+1-r}\|\nabla^{j-r+2}(\psi,\zeta)\|_2^2\bigr),
		\end{align*}
		all of which are bounded by $\sigma$ times the right-hand side of \eqref{eq:high-order-sum-bound}.
		{\color{black}For a hyperbolic coefficient, a background derivative of order $r\ge1$ corresponds to the index $r-1$ in the first two estimates above; the remaining products are combined by \eqref{eq:homogeneous-product}. In a mixed parabolic term with $v\in\{\psi,\zeta\}$, the only allocation at $j=2$ gives
			\begin{align*}
				&C\sigma\tau^{a-1/2}\|\nabla\phi\|_3\|\nabla^2v\|_2\|\nabla^2v\|_6  \le C\sigma\mathcal M\bigl(\tau^{a-1}\|\nabla^2v\|_2^2+\tau^a\|\nabla^3v\|_2^2\bigr).
			\end{align*}
			For $j=3,4$, put a first density derivative in $L^\infty$, or put $\nabla^2v$ in $L^\infty$ and use \eqref{eq:homogeneous-product}. The remaining allocation uses $(3,6)$:
			\begin{align*}
				&(\|\nabla^2\phi\|_3+\|\nabla\phi\|_\infty\|\nabla\phi\|_3)\|\nabla^3v\|_6
				\le C\mathcal M\|\nabla^4v\|_2.
			\end{align*}
			Young's inequality with weights $\tau^a,\tau^{a-r}$ gives $C\sigma\mathcal M(\mathcal D_0+a\tau\mathcal D_1+\tau^{-1}\mathcal N_0)$: the remaining first and zeroth derivatives have weights at most $\tau^{-1},\tau^{-2}$ and are controlled by \eqref{eq:high-order-low-bound}.}
		
		\medskip\noindent
		\emph{Step 2.}
		We next bound the nonlinear parts of the convective and hyperbolic commutators in \eqref{eq:comm-source-def}, which consist of products of derivative orders $(r,s)$ with $r,s\ge 1$ and $r+s=j+1$, multiplied by $\partial^\alpha f$ for $f\in\{\phi,\psi,\zeta\}$.
		
		In the cases $(a,j)\in\{(0,1),(0,2),(1,2)\}$, whenever $\partial_i\psi_k$ multiplies two $j$-th order derivatives of $f\in\{\phi,\zeta\}$ ($|\alpha|=|\alpha'|=j$), we integrate $\partial_i\psi_k$ by parts once: for any smooth bounded function $b_0(\rho,\theta)$,
		\begin{align}
			\tau^a\left|\int_{\R^3} b_0(\rho,\theta)\,\partial_i \psi_k\,
			\partial^\alpha f \partial^{\alpha'}f\,dx\right|
			&\le C\tau^a\|\psi\|_3\|\nabla^j f\|_6\|\nabla^{j+1}f\|_2
			+C\tau^a\|\psi\|_3\|\nabla U\|_3\|\nabla^j f\|_6^2\notag\\
			&\quad +C\sigma\tau^{a-1/2}\|\psi\|_3\|\nabla^j f\|_2\|\nabla^j f\|_6\notag\\
			&\le C\mathcal M(t)\tau^a\|\nabla^{j+1}f\|_2^2
			+C\sigma\mathcal M(t)\tau^{a-1}\|\nabla^j f\|_2^2.
			\label{eq:low-square-ibp}
		\end{align}
		{\color{black}At $j=2$, we use \eqref{eq:low-square-ibp} when only one derivative falls on $\psi$, also for products of different scalar components. }For all remaining $(r,s)$ with $r+s=j+1$ in the convective and hyperbolic commutators, H\"older's inequality and Sobolev embedding give
		\begin{alignat*}{2}
			&j=1\ (r=1,s=1):\quad
			&&\int_{\R^3}|\nabla U|\,|\nabla\psi|\,|\nabla U|\,dx
			\le\|\nabla U\|_3{\color{black}\|\nabla U\|_6\|\nabla\psi\|_2}\\*[2pt]
			&&&\quad\le C\mathcal M(t){\color{black}\big(\|\nabla^2 U\|_2^2+\|\nabla\psi\|_2^2\big)},\\[2pt]
			&j=2\ (\{r,s\}=\{1,2\}):\quad
			&&\tau^a\int_{\R^3}|\nabla U|\,|\nabla^2 U|\,{\color{black}|\nabla^2\psi|}\,dx
			\le\tau^a\|\nabla U\|_3\|\nabla^2U\|_6{\color{black}\|\nabla^2\psi\|_2}\\*[2pt]
			&&&\quad\le C\mathcal M(t)\tau^a\|\nabla^3U\|_2{\color{black}\|\nabla^2\psi\|_2},\\[2pt]
			&j=3\ (\{r,s\}=\{1,3\}):\quad
			&&\tau^a\int_{\R^3}|\nabla U|\,|\nabla^3 U|\,|\nabla^3 U|\,dx
			\le C\mathcal{M}(t)\tau^a\|\nabla^3 U\|_2^2,\\[2pt]
			&j=3\ (r=2,s=2):\quad
			&&\tau^a\int_{\R^3}|\nabla^2 U|\,|\nabla^2 U|\,|\nabla^3 U|\,dx
			\le\tau^a\|\nabla^2 U\|_3\|\nabla^2 U\|_6\|\nabla^3 U\|_2\\*[2pt]
			&&&\quad\le C\mathcal{M}(t)\tau^a\|\nabla^3 U\|_2^2,\\[2pt]
			&j=4\ (r\in\{1,2\},\,s=5-r):\quad
			&&\tau^a\int_{\R^3}|\nabla^{r}U|\,|\nabla^{5-r}U|\,|\nabla^4 U|\,dx\\*
			&&&\quad\le\tau^a\|U\|_{W^{2,\infty}}\bigl(\|\nabla^3 U\|_2+\|\nabla^4 U\|_2\bigr)\|\nabla^4 U\|_2\\*
			&&&\quad\le C\mathcal{M}(t)\tau^a\bigl(\|\nabla^3 U\|_2^2+\|\nabla^4 U\|_2^2\bigr).
		\end{alignat*}
		{\color{black}At $j=2$, the chain rule and \eqref{eq:homogeneous-product}, retaining one and two background derivatives, give the following bound. The H\"older exponents are $(6,3,2)$, $(3,6,2)$ and $(6,3,2)$, respectively:
			\begin{align*}
				&\tau^a\int_{\R^3}\left|\nabla^2\left(\frac{R\theta}{\rho}-\frac{R\bar\theta}{\bar\rho}\right)\right|
				|\nabla\phi|\,|\nabla^2\psi|\,dx\\
				&\le C\mathcal M\tau^a\|\nabla^3U\|_2\|\nabla^2\psi\|_2+C\sigma\mathcal M\bigl(\tau^{a-1/2}\|\nabla^2\phi\|_2
				+\tau^{a-1}\|\nabla U\|_2\bigr)\|\nabla^2\psi\|_2\\
				&\le C\mathcal M(\mathcal D_0+a\tau\mathcal D_1)
				+C\sigma\mathcal M\tau^{-1}\mathcal N_0,
			\end{align*}
			where Young's inequality and \eqref{eq:high-order-low-bound} were used in the last step.}
		
		\medskip\noindent
		\emph{Step 3.}
		For $v\in\{\psi,\zeta\}$, the parabolic commutators in $F_{u,\alpha}$ and $F_{\theta,\alpha}$ involve products of derivative orders $(r,j+2-r)$ ($1\le r\le j$) multiplied by $\partial^\alpha v$. In the two cases $(j,r)=(2,2)$ and $(j,r)=(4,3)$, we integrate by parts once: for any smooth bounded function $b_1(\rho,\theta)$,
		\begin{align*}
			\int_{\R^3} b_1(\rho,\theta)\,\partial^{\alpha'}(\rho^{-1})\,\partial_i\partial_kv\,\partial^\alpha v\,dx
			&=-\int_{\R^3}\partial_i\bigl(b_1(\rho,\theta)\partial^{\alpha'}(\rho^{-1})\partial^\alpha v\bigr)\partial_kv\,dx,\\
			\int_{\R^3} b_1(\rho,\theta)\,\partial_i\partial^{\alpha'}(\rho^{-1})\,
			\partial^\beta v\,\partial^{\beta'}v\,dx
			&=-\int_{\R^3}\partial^{\alpha'}(\rho^{-1})
			\partial_i\bigl(b_1(\rho,\theta)\partial^\beta v\,\partial^{\beta'}v\bigr)\,dx,
		\end{align*}
		where $|\alpha|=|\alpha'|=2$, $|\beta|=3$, and $|\beta'|=4$.
		Applying \eqref{eq:homogeneous-product} together with these two identities for $(j,r)=(2,2),(4,3)$ and H\"older's inequality for the remaining $(j,r)$ pairs gives
		\begin{alignat*}{2}
			&(j,r)=(1,1),(2,1):\quad
			&&\tau^a\|\nabla\phi\|_3\|\nabla^{j+1}v\|_2\|\nabla^j v\|_6
			\le C\mathcal{M}(t)\tau^a\|\nabla^{j+1}v\|_2^2,\\[2pt]
			&(j,r)=(2,2):\quad
			&&\tau^a\bigl(\|\nabla^3\phi\|_2\|\nabla^2 v\|_6+\|\nabla^2\phi\|_6\|\nabla^3 v\|_2\bigr)\|\nabla v\|_3\\*
			&&&\quad\le C\mathcal{M}(t)\tau^a\bigl(\|\nabla^3 U\|_2^2+\|\nabla^3 v\|_2^2\bigr),\\[2pt]
			&j=3\ (r=1,2,3):\quad
			&&\tau^a\bigl(\|\nabla\phi\|_\infty\|\nabla^4 v\|_2
			+\|\nabla^2\phi\|_3\|\nabla^3 v\|_6
			+\|\nabla^3\phi\|_2\|\nabla^2 v\|_\infty\bigr)\|\nabla^3 v\|_2\\*
			&&&\quad\le C\mathcal{M}(t)\tau^a\bigl(\|\nabla^3 U\|_2^2+\|\nabla^4 v\|_2^2\bigr),\\[2pt]
			&j=4\ (r=1,2,3,4):\quad
			&&\tau^a\bigl(\|\nabla\phi\|_{W^{1,\infty}}(\|\nabla^4 v\|_2+\|\nabla^5 v\|_2)\|\nabla^4 v\|_2\\*
			&&&\quad +\|\nabla^4\phi\|_2\|\nabla^2 v\|_\infty\|\nabla^4 v\|_2{\color{black}+\|\nabla^2\phi\|_3\|\nabla^3v\|_6\|\nabla^5v\|_2}\bigr)\\*
			&&&\quad\le C\mathcal{M}(t)\tau^a\bigl(\|\nabla^4 U\|_2^2+\|\nabla^4 v\|_2^2+\|\nabla^5 v\|_2^2\bigr).
		\end{alignat*}
		
		{\color{black}For the $(4,3)$ integration by parts, derivatives of the multiplier contribute at most
			$C(\mathcal M(t)+\sigma\tau^{-1/2})\tau^a
			\|\nabla^2\phi\|_3\|\nabla^3v\|_6\|\nabla^4v\|_2$,
			which is also bounded by the preceding right-hand side. Background factors are estimated as in Step~1.}
		
		{\color{black}For $(j,r)=(2,2)$, the derivative of the multiplier also contributes
			\begin{align*}
				&C\tau^a\|\nabla U\|_3\|\nabla^2\phi\|_6\|\nabla v\|_3\|\nabla^2v\|_6+C\sigma\tau^{a-1/2}\|\nabla^2\phi\|_6\|\nabla v\|_3\|\nabla^2v\|_2\\
				&\le C\mathcal M^2\tau^a(\|\nabla^3\phi\|_2^2+\|\nabla^3v\|_2^2)+C\sigma\mathcal M(\tau^a\|\nabla^3\phi\|_2^2+\tau^{a-1}\|\nabla^2v\|_2^2)\\
				&\le C(\mathcal M^2+\sigma\mathcal M)(\mathcal D_0+a\tau\mathcal D_1).
			\end{align*}
			Coefficient compositions have the same bounds by \eqref{eq:homogeneous-product}.}
		
		Similarly, in the differentiated viscous heating integral, when all $j$ derivatives act on $\mathcal{Q}(\psi)$ and we integrate one derivative by parts in the symmetric case $(j,r)=(4,3)$, we obtain
		{\color{black}the following estimates, with the last product in each pair present only for $j=2$ and $j=4$, respectively:}
		\begin{align*}
			j=1,2:\quad
			&\tau^a\|\nabla\psi\|_3\|\nabla^{j+1}\psi\|_2\|\nabla^j\zeta\|_6+\tau^a\|\nabla^2\psi\|_3\|\nabla^2\psi\|_2\|\nabla^2\zeta\|_6\\
			&\le C\mathcal{M}(t)\tau^a\bigl(\|\nabla^{j+1}(\psi,\zeta)\|_2^2+\|\nabla^2\psi\|_2^2\bigr),\\[2pt]
			j=3,4:\quad
			&\tau^a\bigl(\|\nabla\psi\|_\infty\|\nabla^{j+1}\psi\|_2+\|\nabla^2\psi\|_\infty\|\nabla^j\psi\|_2\bigr)\|\nabla^j\zeta\|_2+\tau^a\|\nabla^2\psi\|_\infty\|\nabla^3\psi\|_2\|\nabla^5\zeta\|_2\\
			&\le {\color{black}C\mathcal{M}(t)\bigl(\mathcal D_0(t)+a\tau\mathcal D_1(t)\bigr)}.
		\end{align*}
		{\color{black}The multiplier of $\partial^\alpha\mathcal Q(\psi)\partial^\alpha\zeta$ is $\theta^{-1}$. At $j=4$, differentiating this multiplier in the above integration by parts also gives terms bounded by
			\begin{align*}
				&C\tau^a\int_{\R^3}|\nabla(\theta^{-1})|\,|\nabla^2\psi|\,
				|\nabla^3\psi|\,|\nabla^4\zeta|\,dx\\
				&\quad\le C\bigl(\mathcal M(t)^2+\sigma\mathcal M(t)\bigr)\tau^a
				\bigl(\|\nabla^3\psi\|_2^2+\|\nabla^4\zeta\|_2^2\bigr)\\
				&\quad\le C\mathcal M(t)\bigl(\mathcal D_0(t)+a\tau\mathcal D_1(t)\bigr),
			\end{align*}
			where we used $|\nabla(\theta^{-1})|\le C|\nabla\zeta|+C\sigma\tau^{-1/2}$, $\tau\ge1$, and $\mathcal M(t),\sigma\le1$. In particular, the third velocity derivative and the fifth temperature derivative in the above heating estimates are both controlled by the full dissipation for $a=0,1$.}
		{\color{black}When $1\le s\le j$ derivatives fall on $(c_v\rho)^{-1}$, the terms with all derivatives on perturbations are estimated as follows:
			\begin{flalign*}
				&\begin{alignedat}{2}
					&j=s=1: &\quad &C\|\nabla\phi\|_6\|\nabla\psi\|_3^2\|\nabla\zeta\|_6,\\
					&j=2,\ s=1: &\quad &C\tau^a\|\nabla\phi\|_3\|\nabla\psi\|_3\|\nabla^2\psi\|_6\|\nabla^2\zeta\|_6,\\
					&j=s=2: &\quad &C\tau^a\|\nabla^2\phi\|_6\|\nabla\psi\|_3^2\|\nabla^2\zeta\|_6
					+C\tau^a\|\nabla\phi\|_3^2\|\nabla\psi\|_\infty\|\nabla\psi\|_6\|\nabla^2\zeta\|_6 .
				\end{alignedat} &&
			\end{flalign*}
			By Sobolev embedding and \eqref{eq:M-def}, these products are bounded by
			$C\mathcal M(t)\mathcal D_0(t)$ for $a=0$ and by $C\mathcal M(t)\mathcal D(t)$ for $a=1$.
			For $j=3,4$, the product estimate used in \eqref{eq:homogeneous-product}, applied to the density factors and to $\nabla\psi$, gives the bound
			\[
			C\tau^a\big(\|\nabla\psi\|_\infty\|\nabla^{j+1}\psi\|_2+
			\|\nabla\psi\|_\infty^2\|\nabla^j\phi\|_2\big)\|\nabla^j\zeta\|_2 .
			\]
			Each group of $r\ge1$ derivatives falling on the background contributes at most
			$C\sigma\tau^{-r/2}$ by \eqref{eq:bg-thermal}. For $j=1,2$, put $\nabla^j\zeta$ in $L^6$, one velocity gradient in $L^3$, and the remaining velocity factor in $L^2$; any remaining first derivative of $\phi$ is bounded in $L^\infty$. For $j=3,4$, apply the same product estimate to the remaining $j-r$ derivatives. Young's inequality and \eqref{eq:high-order-sum-bound} give the required bound in each case.}
		{\color{black}For $j=3,4$, apply H\"older with $(2j/s,2j/(j-s))$ (or $(2,\infty)$ when $s=j$) to the pure coefficient and heating derivatives, then use homogeneous interpolation and Young's inequality. The last low-order product uses $\tau\|\nabla^2\psi\|_2^2$ in $\tau\mathcal D_1$.
			
			The background terms at $j=1,2$ are bounded, respectively, by
			\begin{align*}
				j=1:\quad &C\sigma\mathcal M\tau^{-1/2}\|\nabla\psi\|_2\|\nabla^2\zeta\|_2,\\
				j=2:\quad &C\sigma\mathcal M\tau^{a-1/2}\|\nabla^2\psi\|_2\|\nabla^3\zeta\|_2+C\sigma\mathcal M(\tau^{a-1}+\mathcal M\tau^{a-1/2})\|\nabla\psi\|_2\|\nabla^3\zeta\|_2.
			\end{align*}
			Both are at most $C\sigma\mathcal M(\mathcal D_0+a\tau\mathcal D_1)$. For $j=3,4$, a background group of order $1\le r\le j$ gives
			\begin{align*}
				&C\sigma\mathcal M\tau^{a-r/2}
				(\|\nabla^{j-r+1}\psi\|_2+\mathcal M\|\nabla^{j-r}\phi\|_2)\|\nabla^j\zeta\|_2\\
				&\quad\le C\sigma\mathcal M(\mathcal D_0+a\tau\mathcal D_1+\tau^{-1}\mathcal N_0),
			\end{align*}
			by Young's inequality with weights $\tau^a,\tau^{a-r}$ and the mixed-background bound in Step~1.}
		
		\medskip\noindent
		\emph{Step 4.}
		In \eqref{eq:Ej-identity}, we write the heating terms inside $B_\rho$ and $B_\theta$ as
		\[
		\kappa\Delta\theta+\mathcal Q(u)
		=\bigl(\kappa\bar\theta_{x_1x_1}+\nu\bar u_{x_1}^2\bigr)+\kappa\Delta\zeta 
		+J_b(\psi)+\mathcal Q(\psi).
		\]
		The background parts of $B_\rho$ and $B_\theta$ are bounded by $\sigma$ times the right-hand side of \eqref{eq:high-order-sum-bound}. For $j=3,4$ and $(a,j)=(0,2)$, the remaining parts of $B_\rho$ and $B_\theta$ are bounded in $L^\infty(\R^3)$ by $C\mathcal{M}(t)$, while $\|\nabla^j(\phi,\zeta)\|_2^2$ is part of the dissipation. In the cases $(a,j)\in\{(0,1),(1,2)\}$, the $\diver\psi$ term is bounded by \eqref{eq:low-square-ibp}, and for the $\Delta\zeta$ term we integrate by parts once: for $f\in\{\phi,\zeta\}$ and $b_2(\rho,\theta)\in\{R/(c_v\rho^2),-1/\theta^2\}$,
		\begin{align*}
			&\tau^a\left|\int_{\R^3} b_2(\rho,\theta)\Delta\zeta\,(\partial^\alpha f)^2\,dx\right|\\
			& =\tau^a\left|-2\int_{\R^3} b_2(\rho,\theta)\nabla\zeta\,\partial^\alpha f\,\nabla\partial^\alpha f\,dx
			-\int_{\R^3}\nabla b_2(\rho,\theta)\cdot\nabla\zeta\,(\partial^\alpha f)^2\,dx\right|\\
			&\le C\tau^a\|\nabla\zeta\|_3\|\partial^\alpha f\|_6\|\nabla\partial^\alpha f\|_2
			+C\tau^a\|\nabla U\|_3\|\nabla\zeta\|_3\|\partial^\alpha f\|_6^2
			+C\sigma\tau^{a-1/2}{\color{black}\|\nabla\zeta\|_3\|\partial^\alpha f\|_6\|\partial^\alpha f\|_2}\\
			& \le C\mathcal{M}(t)\tau^a\|\nabla^{j+1}f\|_2^2+C\sigma\mathcal{M}(t)\tau^{a-1}\|\nabla^j f\|_2^2.
		\end{align*}
		{\color{black}For $(a,j)=(0,1),(1,2)$ and $f\in\{\phi,\zeta\}$, the remaining heating parts of $B_\rho,B_\theta$ satisfy
			\begin{align*}
				&\tau^a\int_{\R^3}|\nabla\psi|^2|\partial^\alpha f|^2\,dx
				\le\tau^a\|\nabla\psi\|_3^2\|\partial^\alpha f\|_6^2\le C\mathcal M^2\tau^a\|\nabla^{j+1}f\|_2^2,\\
				&\tau^a\int_{\R^3}|J_b(\psi)|\,|\partial^\alpha f|^2\,dx
				\le C\sigma\mathcal M\tau^{a-1}\|\nabla^jf\|_2^2.
			\end{align*}
			For the differentiated linear heating source, Leibniz' rule gives
			\begin{align*}
				&\left\|\partial^\alpha\!\left(\frac{J_b(\psi)}{c_v\rho}\right)\right\|_2\le C\sigma\sum_{r=0}^j\tau^{-1-r/2}
				(\|\nabla^{j-r+1}\psi\|_2+\mathcal M\|\nabla^{j-r}\phi\|_2).
			\end{align*}
			After pairing with $\tau^a\|\nabla^j\zeta\|_2$, Young's inequality with weights $\tau^{a-1},\tau^{a-1-r}$ and \eqref{eq:high-order-sum-bound} control these terms; the $r=0$ velocity term is covered by its last sum.}
		
		Similarly, splitting the density and temperature gradients into background and perturbation parts in the last three flux terms of \eqref{eq:Ej-identity} and applying H\"older's inequality ($L^3\times L^6\times L^2$ for $1\le j\le 3$ and $L^\infty\times L^2\times L^2$ for $j=4$) together with Step~1 and \eqref{eq:homogeneous-product} controls all remaining flux and source terms. In particular, because $j\ge 1$, the differentiated momentum coupling with $d$ satisfies
		\[
		\nabla(\rho,\theta)=\nabla(\phi,\zeta)+(\bar\rho_{x_1},\bar\theta_{x_1})e_1,\qquad
		\|\partial^\alpha(\rho^{-1}\phi d)\|_2
		\le C\sigma\sum_{i=0}^j\tau^{-1-(j-i)/2}\|\nabla^i U\|_2,
		\]
		and is controlled directly by Step~1.

		Combining Steps~1--4 with \eqref{eq:high-order-sum-bound} and using $\Gamma(t)\ge\sigma\tau^{-1}$ completes the proof of \eqref{eq:Ej-unweighted}--\eqref{eq:Ej-weighted}.
	\end{proof}

	\subsection{Dissipation estimates for the density derivatives}
	
	To obtain dissipation for the derivatives of the density perturbation $\phi$, we define
	\[
	X_j(t)=\sum_{|\alpha|=j}\int_{\R^3}
	\partial^\alpha\psi\cdot\nabla\partial^\alpha(\phi /\bar\rho)\,dx,
	\qquad 1\le j\le3.
	\]
	
	\begin{lemma}
		\label{lem:density-dissipation}
		Under the assumptions of Lemma~\ref{lem:background-bounds} and \eqref{eq:apriori-assumption}, there exist constants $c,C>0$ such that
		\begin{align}
			&X_j'(t)+c\|\nabla^{j+1}\phi\|_2^2
			\le C\bigl(\|\nabla^{j+1}(\psi,\zeta)\|_2^2
			+\|\nabla^{j+2}\psi\|_2^2\bigr)\notag\\
			&\quad+C\bigl(\sigma+\mathcal M(t)\bigr)\mathcal D_0(t)
			+C\Gamma(t)\mathcal N_0(t),
			\qquad 1\le j\le3,\label{eq:Xj-unweighted}\\
			&(1+t)X_j'(t)+c(1+t)\|\nabla^{j+1}\phi\|_2^2
			\le C(1+t)\bigl(\|\nabla^{j+1}(\psi,\zeta)\|_2^2
			+\|\nabla^{j+2}\psi\|_2^2\bigr)\notag\\
			&\quad+C\bigl(\sigma+\mathcal M(t)\bigr)\mathcal D(t)
			+C\Gamma(t)\mathcal N(t),
			\qquad 2\le j\le3.\label{eq:Xj-weighted}
		\end{align}
	\end{lemma}
	
	\begin{proof}
		Set $\tau=1+t$. Dividing the continuity equation in \eqref{eq:physical-pert-system} by $\bar\rho$ gives the transport equation for $\phi/\bar\rho$:
		\[
		D_u\!\left(\frac{\phi}{\bar\rho}\right)=-\left(1+\frac{\phi}{\bar\rho}\right)\left(\diver\psi+\frac{\bar\rho_{x_1}}{\bar\rho}\psi_1\right).
		\]
		For any multi-index $\alpha$ with $|\alpha|=j$ ($1\le j\le 3$), set
		\begin{align*}
			&r_\alpha =\nabla\partial^\alpha\!\left(\frac{\phi}{\bar\rho}\right)
			=\frac{1}{\bar\rho}\nabla\partial^\alpha\phi+C_\alpha,
			\qquad C_\alpha=[\nabla\partial^\alpha,\bar\rho^{-1}]\phi,\\[2pt]
			&F_{0,\alpha} =-\left(1+\frac{\phi}{\bar\rho}\right)\diver \partial^\alpha\psi+R_{0,\alpha},\\[2pt]
			&R_{0,\alpha} =-\mathcal C_\alpha\!\left(\frac{\phi}{\bar\rho},\diver\psi\right)
			-\partial^\alpha\!\left(\left(1+\frac{\phi}{\bar\rho}\right)\frac{\bar\rho_{x_1}}{\bar\rho}\psi_1\right)
			-\mathcal C_\alpha^u\!\left(\frac{\phi}{\bar\rho}\right),
		\end{align*}
		and let $F_{u,\alpha}$ be as defined in the proof of Lemma~\ref{lem:higher-order-E}. Differentiating $X_j(t)$ in time along \eqref{eq:physical-pert-system} and integrating by parts over $\R^3$ yields the identity
		\begin{align}
			X_j'(t)+\sum_{|\alpha|=j}\int_{\R^3} \frac{R\theta}{\rho\bar\rho}|\nabla\partial^\alpha\phi|^2\,dx
			&=\sum_{|\alpha|=j}\int_{\R^3}\Bigl[
			(-R\nabla\partial^\alpha\zeta+\rho^{-1}\mathcal L\partial^\alpha\psi
			+F_{u,\alpha})\cdot r_\alpha\notag\\
			&\quad -\frac{R\theta}{\rho}\nabla\partial^\alpha\phi\cdot C_\alpha
			+\left(1+\frac{\phi}{\bar\rho}\right)|\diver \partial^\alpha\psi|^2
			-\diver \partial^\alpha\psi\,R_{0,\alpha}\notag\\
			&\quad -\partial^\alpha\psi\cdot(\nabla u)^{\top}r_\alpha
			+(\diver u)\partial^\alpha\psi\cdot r_\alpha\Bigr]dx.
			\label{eq:Xj-identity}
		\end{align}
		
		To treat both \eqref{eq:Xj-unweighted} and \eqref{eq:Xj-weighted}, let $(j,a)\in\{(1,0),(2,0),(3,0),(2,1),(3,1)\}$. By \eqref{eq:high-order-low-bound} and the definition of $\mathcal{D}_0(t)$ in \eqref{eq:dissipation}, the lower-order background derivative sum satisfies
		\begin{equation}
			\tau^a\sum_{i=0}^j\tau^{-(j+1-i)}\|\nabla^iU(t)\|_2^2
			\le C\mathcal D_0(t)+C\tau^{-1}\mathcal N_0(t).
			\label{eq:density-sum-bound}
		\end{equation}
		By \eqref{eq:bg-thermal}, the background density derivatives satisfy
		\begin{gather*}
			\|\partial_{x_1}^i(\bar\rho^{-1})\|_\infty\le C\sigma\tau^{-i/2}\quad {\color{black}(1\le i\le j+1\le4)},\\[2pt]
			\|\partial_{x_1}^i(\bar\rho_{x_1}/\bar\rho)\|_\infty\le C\sigma\tau^{-(i+1)/2}\quad {\color{black}(0\le i\le j\le3)},
		\end{gather*}
		so that the commutator $C_\alpha$ obeys
		\[
		\|C_\alpha\|_2\le C\sigma\sum_{i=0}^j\tau^{-(j+1-i)/2}\|\nabla^i U\|_2.
		\]
		Moreover, we claim that the remainder quantities in \eqref{eq:Xj-identity} satisfy
		\begin{equation}
			\begin{aligned}
				&\tau^a\sum_{|\alpha|=j}\left(
				\|F_{u,\alpha}\|_2^2+\|R_{0,\alpha}\|_2^2
				+\||\nabla u|\,|\partial^\alpha\psi|\|_2^2+\|C_\alpha\|_2^2\right)\\
				&\quad\le C\bigl(\sigma+\mathcal M(t)\bigr)^2\bigl(\mathcal D_0(t)+a\tau\mathcal D_1(t)\bigr)
				+C\sigma^2\tau^{-1}\mathcal N_0(t).
			\end{aligned}
			\label{eq:density-remainder-bound}
		\end{equation}
		To verify \eqref{eq:density-remainder-bound}, we check each term on the left-hand side:
		\begin{enumerate}[(i)]
			\item In the convective and pressure commutators inside $F_{u,\alpha}$, as well as in $R_{0,\alpha}$, every nonlinear product of perturbation derivatives has total order at most $j+1\le 4$ and is bounded in $L^2(\R^3)$ by $C\mathcal{M}(t)\|\nabla^{j+1}U\|_2$ via \eqref{eq:homogeneous-product}. All terms in which one or more derivatives hit the background state are bounded via \eqref{eq:density-sum-bound}.
			\item For the commutator $\mathcal{C}_\alpha(\rho^{-1},\mathcal{L}\psi)$ inside $F_{u,\alpha}$, H\"older's inequality ($L^3\times L^6$ for $j=1,2$ and $L^\infty\times L^2$ for $j=3$) and Sobolev embedding give
			\begin{align*}
				j=1:\quad
				&\|\nabla\phi\|_3\|\nabla^2\psi\|_6
				\le C\mathcal{M}(t)\|\nabla^3\psi\|_2,\\[2pt]
				j=2:\quad
				&\|\nabla\phi\|_3\|\nabla^3\psi\|_6
				+\|\nabla^2\phi\|_6\|\nabla^2\psi\|_3
				\le C\mathcal{M}(t)\bigl(\|\nabla^4\psi\|_2+\|\nabla^3 U\|_2\bigr),\\[2pt]
				j=3:\quad
				&\|\nabla\phi\|_{W^{1,\infty}}\bigl(\|\nabla^4\psi\|_2+\|\nabla^3\psi\|_2\bigr)
				+\|\nabla^3\phi\|_2\|\nabla^2\psi\|_\infty\\
				&\quad\le C\mathcal{M}(t)\bigl(\|\nabla^4\psi\|_2+\|\nabla^3\psi\|_2+\|\nabla^3\phi\|_2\bigr),
			\end{align*}
			while terms with $i\ge 1$ derivatives on $\bar\rho^{-1}$ satisfy
			\[
			\tau^a\sum_{i=1}^j\sigma^2\tau^{-i}\|\nabla^{j+2-i}\psi\|_2^2
			\le C\sigma^2\mathcal{D}_0(t).
			\]
			\item For the differentiated momentum source $\partial^\alpha f_u$ and the geometric velocity products, Leibniz's rule and \eqref{eq:bg-thermal}--\eqref{eq:bg-target} give
			\begin{align*}
				&\|\partial^\alpha f_u\|_2
				\le C\sigma\sum_{i=0}^j\tau^{-(j+1-i)/2}\|\nabla^i U\|_2,\\[2pt]
				&\||\nabla u|\,|\partial^\alpha\psi|\|_2
				\le C\mathcal{M}(t)\bigl(\|\nabla^{j+1}\psi\|_2+\|\nabla^j\psi\|_2\bigr)+C\sigma\tau^{-1}\|\nabla^j\psi\|_2.
			\end{align*}
			Squaring and multiplying by $\tau^a$ via \eqref{eq:density-sum-bound} proves \eqref{eq:density-remainder-bound}.
		\end{enumerate}
		
		Applying Cauchy--Schwarz and Young's inequality to \eqref{eq:Xj-identity} and substituting \eqref{eq:density-remainder-bound} establishes \eqref{eq:Xj-unweighted} and \eqref{eq:Xj-weighted}.
	\end{proof}
	
	\subsection{Closure of the a priori energy inequality}
	
	\begin{proposition}
		\label{prop:apriori-closure}
		There exist constants $c,C,\delta_1>0$, independent of $\delta$ and $t$, and an energy functional $\mathcal{E}(t)$ such that every smooth solution satisfying \eqref{eq:apriori-assumption} on $[0,T]$ obeys
		\begin{equation}
			\begin{aligned}
				&C^{-1}\mathcal N(t)\le\mathcal{E}(t)\le C\mathcal N(t),
				\qquad \mathcal{E}'(t)+c\mathcal D(t)\le C\Gamma(t)\mathcal{E}(t),\\[2pt]
				&\Gamma(t)=\frac{\sigma}{1+t}
				+\sqrt{\frac{\Lambda(t)}{1+t}}+\Lambda(t).
			\end{aligned}
			\label{eq:apriori-closed}
		\end{equation}
	\end{proposition}
	
	\begin{proof}
		Set $\tau=1+t$ and recall from \eqref{eq:Nstar} that
		\[
		\mathcal N_1(t)=\|\nabla(q,\psi)(t)\|_2^2
		+\sum_{j=2}^4\|\nabla^jU(t)\|_2^2,
		\qquad \mathcal N(t)=\mathcal N_0(t)+\tau\mathcal N_1(t).
		\]
		With $\alpha_0>0$ fixed as in Lemma~\ref{lem:entropy-estimate}, we define
		\begin{align*}
			&\mathcal F_0(t) =E_p(t)+\lambda_0E_h(t)-\frac{\kappa}{\gamma c_p}Y_0(t)
			+\lambda_1\sum_{j=1}^4E_j(t)+\lambda_2\sum_{j=1}^3X_j(t),\\
			&\mathcal G_1(t) =E_d(t)+E_c(t)-\frac{\kappa}{\gamma c_p}Y_1(t)
			+\mu_1\sum_{j=2}^4E_j(t)
			+\mu_2\sum_{j=2}^3X_j(t),\\
			&\mathcal{E}(t) =K_0\mathcal F_0(t)+\tau \mathcal G_1(t),
		\end{align*}
		where the positive constants $\lambda_0,\lambda_2,\lambda_1,\mu_2,\mu_1,K_0$ are chosen as follows.
		
		First, we check that $\mathcal{E}(t)\sim\mathcal{N}(t)$. By \eqref{eq:qh-chain} and $\tau^{-1}\|U\|_2^2\le C E_0(t)$, for any $\eta_0>0$ the functionals $Y_0(t)$ and $Y_1(t)$ satisfy
		\begin{align}
			\frac{\kappa}{\gamma c_p}|Y_0(t)|&\le\eta_0\|\psi\|_2^2
			+C_{\eta_0}\|\nabla U\|_2^2
			+C_{\eta_0}\sigma^2E_0(t),\label{eq:Y0-coercivity}\\
			\frac{\kappa}{\gamma c_p}|Y_1(t)|&\le\eta_0\|\diver\psi\|_2^2
			+C_{\eta_0}\|\nabla^2U\|_2^2
			+C_{\eta_0}\sigma^2\bigl(
			\tau^{-1}\|\nabla U\|_2^2
			+\tau^{-2}\|U\|_2^2\bigr).
			\label{eq:Y1-coercivity}
		\end{align}
		Choosing $\eta_0>0$ small relative to the kinetic energy in $E_p(t)$ and the divergence term in $E_d(t)$ fixes $C_{\eta_0}$. Similarly, $X_j(t)$ and the curl difference in $E_c(t)$ satisfy
		\begin{align*}
			&|X_j(t)| \le C\|\nabla^j\psi\|_2
			\left(\|\nabla^{j+1}\phi\|_2
			+\sigma\sum_{i=0}^j
			\tau^{-(j+1-i)/2}\|\nabla^i\phi\|_2\right),\\[2pt]
			&\|\nabla\rho\times \psi\|_2
			\le C\|\nabla\phi\|_3\|\psi\|_6
			+C\sigma\tau^{-1/2}\|\psi\|_2
			\le C\mathcal M(t)\|\nabla \psi\|_2
			+C\sigma\tau^{-1/2}\|\psi\|_2.
		\end{align*}
		Using $\|\nabla\psi\|_2^2=\|\diver\psi\|_2^2+\|\curl\psi\|_2^2$, we see that once $\lambda_1\gg\lambda_2+C_{\eta_0}$ and $\mu_1\gg\mu_2+C_{\eta_0}$ are chosen large and $\sigma,\mathcal{M}(t)$ are small,
		\begin{equation}
			\mathcal F_0(t)\sim\mathcal N_0(t),\qquad
			c\mathcal N_1(t)-C\sigma^2\tau^{-1}\mathcal N_0(t)
			\le \mathcal G_1(t)\le C\mathcal N_1(t)+C\sigma^2\tau^{-1}\mathcal N_0(t),
			\label{eq:assembled-coercivity}
		\end{equation}
		and hence $\mathcal{E}(t)\sim\mathcal{N}(t)$ for any $K_0\ge 1$.
		
		Second, let $c_1,c_2,c_3,c_4,c_5>0$ be the dissipation constants on the left-hand sides of Lemmas~\ref{lem:zero-acoustic}, \ref{lem:entropy-estimate}, \ref{lem:first-order-div-curl}, \ref{lem:higher-order-E}, and \ref{lem:density-dissipation}, respectively, and let $C_0$ be the coefficient of $\|\nabla q\|_2^2$ on the right-hand side of \eqref{eq:entropy-ineq}. We fix
		\[
		0<\lambda_0\le\min\left\{1,\,\frac{c_1}{4C_0}\right\},\qquad
		\eta_0=\frac{1}{4}\lambda_0 c_2,
		\]
		and apply Lemma~\ref{lem:zero-acoustic} with this $\eta_0$, which fixes the constant $C_2=C_{\eta_0}$ in \eqref{eq:zero-acoustic-ineq}. Denoting by $C_3$ the constant on the right-hand side of Lemma~\ref{lem:density-dissipation} and by $C_1$ the constant in \eqref{eq:div-curl-ineq}, we choose $\lambda_2,\lambda_1,\mu_2,\mu_1$ so that
		\[
		\begin{aligned}
			\lambda_2c_5&\ge 4C_2,&\qquad
			\lambda_1 c_4&\ge 4(C_2+3\lambda_2C_3+1),\\[2pt]
			\mu_2c_5&\ge 4C_1,&\qquad
			\mu_1 c_4&\ge 4(C_1+2\mu_2C_3+1),
		\end{aligned}
		\]
		increasing $\lambda_1$ and $\mu_1$ if needed so that \eqref{eq:assembled-coercivity} holds. Summing the inequalities in Lemmas~\ref{lem:entropy-estimate}--\ref{lem:density-dissipation} then gives constants $a_0,a_1,C_4>0$ such that
		\begin{align*}
			&\mathcal F_0'(t)+a_0\mathcal D_0(t)
			\le C_4\bigl(\sigma+\mathcal M(t)\bigr)\mathcal D_0(t)
			+C_4\Gamma(t)\mathcal N_0(t),\\
			&\tau \mathcal G_1'(t)+a_1\tau\mathcal D_1(t)
			\le C_4\bigl(\sigma+\mathcal M(t)\bigr)\mathcal D(t)
			+C_4\Gamma(t)\mathcal N(t).
		\end{align*}
		
		Finally, when differentiating $\mathcal{E}(t)=K_0\mathcal{F}_0(t)+\tau\mathcal{G}_1(t)$, the derivative of $\tau=1+t$ produces $\mathcal{G}_1(t)$, which by \eqref{eq:Y1-coercivity} and \eqref{eq:density-sum-bound} satisfies
		\begin{equation}
			|\mathcal G_1(t)|\le C_5\mathcal D_0(t)
			+C_5\sigma^2\tau^{-1}\mathcal N_0(t)
			\le C_5\mathcal D_0(t)+C_5\Gamma(t)\mathcal N_0(t).
			\label{eq:G1-time-weight-cost}
		\end{equation}
		Fixing $K_0=\max\{1,2C_5/a_0\}$ and choosing $\delta_0,\delta_1>0$ small enough that
		\[
		C_5\mathcal{D}_0(t)\le\frac{1}{2}K_0 a_0\mathcal{D}_0(t),\qquad
		(K_0+1)C_4(\sigma+\delta_1)\le\frac{1}{4}\min\{K_0 a_0,a_1\},
		\]
		we obtain \eqref{eq:apriori-closed}.
	\end{proof}
	
	\section{Proof of Theorem~\ref{thm:main}}
	\label{sec:proof-main}
	
	In this final section, we combine the local existence theory with the closed a priori energy inequality of Proposition~\ref{prop:apriori-closure} and homogeneous three-dimensional Gagliardo--Nirenberg interpolation to complete the proof of Theorem~\ref{thm:main}.
	
	\subsection{Local existence and continuation criterion}
	
	We first record the standard local well-posedness and continuation criterion for strong solutions around the smooth planar background $(\bar\rho,\bar u e_1,\bar\theta)$.
	
	\begin{lemma}
		\label{lem:3d-local}
		Let $I_0=[t_0,t_0+T_0]$ be a finite time interval, and let $(\bar\rho,\bar u e_1,\bar\theta)$ be the one-dimensional solution of \eqref{eq:NSF} from Lemma~\ref{lem:background-bounds}, satisfying
		\[
		\sup_{t\in I_0}\|(\bar\rho,\bar u,\bar\theta)(t)\|_{W^{6,\infty}(\R)}\le M_0.
		\]
		Suppose that $U(t_0)=(\phi,\psi,\zeta)(t_0)\in H^4(\R^3)$ satisfies
		\[
		\inf_{x\in\R^3}\rho(t_0,x)>0,\qquad
		\inf_{x\in\R^3}\theta(t_0,x)>0.
		\]
		Then there exists $\tau_*>0$, depending only on the infima of $\rho(t_0,\cdot)$ and $\theta(t_0,\cdot)$, $\|U(t_0)\|_{H^4(\R^3)}$, the physical constants, and $M_0$, such that on $[t_0,t_0+t_*]$ with $t_*=\min\{T_0,\tau_*\}$, the Cauchy problem for \eqref{eq:NSF} admits a unique strong solution $U=(\phi,\psi,\zeta)$ satisfying
		\begin{align*}
			&U\in C([t_0,t_0+t_*];H^4(\R^3)), \qquad
			(\psi,\zeta)\in L^2(t_0,t_0+t_*;H^5(\R^3)), \\[2pt]
			&\nabla\phi\in L^2(t_0,t_0+t_*;H^3(\R^3)), \qquad
			U_t\in L^2(t_0,t_0+t_*;H^3(\R^3)),
		\end{align*}
		with $\rho,\theta>0$ on $[t_0,t_0+t_*]\times\R^3$. Moreover, the local solution can be extended beyond any finite time $T>0$ whenever
		\begin{equation}
			\sup_{0\le t<T}\|U(t)\|_{H^4(\R^3)}<\infty,\qquad
			\inf_{0\le t<T,\,x\in\R^3}\rho(t,x)>0,\qquad
			\inf_{0\le t<T,\,x\in\R^3}\theta(t,x)>0.
			\label{eq:continuation}
		\end{equation}
	\end{lemma}
	
	{\color{black}\begin{proof}
			We apply the local iteration of \cite[Sections~3--5]{MN80} to \eqref{eq:physical-pert-system}, using the background and its derivatives as $L^\infty$ multipliers. At each step we freeze the transport velocity and the diffusion coefficients at the preceding iterate. We first solve the momentum and temperature equations with the remaining terms evaluated at that iterate, and then solve the density equation using the new $\diver\psi$ in the principal coupling.
			
			On a fixed $H^4$ ball about the initial data, with density and temperature bounded below by half their initial infima, the parabolic part is uniformly elliptic and its sources are bounded in $H^3$. Differentiating up to order three and testing with the negative Laplacian of the differentiated solution gives, together with the zeroth-order estimate,
			\[
			\frac{d}{dt}\|(\psi,\zeta)\|_{H^4}^2
			+c\|\nabla(\psi,\zeta)\|_{H^4}^2\le C\big(1+\|(\psi,\zeta)\|_{H^4}^2\big)
			\]
			in the equivalent derivative norm. Here the highest coefficient derivative is paired with a second derivative in $L^\infty$, using $H^4(\R^3)\hookrightarrow W^{2,\infty}(\R^3)$. The fourth-order transport estimate gives
			\[
			\frac{d}{dt}\|\phi\|_{H^4}\le C\big(1+\|\phi\|_{H^4}+\|(\psi,\zeta)\|_{H^5}\big).
			\]
			The last term has time integral at most $\sqrt{t_*}\|(\psi,\zeta)\|_{L^2(t_0,t_0+t_*;H^5)}$. The equations also bound $U_t$ in $L^2(t_0,t_0+t_*;H^3)$, so the change in $U$ is at most $C\sqrt{t_*}$ in $L^\infty$. Since the background has bounded time derivatives, these estimates preserve the positive lower bounds and a uniform $H^4$ bound for sufficiently small $t_*$.
			
			For successive iterates, the difference of diffusion coefficients times the preceding $\mathcal L\psi$ or $\Delta\zeta$ is bounded in $H^2$ by the $H^3$ norm of the preceding difference, since $H^2(\R^3)$ is an algebra. The same bound holds for the remaining source differences. With zero initial difference, the parabolic estimate one order lower bounds the new velocity and temperature difference in $L^\infty H^3\cap L^2H^4$ by $C\sqrt{t_*}$ times the preceding $L^\infty H^3$ difference. The density estimate uses the new velocity divergence in $L^1H^3$, with the additional factor $\sqrt{t_*}$ from Cauchy--Schwarz. For $t_*$ small this gives contraction and uniqueness. The uniform estimates pass to the limit by weak compactness. The parabolic bounds and the transport equation for the fourth derivatives of $\phi$ give continuity into $H^4$. Under \eqref{eq:continuation} the same bounds give a uniform local lifespan as $t_0$ approaches $T$, which proves the continuation criterion.
	\end{proof}}
	
	\subsection{Time-weighted energy bounds and chain-rule estimates for \texorpdfstring{$(q,h)$}{(q,h)}}
	
	On any time interval $[0,T]$ on which \eqref{eq:apriori-assumption} holds, Proposition~\ref{prop:apriori-closure} gives \eqref{eq:apriori-closed}. Using the logarithmic bound \eqref{eq:Lambda-int} from Lemma~\ref{lem:weighted-multiplier} and the Cauchy--Schwarz inequality, we have
	\[
	\int_0^t\sqrt{\frac{\Lambda(s)}{1+s}}\,ds
	\le\left(\int_0^t\Lambda(s)\,ds\right)^{1/2}
	\bigl[\log(1+t)\bigr]^{1/2}
	\le C\sqrt{\delta}\log(1+t),
	\]
	and hence, by Lemma~\ref{lem:background-bounds} and Lemma~\ref{lem:weighted-multiplier},
	\[
	C\int_0^t\Gamma(s)\,ds
	\le C(\sigma+\sqrt{\delta}+\delta)\log(1+t)
	\le 2\eta\log(1+t),
	\]
	where $\eta=C\sqrt{\delta}<1/8$. Multiplying \eqref{eq:apriori-closed} by $\exp\bigl(-C\int_0^t\Gamma(s)\,ds\bigr)$ and integrating over $[0,t]$ yields
	\begin{equation}
		(1+t)^{-2\eta}\mathcal{E}(t)
		+c\int_0^t(1+s)^{-2\eta}\mathcal{D}(s)\,ds
		\le C\eps^2\qquad \text{for all }t\in[0,T].
		\label{eq:weighted-energy-bound}
	\end{equation}
	
	Since $\mathcal{E}(t)\sim\mathcal{N}(t)$ and $\omega(t,x_1)\ge(1+t)^{-1}$, \eqref{eq:weighted-energy-bound} implies the $L^2(\R^3)$ bounds
	\begin{align}
		&\|U(t)\|_2+\|h(t)\|_2\le C\eps(1+t)^{1/2+\eta},\notag \\[2pt]
		&\|(q,\psi)(t)\|_2+\|\sqrt{\omega}\,h(t)\|_2+\|\nabla U(t)\|_2
		\le C\eps(1+t)^\eta,\label{eq:L2-low} \\[2pt]
		&\|\nabla(q,\psi)(t)\|_2+\sum_{j=2}^4\|\nabla^j U(t)\|_2
		\le C\eps(1+t)^{-1/2+\eta}.\label{eq:L2-high}
	\end{align}
	Substituting \eqref{eq:L2-low}--\eqref{eq:L2-high} into the chain-rule estimate \eqref{eq:qh-chain} gives
	\begin{equation}
		\|\nabla^\ell(q,h)(t)\|_2\le C\eps(1+t)^{-1/2+\eta},\qquad 2\le\ell\le 4.
		\label{eq:qh-high}
	\end{equation}
	
	\subsection{\texorpdfstring{$L^p$}{Lp} decay estimates and global continuation}
	
	To derive $L^\infty$ decay from the $L^2(\R^3)$ bounds \eqref{eq:L2-low}--\eqref{eq:qh-high} without using the growing zeroth-order norm of $(U,h)$, we use the homogeneous three-dimensional Gagliardo--Nirenberg inequalities
	\begin{equation}
		\|f\|_\infty\le C\|\nabla f\|_2^{1/2}\|\nabla^2f\|_2^{1/2}, \qquad
		\|f\|_3\le C\|f\|_2^{1/2}\|\nabla f\|_2^{1/2}.
		\label{eq:3d-GN}
	\end{equation}
	Applying the first inequality in \eqref{eq:3d-GN} to $f=U,h$ and to $f=q,\psi,\nabla U,\nabla^2 U$ via \eqref{eq:L2-low}--\eqref{eq:qh-high} gives
	\begin{align*}
		&\|U(t)\|_\infty+\|h(t)\|_\infty\le C\eps(1+t)^{-1/4+\eta},\\[2pt]
		&\|(q,\psi)(t)\|_\infty+\|(p-\bar p)(t)\|_\infty+\|\nabla U(t)\|_\infty+\|\nabla^2 U(t)\|_\infty
		\le C\eps(1+t)^{-1/2+\eta},
	\end{align*}
	which proves \eqref{eq:main-decay}.
	
	Similarly, applying the second inequality in \eqref{eq:3d-GN} to $\psi$, $\nabla U$, and $\nabla^2 U$ gives
	\begin{align*}
		&\|\psi(t)\|_3+\|\nabla U(t)\|_3\le C\eps(1+t)^{-1/4+\eta}, \\[2pt]
		&\|\nabla^2U(t)\|_3\le C\eps(1+t)^{-1/2+\eta}, \\[2pt]
		&\sqrt{1+t}\,\|\psi(t)\|_\infty\|\nabla U(t)\|_3
		\le C\eps^2(1+t)^{-1/4+2\eta},
	\end{align*}
	and therefore $\mathcal{M}(t)$ defined in \eqref{eq:M-def} satisfies
	\begin{equation}
		\mathcal{M}(t)\le C\eps(1+t)^{-1/4+\eta}+C\eps^2(1+t)^{-1/4+2\eta}\qquad \text{for all }t\in[0,T].
		\label{eq:M-decay}
	\end{equation}
	{\color{black}By the homogeneous Sobolev inequality and \eqref{eq:L2-high},
		\[
		\|(q,\psi,\nabla U)(t)\|_6
		\le C\big(\|\nabla(q,\psi)(t)\|_2+\|\nabla^2U(t)\|_2\big)
		\le C\eps(1+t)^{-1/2+\eta}.
		\]
		Interpolating for $(q,\psi,\nabla U)$ between $L^2,L^6$ and then $L^6,L^\infty$, and for $(U,h)$ and $\nabla^2 U$ between $L^2,L^\infty$, gives}
	\begin{align*}
		&\|U(t)\|_p+\|h(t)\|_p \le C_p\eps(1+t)^{-1/4+3/(2p)+\eta},\\[2pt]
		&\|(q,\psi)(t)\|_p+\|\nabla U(t)\|_p \le C_p\eps
		\begin{cases}
			(1+t)^{-3/4+3/(2p)+\eta}, & 2\le p\le 6,\\[2pt]
			(1+t)^{-1/2+\eta}, & 6\le p\le\infty,
		\end{cases}\\[2pt]
		&\|\nabla^2 U(t)\|_p \le C_p\eps(1+t)^{-1/2+\eta}.
	\end{align*}
	
	\begin{proof}[Completion of the proof of Theorem~\ref{thm:main}]
		Choose $\delta_0>0$ sufficiently small so that all background and energy estimates hold for $0<\delta\le\delta_0$ with $\eta={\color{black}C\sqrt\delta\le }C\sqrt{\delta_0}<1/8$. The constants in all preceding estimates are then uniform with respect to $\delta\in(0,\delta_0]$. By the Sobolev embedding $H^4(\R^3)\hookrightarrow W^{2,\infty}(\R^3)$ and \eqref{eq:3d-GN}, the initial data satisfy
		\[
		\|U_0\|_{W^{2,\infty}(\R^3)}\le C\eps,\qquad
		\mathcal{M}(0)\le C\eps+C\eps^2.
		\]
		We may therefore choose $\eps_0>0$ small enough, independently of $\delta$, that the initial data satisfy \eqref{eq:apriori-assumption} with strict inequalities at $t=0$.
		
		Let $T_{\max}>0$ be the maximal existence time of the strong solution given by Lemma~\ref{lem:3d-local}, and let $T_*\in(0,T_{\max}]$ be the maximal time on which \eqref{eq:apriori-assumption} holds. On any interval $[0,T]$ with $T<T_*$, the estimates \eqref{eq:main-decay} and \eqref{eq:M-decay} hold with constants independent of $T$. Since $\eta<1/8$, reducing $\eps_0>0$ if necessary ensures that
		\[
		\|\phi(t)\|_\infty\le\frac{\rho_*}{4},\qquad
		\|\zeta(t)\|_\infty\le\frac{\theta_*}{4},\qquad
		\mathcal M(t)\le\frac{\delta_1}{2}\qquad \text{for all }t\in[0,T],
		\]
		which strictly improves the bounds in \eqref{eq:apriori-assumption}. By continuity in $H^4(\R^3)$, we obtain $T_*=T_{\max}$.
		
		Finally, if $T_{\max}<\infty$, then \eqref{eq:L2-low}--\eqref{eq:L2-high} and the improved $L^\infty$ bounds imply
		\[
		\sup_{0\le t<T_{\max}}\|U(t)\|_{H^4(\R^3)}
		\le C\eps(1+T_{\max})^{1/2+\eta}<\infty,\qquad
		\rho(t,x)\ge\frac{3\rho_*}{4},\qquad
		\theta(t,x)\ge\frac{3\theta_*}{4},
		\]
		on $[0,T_{\max})\times\R^3$, so the continuation criterion \eqref{eq:continuation} in Lemma~\ref{lem:3d-local} extends the solution beyond $T_{\max}$, a contradiction. Thus $T_{\max}=\infty$, and \eqref{eq:main-decay} and \eqref{eq:main-energy} hold for all $t\ge 0$. This completes the proof of Theorem~\ref{thm:main}.
	\end{proof}
	
	\medskip
	\noindent {\bf Acknowledgment:}\,
	The research of Feimin Huang was supported in part by the National Key
	R\&D Program of China, grant no. 2021YFA1000800, and the National Natural Sciences Foundation
	of China, grant no. 12288201.
	
	\medskip
	\noindent {\bf Disclosure of AI assistance:}
	The key ideas of this paper are the formulation in the velocity, pressure, and entropy variables $(u,p,s)$, the separate treatment of the estimates for $(u,p)$ and $s$, and the construction of the weight functions used in the proof. These ideas were developed by the authors through their discussions from 2025 to the present. During the preparation of the manuscript, the authors used ChatGPT to assist with language polishing, checking calculations and proofs, and carrying out auxiliary computations. The authors take full responsibility for the accuracy, originality, and integrity of the entire paper.
	
	\medskip
	\noindent{\bf Data availability:} The manuscript contains no associated data.

	\medskip
	\noindent{\bf Conflict of Interest:} The authors declare that they have no conflict of interest.

\end{document}